\def\wlog#1{}
 \makeatletter
 \def\@latex@info#1{}
 \makeatother

 \makeatletter
 \def\@font@info#1{}
 \makeatother

\documentclass[12pt]{amsart}

\usepackage{amssymb}

\usepackage{graphicx}

\usepackage{upref}

\newcommand{\conservepaper}{
\setlength{\textwidth}{6.5 in}
\setlength{\textheight}{9.0 in}
\hoffset=-0.75in
\voffset=-0.5in
 }
 \conservepaper

\swapnumbers

\begingroup \makeatletter 

\gdef\swappedhead#1#2#3{%
  \thmnumber{\@upn{\the\thm@headfont#2}}
  \thmname{\@ifnotempty{#2}{.~}#1}%
  \thmnote{ {\the\thm@notefont#3}}}

\gdef\swappedhead@plain#1#2#3{%
  \thmnumber{\@upn{#2}}\thmname{\@ifnotempty{#2}{. }#1}%
  \thmnote{ \textmd{\upshape#3}}}

\endgroup 

\newenvironment{thmref}{\thmrefer}{}
\newcommand{\thmrefer}[1]{\renewcommand\theequation
  {\protect\ref{#1}$'$}\addtocounter{equation}{-1}}

\newenvironment{thmrefref}{\thmreferref}{}
\newcommand{\thmreferref}[1]{\renewcommand\theequation
  {\protect\ref{#1}$''$}\addtocounter{equation}{-1}}

\numberwithin{figure}{section}

\numberwithin{equation}{section}
\newtheorem{thm}[equation]{Theorem}

\newtheorem{lem}[equation]{Lemma}
\newtheorem{cor}[equation]{Corollary}
\newtheorem{prop}[equation]{Proposition}

\theoremstyle{definition}

\newtheorem{defn}[equation]{Definition}
\newtheorem{conj}[equation]{Conjecture}

\newtheorem{eg}[equation]{Example}

\newtheorem{ques}[equation]{Question}
\newtheorem{standing}[equation]{Standing assumptions}
\newtheorem{notation}[equation]{Notation}

\theoremstyle{remark}

\newtheorem{rem}[equation]{Remark}

\newtheorem{ack}[equation]{Acknowledgment} 

\newcommand{\pref}[1]{{\upshape(}\ref{#1}{\upshape)}}
\renewcommand{\see}[1]{{\upshape(}see~\ref{#1}{\upshape)}}
\newcommand{\seesee}[2]{{\upshape(}see~\ref{#1} {#2}{\upshape)}}
\newcommand{\seeand}[2]{\seesee{#1}{and~\ref{#2}}}
\newcommand{\fullsee}[2]{{\upshape(}see~\fullref{#1}{#2}{\upshape)}}
\newcommand{\cf}[1]{(cf.~\ref{#1})}
\newcommand{\fullcf}[2]{{\upshape(}cf.~\fullref{#1}{#2}{\upshape)}}
\newcommand{\fullref}[2]{\ref{#1}\pref{#1-#2}}
\newcommand{\bibURL}[1]{{\unskip\nobreak\hfil\penalty50{\tt#1}}}

\makeatletter
 \def\l@subsection{\@tocline{2}{0pt}{4pc}{5pc}{}}
 \makeatother

\newcommand{\Lie}[1]{\mathfrak{\lowercase{#1}}}
\newcommand{\rsp}{\Lie N}
\newcommand{\rsg}{N}
\newcommand{\Dn}{\Lie D}
\newcommand{\Dh}{\Dn_{\Lie H}}
\newcommand{\so}{\operatorname{\Lie{SO}}}
\newcommand{\su}{\operatorname{\Lie{SU}}}
\newcommand{\Liesl}{\operatorname{\Lie{SL}}}
\newcommand{\Liesp}{\operatorname{\Lie{Sp}}}

\newcommand{\SO}{\operatorname{SO}}
\newcommand{\SL}{\operatorname{SL}}
\newcommand{\GL}{\operatorname{GL}}
\newcommand{\SU}{\operatorname{SU}}
\newcommand{\Symp}{\operatorname{Sp}}
\newcommand{\Id}{\operatorname{Id}}
\newcommand{\real}{\mathord{\mathbb{R}}}
\newcommand{\complex}{\mathord{\mathbb{C}}}
\newcommand{\integer}{\mathord{\mathbb{Z}}}
\newcommand{\quaternion}{\mathord{\mathbb{H}}}
\newcommand{\homology}{\mathord{\mathcal{H}}}
\newcommand{\hilbert}{\mathord{\mathcal{H}}}

\newcommand{\Rrank}{\operatorname{\hbox{\upshape$\real$-rank}}}
\newcommand{\projX}{\mathord{\mathbb{P}X}}
\newcommand{\F}{\mathbb{F}}
\newcommand{\df}{q} 
\newcommand{\Fim}{\F_{\text{imag}}}
\newcommand{\dimF}{\operatorname{\dim_{\F}}}
\newcommand{\dimR}{\operatorname{\dim_{\real}}}
\newcommand{\dimC}{\operatorname{\dim_{\complex}}}
\newcommand{\Mat}{\operatorname{Mat}}
\newcommand{\cjg}[1]{\overline{#1}}
\newcommand{\Hc}{\Hcc{c}}
\newcommand{\Hcc}[1]{H_{[#1]}}
\newcommand{\LieHc}{\LieHcc{c}}
\newcommand{\LieHcc}[1]{{\Lie H}_{[#1]}}
\newcommand{\iso}{\cong}
\newcommand{\diffeo}{\simeq}
\newcommand{\Ad}{\operatorname{Ad}\nolimits}
\newcommand{\ad}{\operatorname{ad}\nolimits}

\newcommand{\bigset}[2]{\left\{\, #1
 \mathrel{\left| \vphantom {\left\{ #1 \mid #2 \right\} } \right.}
 #2 \,\right\} }

 \newcounter{case}
 \newenvironment{case}[1][\unskip]{\refstepcounter{case}\em
 \medskip \noindent Case \thecase\ #1.\ }{\unskip\upshape}
 \renewcommand{\thecase}{\arabic{case}}

 \newcounter{subcase}
 \newenvironment{subcase}[1][\unskip]{\refstepcounter{subcase}\em
 \medskip \noindent Subcase \thesubcase\ #1.\ }{\unskip\upshape}
\numberwithin{subcase}{case}

 \newcounter{subsubcase}
 
\numberwithin{subsubcase}{subcase}

 \newcounter{step}
\newenvironment{step}[1][\unskip]{\refstepcounter{step}
\em
 \medskip \noindent Step \thestep\
 #1.\ }{\unskip\upshape}

 \newcounter{stepinsubcase}
\newenvironment{stepinsubcase}[1][\unskip]{\refstepcounter{stepinsubcase}
\em
 \medskip \noindent Step \thestepinsubcase\
 #1.\ }{\unskip\upshape}
 \numberwithin{stepinsubcase}{subcase}

 \newcounter{substep}
 \newenvironment{substep}[1][\unskip]{\refstepcounter{substep}\em
 \medskip \noindent Substep \thesubstep\ #1.\ }{\unskip\upshape}
\numberwithin{substep}{step}

\newenvironment{pfassump}[1][\unskip]{\em
 \medskip \noindent Assumption.\ }{\unskip\upshape}

\renewcommand{\Re}{\operatorname{Re}}
\renewcommand{\Im}{\operatorname{Im}}

\newcommand{\s}{\vphantom{\vrule height 15pt depth 5pt}}

\newcommand{\xx}{{\mathord{\mathsf{x}}}}
\newcommand{\yy}{{\mathord{\mathsf{y}}}}

\begin{document}

\title[Tessellations of homogeneous spaces]{Tessellations of homogeneous
spaces \\
 of classical groups of real rank two}

\author{Alessandra Iozzi}

 \address
 {Department of Mathematics,
 University of Maryland,
 College Park, MD 20910
 USA}

\curraddr
 {FIM,
 ETH Zentrum,
 CH--8092 Z\"urich,
 Switzerland}

 \email{iozzi@math.ethz.ch}

\author{Dave Witte}
 \address
 {Department of Mathematics,
 Oklahoma State University,
 Stillwater, OK 74078
 USA}

 \email{dwitte@math.okstate.edu}

\date{\today}

\begin{abstract}
  Let $H$ be a closed, connected subgroup of a connected, simple Lie group $G$
with finite center.  The homogeneous space $G/H$ has a {\it tessellation} if
there is a discrete subgroup $\Gamma$ of $G$, such that $\Gamma$ acts properly
discontinuously on $G/H$, and the double-coset space $\Gamma\backslash G/H$ is
compact.  Note that if either $H$ or $G/H$ is compact, then $G/H$ has a
tessellation; these are the obvious examples.

It is not difficult to see that if $G$ has real rank one, then only the
obvious homogeneous spaces have tessellations.  Thus, the first interesting
case is when $G$ has real rank two.  In particular, R.~Kulkarni and
T.~Kobayashi constructed examples that are not obvious when $G =
\SO(2,2n)^\circ$ or $\SU(2,2n)$. H.~Oh and D.~Witte constructed additional
examples in both of these cases, and obtained a complete classification when
$G = \SO(2,2n)^\circ$. We simplify the work of Oh-Witte, and extend it to
obtain a complete classification when $G = \SU(2,2n)$. This includes the
construction of another family of examples.

The main results are obtained from methods of Y.~Benoist and T.~Kobayashi: we
fix a Cartan decomposition $G = K A^+ K$, and study the intersection $(KHK)
\cap A^+$. Our exposition generally assumes only the standard theory of
connected Lie groups, although basic properties of real algebraic groups are
sometimes also employed; the specialized techniques that we use are developed
from a fairly elementary level.
 \end{abstract}

\maketitle

\tableofcontents

\section{Introduction}

\begin{defn}[{\cite[pp.~43--44]{KobayashiNomizu1}}] \label{PDDefn}
 A group~$\Gamma$ of homeomorphisms of a topological space~$M$ acts
\emph{properly discontinuously} on~$M$ if, for every compact subset~$C$
of~$M$,
 $$ \mbox{ $\{\, \gamma \in \Gamma \mid C \cap \gamma C \neq \emptyset \,\}$
 is finite.} $$
 \end{defn}

Classically, a discrete group~$\Gamma$ of isometries of a Riemannian
manifold~$M$ is a crystallographic group if $\Gamma$ acts properly
discontinuously on~$M$, and the quotient $\Gamma \backslash M$ is compact. The
$\Gamma$-translates of any fundamental domain for $\Gamma \backslash M$ form a
tessellation of~$M$. 

These notions generalize to any homogeneous space, even without an invariant
metric.

\begin{defn} \label{TessDefn}
 Let
 \begin{itemize}
 \item $G$ be a Lie group and
 \item $H$ be a closed subgroup of~$G$.
 \end{itemize}
 A discrete subgroup~$\Gamma$ of~$G$ is a \emph{crystallographic group} for
$G/H$ if
 \begin{enumerate}
 \item $\Gamma$ acts properly discontinuously on~$G/H$; and
 \item $\Gamma \backslash G/H$ is compact.
 \end{enumerate}
 We say that $G/H$ has a \emph{tessellation} if there exists a
crystallographic group~$\Gamma$ for $G/H$.
 \end{defn}

Crystallographic groups and the corresponding tessellations have been studied
for many groups~$G$. (A brief recent introduction to the subject is given in
\cite{Kobayashi-survey00}.) The classical Bieberbach Theorems
\cite[Chap.~1]{CharlapBook} deal with the case where $G$ is the group of
isometries of Euclidean space $\real^n = G/H$. As another example, the
Auslander Conjecture \cite{Abels, AbelsMargulisSoifer, FriedGoldman,
Margulis-Auslander, Tomanov} asserts that if $G$ is the group of all affine
transformations of~$\real^n$, then every crystallographic group has a
solvable subgroup of finite index. In addition, the case where $G$ is solvable
has been discussed in \cite{Witte-Solvtess}. 

In this paper, we focus on the case where $G$ is a simple Lie group, such as
$\SL(n,\real)$, $\SO(m,n)$, or $\SU(m,n)$.

\begin{standing} \label{standing}
 Throughout this paper:
 \begin{enumerate}
 \item \label{standing-G}
 $G$ is a  linear, semisimple Lie group with only finitely many connected
components; and
 \item \label{standing-H}
 $H$ is a closed subgroup of~$G$ with only finitely many connected components.
 \end{enumerate}
 \end{standing}

\begin{rem} \label{Hdisconnected}
 Because $H/H^\circ$ is finite (hence compact), it is easy to see that $G/H$
has a tessellation if and only if $G/H^\circ$ has a tessellation. Also, if
$G/H$ has a tessellation, then $G^\circ / H^\circ$ has a tessellation.
Furthermore, the converse holds in many situations. (See~\S\ref{disconnected}
for a discussion of this issue.) Thus, there is usually no harm in assuming
that both $G$ and~$H$ are connected; we will feel free to do so whenever it
is convenient. On the other hand, because $\SO(m,n)$ is usually not connected
(it usually has two components \cite[Lem.~10.2.4, p.~451]{HelgasonBook}), it
would be somewhat awkward to make this a blanket assumption.
 \end{rem}

\begin{eg} \label{classical}
 There are two classical cases in which $G/H$ is well known to have a
tessellation.
 \begin{enumerate}
 \item If $G/H$ is compact, then we may let $\Gamma = e$ (or any finite
subgroup of~$G$).
 \item \label{classical-Borel}
 If $H$ is compact, then we may let $\Gamma$ be any cocompact lattice in~$G$.
(A.~Borel \cite{Borel-CK} proved that every connected, simple Lie group has a
cocompact lattice.)
 \end{enumerate}
  Thus, the existence of a tessellation is an interesting question only when
neither $H$ nor~$G/H$ is compact. (In this case, any crystallographic
group~$\Gamma$ must be infinite, and cannot be a lattice in~$G$.)
 \end{eg}

Given~$G$ (satisfying~\fullref{standing}{G}), we would like to find all the
subgroups~$H$ (satisfying~\fullref{standing}{H}), such that $G/H$ has a
tessellation. This seems to be a difficult problem in general. (See the
surveys \cite{Kobayashi-survey97} and \cite{Labourie-survey} for a discussion
of the many partial results that have been obtained, mainly under the
additional assumption that $H$ is reductive.) However, it can be solved in
certain cases of low real rank. In particular, as we will now briefly
explain, the problem is very easy if $\Rrank G = 0$ or~$1$. Most of this
paper is devoted to solving the problem for certain cases where $\Rrank G =
2$.

If $\Rrank G = 0$ (that is, if $G$ is compact), then $G/H$ must be compact
(and $H$ must also be compact), so $G/H$ has a tessellation, but this is not
interesting. If $\Rrank G = 1$, then there are some interesting homogeneous
spaces, but it turns out that none of them have tessellations.

\begin{center}
 \begin{figure}
 \includegraphics[scale=0.5]{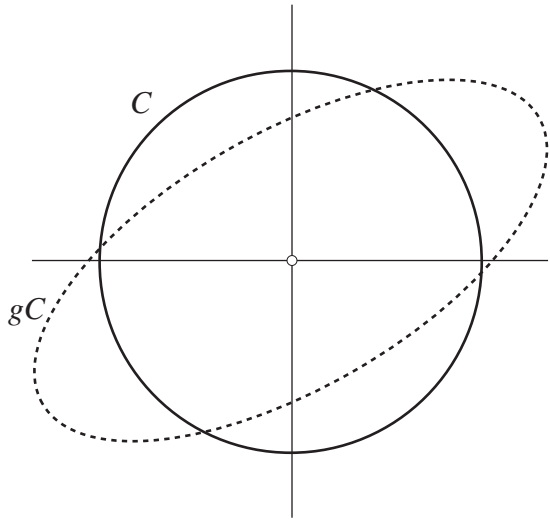}
 \caption{The Calabi-Markus Phenomenon (Example~\ref{CalabiMarkusCircle}): $C
\cap gC \neq \emptyset$, for every $g \in \SL(2,\real)$, so no infinite
subgroup of $\SL(2,\real)$ acts properly discontinuously.}
 \label{circle+ellipse}
 \end{figure}
 \end{center}

\begin{eg} \label{CalabiMarkusCircle}
 $G = \SL(2,\real)$ is transitive on $\real^2 - \{0\}$, so $\real^2 - \{0\}$
is a homogeneous space for~$G$. It does not have a tessellation, for reasons
that we now explain.

Let $C$ be the unit circle, so $C$ is a compact subset of $\real^2 -
\{0\}$.

We claim that $C \cap gC \neq \emptyset$, for every $g \in G$ (cf.\ 
Figure~\ref{circle+ellipse}). To see this, note that, because $\det g = 1$,
the ellipse bounded by $gC$ has the same area as the disk bounded by~$C$, so
$gC$ cannot be contained in the interior of the disk bounded by~$C$, and
cannot contain $C$ in its interior. Thus, $gC$ must be partly inside~$C$ and
partly outside, so $gC$ must cross~$C$, as claimed.

Let $\Gamma$ be any discrete subgroup of~$G$. The preceding paragraph
implies that $C \cap \gamma C \neq \emptyset$, for every $\gamma \in \Gamma$.
If $\Gamma$ acts properly discontinuously on $\real^2 - \{0\}$, then,
because $C$ is compact, this implies that $\Gamma$ is finite. So the
quotient $\Gamma \backslash (\real^2 - \{0\})$ is not compact. Therefore
$\Gamma$ is not a crystallographic group. We have shown that no subgroup
of~$G$ is a crystallographic group, so we conclude that $\real^2 - \{0\}$
does not have a tessellation.
 \end{eg}

This example illustrates the \emph{Calabi-Markus Phenomenon}: if there is a
compact subset~$C$ of $G/H$, such that $C \cap gC \neq \emptyset$, for every
$g \in G$, then no infinite subgroup of~$G$ acts properly discontinuously on
$G/H$ \see{CalabiMarkus}. Thus, $G/H$ does not have a tessellation, unless
$G/H$ is compact \see{CDS->notess}.

We will see in Section~\ref{CartanSect} that the following proposition can be
proved quite easily from basic properties of the Cartan projection.

\begin{thmref}{Rrank1-CDS}
 \begin{prop}[{(cf.~\cite[Lem.~3.2]{Kobayashi-isotropy})}] \label{Rrank1-CM}
 If $\Rrank G = 1$, and $H$ is not compact, then there is a compact subset~$C$
of $G/H$, such that $C \cap gC \neq \emptyset$, for every $g \in G$.
 \end{prop}
 \end{thmref}

\begin{cor}[{(Kulkarni)}] \label{rank1->notess}
 If $\Rrank G = 1$, and neither $H$ nor $G/H$ is compact, then $G/H$ does not
have a tessellation.
 \end{cor}

We now consider groups of real rank two. The obvious
example is $\SL(3,\real)$, but, in this case, once again, none of the
interesting homogeneous spaces have tessellations. Moreover, the same is true
when real numbers are replaced by complex numbers or quaternions. The case
where $\dim H = 1$ relies on beautiful methods of Y.~Benoist and F.~Labourie
\cite{BenoistLabourie} or G.~A.~Margulis \cite{Margulis-CK}, which we
describe in Section~\ref{1DSect}.

\begin{thm}[{(Benoist, Margulis, Oh-Witte, see~\S\ref{SL3notessSect})}]
\label{SL3->notess}
 If
 $$ \hbox{$G = \SL(3,\real)$, $\SL(3,\complex)$, or $\SL(3,\quaternion)$,} $$
 and neither $H$ nor $G/H$ is compact, then $G/H$ does not have a 
tessellation.
 \end{thm}

It is important to note that some interesting homogeneous spaces do have
tessellations.

\begin{eg} \label{G=LxH}
 Suppose $G = L \times H$, and let $\Gamma$ be a cocompact lattice in~$L$.
Then $\Gamma$ acts properly discontinuously on $L \cong G/H$, and $\Gamma
\backslash G/H \cong \Gamma \backslash L$ is compact. So $G/H$ has a
tessellation.
 \end{eg}

The following easy lemma generalizes this example to the situation where $G$
is a more general product of $L$ and~$H$, not necessarily a direct product.

\begin{thmref}{construct-tess}
 \begin{lem} \label{G=LH}
 Let $H$ and~$L$ be closed subgroups of~$G$, such that
 \begin{itemize}
 \item $G = LH$,
 \item $L \cap H$ is compact; and
 \item $L$ has a cocompact lattice~$\Gamma$.
 \end{itemize}
 Then $G/H$ has a tessellation.  {\upshape(}Namely, $\Gamma$ is a
crystallographic group for $G/H$.{\upshape)}
 \end{lem}
 \end{thmref}

For $G = \SO(2,n)$ or $\SU(2,n)$, this lemma leads to some interesting
examples found by R.~Kulkarni \cite[Thm.~6.1]{Kulkarni} and T.~Kobayashi
\cite[Prop.~4.9]{Kobayashi-properaction}. 

 \begin{eg}[(Kulkarni, Kobayashi)]
\label{KulkarniEg}
 There are natural embeddings
 $$ \hbox{$\SO(1,n) \hookrightarrow \SO(2,n)$
 \qquad and \qquad
 $\SU(1,n) \hookrightarrow \SU(2,n)$.} $$
 Furthermore, identifying $\complex^{1+m}$ with $\real^{2+2m}$ yields
an embedding
 $$\SU(1,m) \hookrightarrow \SO(2,2m) .$$
 Similarly, identifying $\quaternion^{1+m}$ with $\complex^{2+2m}$
yields an embedding
 $$ \Symp(1,m) \hookrightarrow \SU(2,2m) .$$
 Thus, we may think of $\SO(1,2m)$ and $\SU(1,m)$ as subgroups of
$\SO(2,2m)$; and we may think of $\SU(1,2m)$ and $\Symp(1,m)$ as subgroups
of $\SU(2,2m)$.

With the above understanding, we see that $\SO(1,2m)$ is the stabilizer
of a vector of norm~$+1$. Since $\SU(1,m)$ is transitive on the set of such
vectors, we have
 $$ \SO(2,2m) = \SO(1,2m) \, \SU(1,m) .$$
 Similarly,
 $$\SU(2,2m) = \SU(1,2m) \, \Symp(1,m) . $$
 Then Lemma~\ref{G=LH} implies that each of the following four
homogeneous spaces has a tessellation:
 \begin{itemize}
 \item $\SO(2,2m)/\SO(1,2m)$,
 \item $\SO(2,2m)/\SU(1,m)$,
 \item $\SU(2,2m)/\SU(1,2m)$, and
 \item $\SU(2,2m)/\Symp(1,m)$.
 \end{itemize}
 \end{eg}

\begin{rem}
 When discussing $\SO(2,n)$ or $\SU(2,n)$, we always assume $n > 2$. This
causes no harm, because $\SO(2,2)$ is locally isomorphic to $\SL(2,\real)
\times \SL(2,\real)$ \cite[(x), p.~520]{HelgasonBook}, and $\SU(2,2)$ is
locally isomorphic to $\SO(2,4)$ \cite[(vi), p.~519]{HelgasonBook}. 
 \end{rem}

When $n$ is even, H.~Oh and D.~Witte \cite{OhWitte-CK} provided a complete
description of all the (closed, connected) subgroups~$H$, such that
$\SO(2,n)/H$ has a tessellation, but their classification is not quite
complete when $n$~is odd.

In this paper, we extend the work of Oh and Witte to obtain analogous results
for homogeneous spaces of $G = \SU(2,n)$. We also give a much shorter proof of
the main results of \cite{OhWitte-CK}. The same techniques should yield
significant results for homogeneous spaces of the other simple groups of real
rank two, although the calculations seem to be difficult. On
the other hand, the groups of higher real rank require different ideas.

Once one knows that a
tessellation of $G/H$ exists, it would be interesting to find \emph{all} of
the crystallographic groups for $G/H$ and, for each crystallographic group,
describe the possible tessellations. These are much more delicate questions,
which we do not address at all. (W.~Goldman \cite{Goldman-nonstandard},
F.~Salein \cite{Salein}, T.~Kobayashi \cite{Kobayashi-deformation}, and
A.~Zeghib \cite{Zeghib-deSitter} have interesting results in some special
cases.)

In the remainder of this introduction, we state the specific results for
homogeneous spaces of $\SO(2,n)$ and $\SU(2,n)$.

\begin{notation}[{(Iwasawa \cite[p.~533]{Iwasawa})}] \label{d(H)-defn}
 For any connected Lie group~$H$, let
 $$ d(H) = \dim H - \dim K_H ,$$
 where $K_H$ is any maximal compact subgroup of~$H$.  This is well defined,
because all the maximal compact subgroups of~$H$ are conjugate
\cite[Thm.~15.3.1(iii), pp.~180--181]{Hochschild-Lie}.
 \end{notation}

\begin{eg} \label{d(G)} If $H$ is semisimple, we have the Iwasawa
decomposition $H = K_H A_H N_H$ \cite[Thm.~6.5.1, pp.~270--271]{HelgasonBook},
from which it is obvious that $d(H) = \dim(A_H N_H)$.

 This yields the following calculations \seeand{d(SU2)}{d(Sp)}):
 \begin{itemize}
 \item $d \bigl( \SO(1,n) \bigr) = n$.
 \item $d \bigl( \SO(2,n) \bigr) = 2n$.
 \item $d \bigl( \SU(1,n) \bigr) = 2n$.
 \item $d \bigl( \SU(2,n) \bigr) = 4n$.
 \item $d \bigl( \Symp(1,n) \bigr) = 4n$.
 \end{itemize}
 \end{eg}

\begin{rem} \label{d(H)=dimH}
 If $H \subset AN$ (for some Iwasawa decomposition $G = KAN$ of~$G$), then
$d(H) = \dim H$ (see~\ref{ANsc} and~\fullref{solvable}{nocpct}).
 \end{rem}

\begin{notation}[{\cite[Defn.~2.1.1]{Kobayashi-criterion}}] \label{simdefn}
 For subgroups $H_1$ and~$H_2$ of~$G$, we write $H_1 \sim H_2$ if there is a
compact subset~$C$ of~$G$, such that $H_1 \subset C H_2 C$ and $H_2 \subset
C H_1 C$.
 \end{notation}

\begin{rem}
 Note that $d$ is not invariant under the equivalence relation~$\sim$. For
example, the Cartan decomposition $G = KAK$ implies that $G \sim A$, but we
have $d(A) = \dim A \neq \dim(AN) = d(G)$.
 \end{rem}

The following two theorems state a version of the main results for even~$n$.

\begin{thmref}{SUF-known}
 \begin{thm}[{(Oh-Witte \cite{OhWitte-CK})}] \label{OW-known}
 Assume $G = \SO(2,2m)$, and let $H$ be a closed, connected, subgroup of
$G$, such that neither $H$ nor $G/H$ is compact.

 The homogeneous space $G/H$ has a tessellation if and only if
 \begin{enumerate}
 \item $d(H) = 2m$; and
 \item either $H \sim \SO(1,2m)$ or $H \sim \SU(1,m)$.
 \end{enumerate}
 \end{thm}
 \end{thmref}

\begin{thmrefref}{SUF-known}
\begin{thm} \label{IW-known}
 Assume $G = \SU(2,2m)$, and let $H$ be a closed, connected, subgroup of $G$,
such that neither $H$ nor $G/H$ is compact.

 The homogeneous space $G/H$ has a tessellation if and only if
 \begin{enumerate}
 \item $d(H) = 4m$; and
 \item either $H \sim \SU(1,2m)$ or $H \sim \Symp(1,m)$.
 \end{enumerate}
 \end{thm}
 \end{thmrefref}

The subgroups~$H$ that arise in Theorems~\ref{OW-known} and~\ref{IW-known} can
also be described more explicitly (cf.~\ref{SOevenTess} and~\ref{SUevenTess}
below).

T.~Kobayashi \cite[1.4]{Kobayashi-deformation} conjectured that if $H$ is
reductive and it is impossible to construct a tessellation of $G/H$ by using
a generalization of Lemma~\ref{G=LH} \see{construct-tess}, then $G/H$ does
not have a tessellation. The following lists three special cases of this
general conjecture.

\begin{conj} \label{notessSU/Sp}
 The homogeneous spaces
 \begin{enumerate}
  \renewcommand{\theenumi}{\alph{enumi}} 

 \item \label{notessSU/Sp-SO2/SU1}
 $\SO(2,2m+1)/\SU(1,m)$,

 \item \label{notessSU/Sp-SU2/Sp1}
  $\SU(2,2m+1)/\Symp(1,m)$, and

 \item \label{notessSU/Sp-SU2/SU1}
  $\SU(2,2m+1)/\SU(1,2m+1)$

 \end{enumerate}
 do not have tessellations.
 \end{conj}

If this conjecture is true, then, for odd~$n$, there is no interesting
example of a homogeneous space of $\SO(2,n)$ or $\SU(2,n)$ that has a
tessellation.

\begin{thmref}{SUF->complete}
 \begin{thm}[{(Oh-Witte \cite[Thm.~1.7]{OhWitte-CK}, Iozzi-Witte)}]
\label{SU1m->complete}
 Assume
 $$ \hbox{$G = \SO(2,2m+1)$ or $\SU(2,2m+1)$,} $$
 and let $H$ be any closed, connected subgroup of $G$, such that neither $H$
nor~$G/H$ is compact.

 If Conjecture~\ref{notessSU/Sp} is true, then $G/H$ does not have a
tessellation.
 \end{thm}
 \end{thmref}

The proof of Theorem~\ref{SU1m->complete} assumes the following special case
proved by R.~Kulkarni \cite[Cor.~2.10]{Kulkarni}. In short, Kulkarni noted
that the Euler characteristic of $\Gamma \backslash G/H$ must both vanish
(because the Euler characteristic of $G/H$ vanishes) and not vanish (by the
Gauss-Bonnet Theorem). (Other results in the same spirit, obtaining a
contradiction from the study of characteristic classes of $\Gamma \backslash
G/H$, appear in \cite{KobayashiOno}.)

\begin{thm}[{(Kulkarni)}] \label{SO2n/SO1odd-notess}
 If $n$~is odd, then $\SO(2,n)/\SO(1,n)$ does not have a
tessellation.
 \end{thm}

Let us give a more explicit description of the closed, connected subgroups~$H$
of $\SO(2,2m)$ or $\SU(2,2m)$, such that $G/H$ has a tessellation. This
shows that if $n$~is even, then the Kulkarni-Kobayashi examples
\pref{KulkarniEg} and certain deformations are essentially the only
interesting homogeneous spaces of $\SO(2,n)$ or $\SU(2,n)$ that have
tessellations.

\begin{notation}[{\cite[Thm.~6.5.1, pp.~270--271]{HelgasonBook},
\cite[p.~180]{Hochschild-Lie}}] \label{KANDefn}
 Fix an Iwasawa decomposition  $G = KAN$. Thus,
 \begin{itemize}
 \item $K$ is a maximal compact subgroup,
 \item $A$ is the identity component of a maximal split torus, and
 \item $N$ is a maximal unipotent subgroup.
 \end{itemize}
 \end{notation}

The following two results are stated only for subgroups of~$AN$, because the
general case reduces to this \see{HcanbeAN}. The reason is basically that $H$
contains a connected, cocompact subgroup that is conjugate to a subgroup
of $AN$. (Clearly, if $H'$ is any cocompact subgroup of~$H$, then $G/H$ has a
tessellation if and only if $G/H'$ has a tessellation.) This is not quite
true in general, but the following lemma provides a satisfactory substitute,
by showing that it becomes true after enlarging $H$ by a compact amount.

\begin{thmref}{HcanbeAN}
 \begin{lem}
 After replacing $H$ by a conjugate subgroup, there is a closed, connected
subgroup~$H^*$ of $G$, such that $H^*/H$ and $H^*/(AN \cap H^*)^\circ$ are
compact, where $(AN \cap H^*)^\circ$ denotes the identity component of $AN
\cap H^*$.
 \end{lem}
 \end{thmref}

\begin{thmref}{SUFevenTess}
 \begin{thm}[{(Oh-Witte \cite[Thm.~1.7]{OhWitte-CK})}] \label{SOevenTess}
 Assume $G = \SO(2,2m)$, and let $H$ be a closed, connected, nontrivial,
proper subgroup of $AN$.

 The homogeneous space $G/H$ has a tessellation if and only if $H$ is
conjugate to a subgroup~$H'$, such that either
 \begin{enumerate}
 \item \label{SOevenTess-SO}
 $H' = \SO(1,2m) \cap AN$; or
 \item \label{SOevenTess-Sp}
 $H'$ belongs to a certain family~$\{H_B\}$ of deformations of $\SU(1,m) \cap
AN$, described explicitly in Theorem~\ref{HBthm} {\upshape(}with $\F
= \real${\upshape)}.
 \end{enumerate}
 \end{thm}
 \end{thmref}

\begin{thmrefref}{SUFevenTess}
\begin{thm} \label{SUevenTess}
 Assume $G = \SU(2,2m)$, and let $H$ be a closed, connected, 
nontrivial, proper
subgroup of $AN$.

 The homogeneous space $G/H$ has a tessellation if and only if $H$ is
conjugate to a subgroup~$H'$, such that either
 \begin{enumerate}
 \item \label{SUevenTess-SU}
 $H'$ belongs to a certain family~$\{\Hc\}$ of deformations of
$\SU(1,2m) \cap AN$, described explicitly in Theorem~\ref{SUegs}; or
 \item \label{SUevenTess-Sp}
 $H'$ belongs to a certain family~$\{H_B\}$ of deformations of $\Symp(1,m) \cap
AN$, described explicitly in Theorem~\ref{HBthm} {\upshape(}with $\F
= \complex${\upshape)}.
 \end{enumerate}
 \end{thm}
 \end{thmrefref}

The proof of Theorem~\ref{SOevenTess} (or~\ref{OW-known}) in \cite{OhWitte-CK}
requires a list \cite{OhWitte-CDS} of all the homogeneous spaces of
$\SO(2,n)$ that admit a proper action of a noncompact subgroup of $\SU(2,n)$.
(The list was obtained by very tedious case-by-case analysis. It was extended
to homogeneous spaces of $\SU(2,n)$ in \cite{IozziWitte-CDS}.)
The following proposition \pref{tessAN->dim>1,2} provides an \emph{a~priori}
lower bound on $\dim H$, and it turns out that the classification of the
interesting subgroups of large dimension can be achieved fairly easily
(see~\S\ref{SUFlargeSect}). This is the main reason that we are able to give
reasonably short complete proofs of Theorems~\ref{OW-known}, \ref{IW-known},
\ref{SU1m->complete}, \ref{SOevenTess}, and~\ref{SUevenTess}.

 \begin{prop}[{(see \ref{tess->dim>1,2}, \ref{mu(SUorSp)}, and~\ref{d(Sp)})}]
\label{tessAN->dim>1,2}
 Suppose $G = \SO(2,n)$ or $\SU(2,n)$, and
 let $H$ be a closed, connected, nontrivial subgroup of $AN$.
 If $G/H$ has a tessellation, then
 $$\dim H \ge
 \begin{cases}
 n & \hbox{if $G = \SO(2,n)$ and $n$ is even}; \\
 n-1 & \hbox{if $G = \SO(2,n)$ and $n$ is odd}; \\
 2n & \hbox{if $G = \SU(2,n)$ and $n$ is even}; \\
 2n-2 & \hbox{if $G = \SU(2,n)$ and $n$ is odd}.
 \end{cases}
 $$
 \end{prop}

\begin{ack}
 This research was partially supported by a grant from the National Science
Foundation (DMS-9801136).
 Much of the work was carried out during productive visits to the Isaac Newton
Institute for Mathematical Sciences (Cambridge, U.K.); we would like to thank
the Newton Institute for the financial support that made the visits possible.
 A.I.\ would like to thank the mathematics department of Oklahoma State
University for their warm and generous hospitality and Marc Burger
both for pointing out a mistake in the original statement 
and proof of Theorem~\ref{SUFevenTess} and for many enlightening
conversations.  
 \end{ack}

\section{Cartan projection and Cartan-decomposition subgroups}
\label{CartanSect}

The main problem in this paper is to determine whether or not a homogeneous
space $G/H$ has a tessellation. This requires some method to determine
whether or not a given discrete subgroup~$\Gamma$ of~$G$ acts properly
discontinuously on $G/H$. Y.~Benoist and T.~Kobayashi (independently)
demonstrated that the Cartan projection~$\mu$ is an effective tool to study
this question. It is the foundation of almost all of our work in later
sections.

In this section, we introduce the Cartan projection, and describe some of its
basic properties. First, however, we recall the notion of a proper action
(a generalization of properly discontinuous actions) and of a
Cartan-decomposition subgroup. At the end of the section, we use the Cartan
projection to briefly discuss the question of when there is a loss of
generality in assuming that $G$ is connected.

\subsection{Proper actions}

\begin{defn}[{(\cite[Defn.~2.11]{Kobayashi-survey97},
\cite[Defn.~1.2.2, (6)]{Palais-Slice}}]
 A topological group~$L$ of homeomorphisms of a topological space~$M$ acts
\emph{properly} on~$M$ if, for every compact subset~$C$
of~$M$,
 $$ \mbox{ $\{\, g \in L \mid C \cap g C \neq \emptyset \,\}$
 is compact.} $$
 \end{defn}

\begin{rem}
 It is important to note that a \emph{discrete} group of homeomorphisms
of~$M$ acts properly on~$M$ if and only if it acts properly discontinuously
on~$M$.
 \end{rem}

For the special case where $M = G/H$ is a homogeneous space, the following
lemma restates the definition of a proper action in more group-theoretic
terms.

\begin{lem}[{\cite[Obs.~2.1.3]{Kobayashi-criterion}}] \label{proper<>CHC}
 A closed subgroup~$L$ of~$G$ acts properly on $G/H$ if and only if, for
every compact subset~$C$ of~$G$, the intersection $L \cap (CHC)$ is compact.
 \end{lem}

\begin{proof} If $C$ is any compact subset of~$G$, then
$\overline{C} = CH/H$ is a compact subset of $G/H$; furthermore, any compact
subset of $G/H$ is contained in one of the form~$\overline{C}$. We have
 \begin{align*}
 \{\, g \in L \mid \overline{C} \cap g \overline{C} \neq
\emptyset \,\}
 &=  \{\, g \in L \mid (CH) \cap (g CH) \neq \emptyset \,\} \\
 &=  \{\, g \in L \mid g \in (CH) (CH)^{-1} \,\} \\
 &= L \cap (C H C^{-1}) .
 \end{align*}
 \end{proof}

This has the following well-known, easy consequence.

\begin{cor}[{(cf.\ \cite[Lem.~2.2(2)]{Kobayashi-criterion})}]
\label{CHCproper}
 Suppose $H$, $H_1$, $L$, and~$L_1$ are closed subgroups of~$G$.
 If 
 \begin{itemize}
 \item $L$ acts properly on $G/H$, and 
 \item there is a compact subset~$C$ of~$G$, such that $H_1 \subset CHC$ and
$L_1 \subset CLC$, 
 \end{itemize}
 then $L_1$ acts properly on $G/H_1$.
 \end{cor}

\subsection{Cartan-decomposition subgroups}

The following definition describes the subgroups to which the
Calabi-Markus Phenomenon applies (cf.\ Example~\ref{CalabiMarkusCircle}).

\begin{defn}
 We say that $H$ is a \emph{Cartan-decomposition subgroup} of~$G$ if $H \sim
G$ (see Notation~\ref{simdefn}).
 \end{defn}

\begin{rem} \label{AisCDS}
 From the Cartan decomposition $G = KAK$, we know that $A$ is a
Cartan-decomposition subgroup.
 \end{rem}

\begin{rem} \label{CDSconj}
 Any conjugate of a Cartan-decomposition subgroup is a Cartan-decomposition
subgroup.
 \end{rem}

\begin{lem}[{(Calabi-Markus Phenomenon, cf.~\cite[pf.~of
Thm.~A.1.2]{Kulkarni})}] \label{CalabiMarkus}
 If $H$ is a Cartan-decomposition subgroup of~$G$, and $\Gamma$ is a
discrete subgroup of~$G$ that acts properly discontinuously on $G/H$, then
$\Gamma$ is finite.
 \end{lem}

\begin{proof}
 Because $H$ is a Cartan-decomposition subgroup, there is a compact
subset~$C$ of~$G$, such that $CHC = G$. However, from
Lemma~\ref{proper<>CHC}, we know that $\Gamma \cap (CHC)$ is finite.
 Therefore
 $$ \Gamma = \Gamma \cap G = \Gamma \cap (CHC) $$
 is finite.
 \end{proof}

The following well-known, easy fact is a direct consequence of the
Calabi-Markus Phenomenon. It is an important first step toward determining
which homogeneous spaces have tessellations.

\begin{cor}  \label{CDS->notess}
 If $H$ is a Cartan-decomposition subgroup of~$G$, such that $G/H$ is not
compact, then $G/H$ does not have a tessellation.
 \end{cor}

\subsection{The Cartan projection}

\begin{notation}[{\cite[\S9.1, p.~402]{HelgasonBook}}] \label{A+Defn}
 \begin{itemize}
 \item If $G$ is connected, let $A^+$ be the (closed) positive Weyl chamber
of~$A$ in which the roots occurring in the Lie algebra of~$N$ are positive
\cf{KANDefn}. Thus, $A^+$ is a fundamental domain for the action of the
(real) Weyl group of~$G$ on~$A$.

 \item In the general case, let $A^+$ be a closed, convex fundamental domain
for the action of the (real) Weyl group of~$G$ on~$A$, such that $A^+$ is
contained in the (closed) positive Weyl chamber of~$A$ in which the roots
occurring in the Lie algebra of~$N$ are positive.
 \end{itemize}
 \end{notation}

\begin{defn}[{\cite[Thm.~9.1.1, p.~402]{HelgasonBook}, \cite{Benoist,
Kobayashi-survey97}}]
 For each element~$g$ of~$G$, the Cartan decomposition $G = K A^+ K$ implies
that there is an element~$a$ of~$A^+$ with $g \in K a K$. In fact, the
element~$a$ is unique, so there is a well-defined function $\mu \colon G \to
A^+$ given by $g \in K \, \mu(g) \, K$. 

We remark that the function $\mu$ is continuous and proper (that is, the
inverse image of any compact set is compact).
 \end{defn}

The following crucial result of Y.~Benoist provides a uniform estimate on the
variation of~$\mu$ over disks of bounded radius. (A related result was proved,
independently and simultaneously, by T.~Kobayashi
\cite[Thm.~3.4]{Kobayashi-criterion}.) The proof is both elementary and
elegant. However, it requires a bit of notation, so we postpone it to
\S\ref{CalcSect-Benoist} (and, for concreteness, we will assume that $G$ is
either $\SO(2,n)$ or $\SU(2,n)$ in the proof).

\begin{prop}[{(Benoist \cite[Prop.~5.1]{Benoist})}] \label{bddchange}
 For any compact subset~$C$ of~$G$, there is a compact subset~$C'$ of~$A$,
such that $\mu(C g C) \subset \mu(g) C'$, for all $g \in G$.
 \end{prop}

\begin{notation}
 For subsets $U$ and~$V$ of~$A^+$, we write $U \approx V$ if there 
is a compact
subset~$C$ of~$A$, such that $U \subset VC$ and $V \subset UC$.
 \end{notation}

\begin{cor}[{(Benoist \cite[Prop.~5.1]{Benoist}, Kobayashi
\cite[Thm.~1.1]{Kobayashi-criterion})}] \label{simvsmu}
 For any subgroups $H_1$ and~$H_2$ of~$G$, we have $H_1 \sim H_2$ if and
only if $\mu(H_1) \approx \mu(H_2)$.
 \end{cor}

\begin{proof}
 ($\Rightarrow$) Let $C$ be a compact subset of~$G$, such that $H_1 \subset
CH_2 C$ and $H_2 \subset CH_1 C$. Choose a corresponding compact subset~$C'$
of~$A$, as in Proposition~\ref{bddchange}. Then 
 $$\mu(H_1) \subset \mu(CH_2 C) \subset \mu(H_2) C'$$
 and, similarly,
 $\mu(H_2) \subset \mu(H_1) C'$.

($\Leftarrow$) Let $C$ be a compact subset of~$A$, such that $\mu(H_1) \subset
\mu(H_2) C$ and $\mu(H_2) \subset \mu(H_1) C$. Then
 $$ H_1 \subset K \, \mu(H_1) \, K
 \subset K \, (\mu(H_2) C) \, K
 \subset K \, \bigl( (KH_2 K) C \bigr) \, K $$
 and, similarly, $H_2 \subset K H_1 (KCK)$.
 \end{proof}

The special case where $H_2 = G$ (and $H_1$ is closed and almost connected)
can be restated as follows.

\begin{cor}[{(Benoist, Kobayashi)}] \label{CDSvsmu}
 $H$ is a Cartan-decomposition subgroup of~$G$ if and only if $\mu(H) \approx
A^+$.
 \end{cor}

\begin{cor}[{(cf.~\cite[Lem.~3.2]{Kobayashi-isotropy})}]
\label{Rrank1-CDS}
 Assume that $\Rrank G = 1$. The subgroup~$H$ is a Cartan-decomposition
subgroup of~$G$ if and only if $H$ is noncompact.
 \end{cor}

\begin{proof}
 ($\Leftarrow$) We have $\mu(e) = e$, and, because $\mu$ is a proper map,
we have $\mu(h) \to \infty$ as $h \to \infty$ in~$H$. Because $\Rrank G =
1$, we know that $A^+$ is homeomorphic to the half-line $[0,\infty)$
(with the point~$e$ in~$A^+$ corresponding to the endpoint~$0$ of the
half-line), so, by continuity, it must be the case that $\mu(H) = A^+$. Then
Corollary~\ref{CDSvsmu} implies that $H$ is a Cartan-decomposition subgroup,
but we provide the following direct proof that avoids any appeal to
Proposition~\ref{bddchange}.

From the definition of~$\mu$, we have $KHK = K \, \mu(H) \, K$.
 Therefore
 $$ KHK = K \, \mu(H) \, K = K A^+ K = G ,$$
 so $H$ is a Cartan-decomposition subgroup (by taking $C = K$ in
Definition~\ref{simdefn}).
 \end{proof}

By using Lemma~\ref{proper<>CHC}, the proof of Corollary~\ref{simvsmu} also
establishes the following.

\begin{cor}[{(Benoist \cite[Prop.~1.5]{Benoist}, Kobayashi
\cite[Cor.~3.5]{Kobayashi-criterion})}] \label{proper<>mu(L)}
 Suppose $H$ and $L$ are closed subgroups of~$G$. The subgroup $L$ acts
properly on $G/H$ if and only if $\mu(L) \cap \mu(H) C$ is compact, for every
compact subset~$C$ of~$A$.
 \end{cor}

\subsection{Disconnected groups} \label{disconnected} 

As was mentioned in Remark~\ref{Hdisconnected}, we may assume, without loss
of generality, that $H$ is connected. However, it may not be possible to
assume that $G$ is connected, because, although there are no known examples,
it is possible that the following question has an affirmative answer.

\begin{ques} \label{assumeGconn?}
 Does there exist a homogeneous space $G/H$ (satisfying
Assumption~\ref{standing}), such that $G^\circ/H^\circ$ has a tessellation,
but $G/H$ does not have a tessellation?
 \end{ques}

If $\Gamma$ is a crystallographic group for $G^\circ/H^\circ$, then it is
easy to see that $\Gamma \backslash G/H$ is compact. However, the following
example shows that $\Gamma$ may not act properly discontinuously on $G/H$.

\begin{eg}
 Let 
 \begin{itemize}
 \item $L = H = \SL(2,\real)$, 
 \item $\sigma$ be the automorphism of $L \times H$ that interchanges the two
factors (that is, $\sigma(x,y) = (y,x)$), 
 \item $G = (L \times H) \rtimes \langle \sigma \rangle $ (semidirect
product), and
 \item $\Gamma$ be a cocompact lattice in~$L$ \fullcf{classical}{Borel}.
 \end{itemize}
 Then $H = H^\circ$, and $\Gamma$ is a crystallographic group for $G^\circ/H
= (L \times H)/H$ (see Example~\ref{G=LxH}).

However, $\Gamma \subset L = \sigma^{-1} H \sigma$, so $\Gamma$ does not act
properly on $G/H$ \seesee{proper<>CHC}{with $C = \{\sigma, \sigma^{-1}\}$}.

Even so, $G/H$ does have a tessellation, because the diagonal
embedding
 $$\Delta(\Gamma) = \{\, (\gamma,\gamma) \in L \times H \mid \gamma \in
\Gamma \,\} $$
 is a crystallographic group for $G/H$. Thus, this example does not provide
an answer to Question~\ref{assumeGconn?}.
 \end{eg}

In this example, $\sigma$ represents an element of the Weyl group of~$G$ that
does not belong to the Weyl group of~$G^\circ$. The following proposition
shows that this is a crucial ingredient in the construction.

\begin{prop}
 Let $\Gamma$ be a crystallographic group for $G^\circ/H^\circ$.

If the {\upshape(}real{\upshape)} Weyl group of $G$ is same as the
{\upshape(}real{\upshape)} Weyl group of~$G^\circ$, then $\Gamma$ is a
crystallographic group for $G/H$.
 \end{prop}

\begin{proof}
 By assumption, we may choose the same fundamental domain~$A^+$ for the Weyl
groups of $G$ and~$G^\circ$. Let $\mu \colon G \to A^+$ and $\mu^\circ \colon
G^\circ \to A^+$ be the Cartan projections; then $\mu^\circ$ is the
restriction of~$\mu$ to~$G^\circ$. For simplicity, assume, without loss of
generality, that $H \subset G^\circ$ (for example, assume $H$ is connected).
Then, for any compact subset~$C$ of~$A$, we have
 $$ \mu(\Gamma) \cap \mu(H) C = \mu^\circ(\Gamma) \cap \mu^\circ(H) C $$
 is finite \see{proper<>mu(L)}. Thus, $\Gamma$ acts properly discontinuously
on $G/H$ \see{proper<>mu(L)}, as desired.
 \end{proof}

 A.~Borel and J.~Tits \cite[Cor.~14.6, p.~147]{BorelTits-Reductive} proved
that if $G$ is Zariski connected, then every element of the Weyl group of~$G$
has a representative in~$G^\circ$. Also, any element of the Weyl group must
act as an automorphism of the root system. Thus, we have the following
corollary.

\begin{cor}
 Let $\Gamma$ be a crystallographic group for $G^\circ/H^\circ$.

If either
 \begin{itemize}
 \item $G$ is Zariski connected, or
 \item every automorphism of the real root system of $G^\circ$ belongs to the
Weyl group of the root system, 
 \end{itemize}
 then $\Gamma$ is a crystallographic group for $G/H$.
 \end{cor}

\begin{eg}
 \begin{enumerate}
 \item If $G = \SO(2,n)$, then $G$ is Zariski connected (because
$\SO(n+2,\complex)$ is connected \cite[Thm.~2.1.9, p.~60]{GoodmanWallach}), so
$G/H$ has a tessellation if and only if $G^\circ/H^\circ$ has a tessellation.
 \item More generally, if $G^\circ = \SO(2,n)^\circ$ or $\SU(2,n)$ (with $n
\ge 3$), then every automorphism of the real root system of $G^\circ$ belongs
to the Weyl group of the root system (cf.\ Figure~\ref{rootspict}), so $G/H$
has a tessellation if and only if $G^\circ/H^\circ$ has a tessellation.
 \end{enumerate}
 \end{eg}

\begin{eg}
  If $G = \SL(3,\real) \rtimes \langle \sigma \rangle$, where $\sigma$ is the
Cartan involution of $\SL(3,\real)$, then $\sigma$ represents an element of
the Weyl group of~$G$ that does not belong to the Weyl group of~$G^\circ$, so
the proposition does not apply to~$G$. However, this does not matter: if
neither $H$ nor $G/H$ is compact, then Theorem~\ref{SL3->notess} implies that
$G^\circ/H^\circ$ has no tessellations, so $G/H$ has no tessellations either.
 \end{eg}

\section{Preliminaries on subgroups of $AN$} \label{ANSect}

This section recalls a technical result that often
allows us to assume that $H$ is a subgroup of~$AN$. It also recalls some
basic topological properties of such subgroups, and also recalls a simple
observation relating these subgroups to the root spaces of the Lie
algebra~$\Lie G$.

\subsection{Reduction to subgroups of $AN$}

\begin{defn}[{\cite[Thm.~9.7.2, p.~431]{HelgasonBook}}]
 An element~$g$ of~$G$ is:
 \begin{itemize}
 \item \emph{hyperbolic} if $g$ is conjugate to an element of~$A$;
 \item \emph{unipotent} if $g$ is conjugate to an element of~$N$; 
 \item \emph{elliptic} if $g$ is conjugate to an element of~$K$.
 \end{itemize}
 \end{defn}

\begin{lem}[{(Real Jordan Decomposition \cite[Lem.~9.7.1,
p.~430]{HelgasonBook})}] \label{JordanDecomp}
 Each $g \in G$ has a unique decomposition in the form $g = auc$, such that 
 \begin{itemize}
 \item $a$ is hyperbolic, $u$ is unipotent, and $c$~is elliptic; and
 \item $a$, $u$, and~$c$ all commute with each other.
 \end{itemize}
 \end{lem}

\begin{rem} \label{JordanCommute}
 If $g = auc$ is the Real Jordan Decomposition of some element~$g$ of~$G$,
then $a$, $u$, and~$c$ commute, not only with each other, but also with
any element of~$G$ that commutes with~$g$. This is because the Real Jordan
Decomposition of $h^{-1} g h$ is
 $$ h^{-1} g h = (h^{-1} a h) (h^{-1} u h) (h^{-1} c h) :$$
 if $h^{-1} g h = g$, then the uniqueness of the Real Jordan Decomposition
of~$g$ implies $h^{-1} a h = a$, $h^{-1} u h = u$, and $h^{-1} g h = c$.
 \end{rem}

The following observation is a generalization of the fact that a collection of
commuting triangularizable matrices can be simultaneously triangularized.

\begin{lem}[{(cf.\ pf.\ of \cite[Thm.~17.6]{Humphreys-Algic})}]
\label{simultaneous}
 If $H$ is abelian {\upshape(}or, more generally, solvable{\upshape)}, and is
generated by hyperbolic and/or unipotent elements, then $H$ is conjugate to a
subgroup of~$AN$.
 \end{lem}

Because of the following result, we usually assume $H \subset AN$ (by
replacing $H$ with a conjugate of~$H'$).

\begin{lem}[{(cf.~\cite[Lem.~2.9]{OhWitte-CDS})}] \label{HcanbeAN}
 If $H$ is connected, then there is a closed, connected subgroup~$H'$ of~$G$
and a compact, connected subgroup~$C$ of~$G$, such that
 \begin{enumerate}
 \item \label{HcanbeAN-inAN}
 $H'$ is conjugate to a subgroup of~$AN$;
 \item \label{HcanbeAN-CH=CH'}
 $C H = C H'$ is a subgroup of~$G$; and
 \item \label{HcanbeAN-d(H')}
 $d(H') = d(H)$ {\upshape(}see  Notation~\ref{d(H)-defn}{\upshape)}.
 \end{enumerate}
 Moreover, it is easy to see from~\pref{HcanbeAN-CH=CH'} that the homogeneous
space $G/H$ has a tessellation if and only if $G/H'$ has a tessellation.
 \end{lem}

\begin{proof}[Idea of proof.]
 First, let us note that every connected subgroup of $AN$ is closed
(see~\fullref{solvable}{H=Rn} and~\ref{ANsc}), so we do not need to show that
$H'$ is closed.

Second, let us note that \pref{HcanbeAN-d(H')} is a consequence
of~\pref{HcanbeAN-inAN} and~\pref{HcanbeAN-CH=CH'}. To see this, let $K^*$ be
a maximal compact subgroup of~$CH$ that contains~$C$. Then a standard
argument shows that $K^* \cap H$ is a maximal compact subgroup of~$H$.
(Because all maximal compact subgroups of $CH$ are conjugate, there is some
$g \in CH$, such that $(g^{-1} K^* g) \cap H$ is a maximal compact subgroup
of~$H$ that contains $K^* \cap H$. Since $C \subset K^*$, we know that $C$
normalizes~$K^*$, so we may assume $g \in H$; thus, $g$ normalizes~$H$. Then 
 $g^{-1} (K^* \cap H) g = (g^{-1} K^* g) \cap H$
 contains $K^* \cap H$. Because $K^* \cap H$ is compact, this implies that
$g$ normalizes $K^* \cap H$. So $K^* \cap H = (g^{-1} K^* g) \cap H$ is a
maximal compact subgroup of~$H$.) Therefore
 $$ \dim (K^*H/K^*) = \dim \bigl( H/ (K^* \cap H) \bigr)
 = \dim H - \dim(K^* \cap H)
 = d(H) .$$
 Similarly, $ \dim (K^*H'/K) = d(H')$. Since $K^*H = CH = CH' = K^*H'$, we
conclude that $d(H') = d(H)$, as desired.

\setcounter{case}{0}

\begin{case} \label{HcanbeAN-ss}
 Assume $H$ is semisimple.
 \end{case}
 We have an Iwasawa decomposition $H = K_H A_H N_H$; let $H' = A_H N_H$ and
$C = K_H$.

\begin{case}
 Assume $H = \{h^t\}$ is a one-parameter subgroup.
 \end{case}
 Let 
 \begin{itemize}
 \item $h^t = a^t u^t c^t$ be the Real Jordan Decomposition of~$h^t$
\see{JordanDecomp};
 \item $H' = \{ a^t u^t\}$; and
 \item $C = \overline{\{c^t\}}$ be the closure of~$\{c^t\}$.
 \end{itemize}
 (Lemma~\ref{simultaneous} implies that $H'$ is conjugate to a subgroup
of~$AN$.)

\begin{case} \label{HcanbeAN-abel}
 Assume $H$ is abelian.
 \end{case}
 We may write $H$ as a product of one-parameter subgroups:
 $$ H = \{\, h_1^{t_1} h_2^{t_2} \cdots h_r^{t_r} \mid t_1,\ldots,t_r \in
\real \,\} .$$
 Let $h_j^t = a_j^t u_j^t c_j^t$ be the Real Jordan Decomposition of~$h_j^t$
\see{JordanDecomp}. Note that
$a_j^{t_j}$, $u_j^{t_j}$, and~$c_j^{t_j}$ commute, not only with each other,
but also with every  $a_k^{t_k}$, $u_k^{t_k}$, and~$c_k^{t_k}$
\see{JordanCommute}. Let
 $$ H' = \{\, (a_1^{t_1} u_1^{t_1}) (a_2^{t_2} u_2^{t_2}) \cdots (a_k^{t_r}
u_k^{t_r}) \mid t_1,\ldots,t_r \in \real \,\} ,$$
 and let $C = \overline{\{c_1^t\}} \cdots \overline{\{c_1^t\}}$.
 (Lemma~\ref{simultaneous} implies that $H'$ is conjugate to a subgroup
of~$AN$.)

\begin{case}
 The general case.
 \end{case}
 From the Levi decomposition \cite[p.~91]{JacobsonBook}, we know that
there is a connected, semisimple subgroup~$L$ of~$H$ and a connected,
solvable, normal subgroup~$R$ of~$H$, such that $H = LR$ (and $L \cap R$ is
finite). Let $U = [H,R]$, so $U$ is a connected, normal subgroup of~$H$, and
$U$ is conjugate to a subgroup of~$N$ (cf.\ \cite[Cor.~2.7.1,
p.~51]{JacobsonBook}). By modding out~$U$, we (essentially) reduce to the
direct product of Cases~\ref{HcanbeAN-ss} and~\ref{HcanbeAN-abel}.
 \end{proof}

\begin{rem}
 For $H$ and $H'$ as in Lemma~\ref{HcanbeAN}, Proposition~\ref{AN/H=Rd}
(and~\ref{ANsc}) implies that if $H' \neq AN$, then $AN/H'$ is not compact;
also, Proposition~\fullref{solvable}{nocpct} (and~\ref{ANsc}) implies that if
$H' \neq e$, then $H'$ is not compact. Therefore:
 \begin{itemize}
 \item $H' = AN$ if and only if $G/H$ is compact; and
 \item $H' = e$ if and only if $H$ is compact.
 \end{itemize}
 Thus, if neither $H$ nor $G/H$ is compact, then $H'$ is a nontrivial, proper
subgroup of $AN$.
 \end{rem}

\subsection{Topology of solvable groups and their homogeneous spaces}

Everything is this subsection is well known, though somewhat scattered in the
literature. The main results are Propositions~\ref{solvable}
and~\ref{AN/H=Rd}, which, together with Corollary~\ref{ANsc}, show that
connected subgroups of $AN$ and their homogeneous spaces are very well behaved
topologically. Corollary~\ref{fiberbundle}, on the homology of very simple
quotient spaces, is also used in later sections.

We begin with the easy case of abelian groups. This lemma generalizes almost
verbatim to solvable groups \see{solvable}, but the proof in that generality
is not as trivial.

\begin{lem}[{\cite[Thm.~3.6.2, p.~196]{Varadarajan}}] \label{abelian}
 Let $R$ be a $1$-connected, abelian Lie group. 
 \begin{enumerate}
 \item \label{abelian-H=Rn}
 If $H$ is a connected subgroup of~$R$, then $H$ is closed, simply connected,
and isomorphic to~$\real^k$, for some~$k$.
 \item \label{abelian-HcapL}
 If $H$ and $L$ are connected subgroups of~$R$, then $H \cap L$ is connected.
 \item \label{abelian-nocpct}
 If $C$ is a compact subgroup of~$R$, then $C$ is trivial.
 \end{enumerate}
 \end{lem}

\begin{proof}
 Because $R$ is abelian and 1-connected, the exponential map is a Lie group
isomorphism from the additive group of the Lie algebra~$\Lie R$ onto~$R$.

\pref{abelian-H=Rn} Let $k = \dim H$. Because the exponential map is a Lie
group isomorphism (hence a diffeomorphism), and because $\Lie H$ is a closed
$k$-submanifold of~$\Lie R$, we know that $\exp(\Lie H)$ is a closed
$k$-submanifold of $R$. Of course, $\exp(\Lie H)$ is contained in~$H$, which
is also a $k$-submanifold of $R$. Because the dimensions are the same, we
know that $\exp(\Lie H)$ is open in~$H$. Also, because $\exp(\Lie H)$ is
closed in $R$, we know that $\exp(\Lie H)$ is closed in~$H$. Therefore 
 \begin{equation} \label{abelianPf-H=Rn-exp(H)}
 \exp(\Lie H) = H
 \end{equation} (because $H$ is connected). Finally, we know that
$\exp|_{\Lie H}$ is a diffeomorphism from its domain~$\Lie H \diffeo \real^k$
onto its image~$H$.

\pref{abelian-HcapL} From \pref{abelianPf-H=Rn-exp(H)}, we have $\exp(\Lie H)
= H$ and, similarly, $\exp(\Lie L) = L$. Also, because $\exp$ is bijective,
we have 
 $\exp \Lie H \cap \exp \Lie L = \exp( \Lie H \cap \Lie L ) $.
 Therefore
 $$ H \cap L = \exp \Lie H \cap \exp \Lie L = \exp( \Lie H \cap \Lie L ) $$
 is connected.

\pref{abelian-nocpct}  Because $\real^k$ is not compact (for $k > 0$), we
know, from~\fullref{abelian}{H=Rn}, that $C^\circ$ is trivial; so $C$ is
finite. Since $R \iso (\Lie R, +) \iso \real^d$ has no elements of finite
order, we conclude that $C$ is trivial.
 \end{proof}

As is usual in the theory of solvable groups, the main results of this
section are proved by induction, based on modding out some normal
subgroup~$L$. To be effective, this method requires an understanding of the
quotient space $R/L$. The information we need (even if $L$ is not normal)
comes from the following elementary observation, because $R$ is a principal
$L$-bundle over $R/L$.

\begin{lem} \label{trivialbundle}
 Let $P$ be a principal $H$-bundle over a manifold~$M$. 
 \begin{enumerate}
  \item \label{trivialbundle-H=Rn}
 If $H$ is diffeomorphic to $\real^n$, then 
  \begin{enumerate}
   \item $P$ is $H$-equivariantly diffeomorphic to $M \times H$, so
   \item $P$ is homotopy equivalent to~$M$.
  \end{enumerate}
  \item \label{trivialbundle-M=Rn}
 If $M$ is diffeomorphic to $\real^n$, then 
  \begin{enumerate}
   \item $P$ is $H$-equivariantly diffeomorphic to $M \times H$, so
   \item $P$ is homotopy equivalent to~$H$.
  \end{enumerate}
 \end{enumerate}
 \end{lem}

\begin{proof}
 Any principal bundle with a section is trivial \cite[Cor.~4.8.3,
p.~48]{Husemoller}. If either the fiber or the base is contractible, then
there is no obstruction to constructing a section \cite[Thm.~2.7.1(H1),
p.~21]{Husemoller}, so $P$ is trivial: $P \diffeo M \times H$. (The
diffeomorphism can be taken to be $H$-equivariant, with respect to the
natural $H$-action on $M \times H$, given by $(m,h)h' = (m,h h')$.) Then the
conclusions on homotopy equivalence follow from the fact that $\real^n$ is
contractible (that is, homotopically trivial).
 \end{proof}

We recall the long exact sequence of the fibration $H \to R \to R/H$:

\begin{lem}[{\cite[Cor.~IV.8.6, p.~187]{WhiteheadBook}}]
\label{HtpyExactFibration}
 Let $H$ be a closed subgroup of a Lie group~$R$. There is a
{\upshape(}natural{\upshape)} long exact sequence of homotopy groups:
 $$ \cdots \to \pi_1(H) \to \pi_1(R) \to \pi_1(R/H) \to \pi_0(H) \to \pi_0(R)
\to \pi_0(R/H) \to 0 .$$
 \end{lem}

\begin{cor} \label{R/Hsc}
 Let $H$ be a closed subgroup of a $1$-connected Lie group~$R$. The
homogeneous space $R/H$ is simply connected if and only if $H$ is connected.
 \end{cor}

\begin{proof}
 Because $R$ is 1-connected, we have $\pi_1(R) = \pi_0(R) = 0$, so,
from~\pref{HtpyExactFibration}, we know that the sequence
 $$ 0 \to \pi_1(R/H) \to \pi_0(H) \to 0 $$
 is exact.
 Thus, $\pi_1(R/H) \iso \pi_0(H)$, so the desired conclusion is immediate.
 \end{proof}

As a step toward Proposition~\ref{solvable}, we prove two special cases that
describe the topology of normal subgroups.

\begin{lem} \label{R=Rd}
 If $R$ is a $1$-connected, solvable Lie group, then $R$ is diffeomorphic to
$\real^d$, for some~$d$.
 \end{lem}

\begin{proof}
 We may assume the group~$R$ is nonabelian (otherwise, the desired conclusion
is given by Lemma~\fullref{abelian}{H=Rn}).  Then, because $R$ is solvable,
there is a nontrivial, connected, proper, closed, normal subgroup~$L$ of~$R$.
Since $R/L$ is simply connected \see{R/Hsc}, and $\dim(R/L) < \dim R$, we may
assume, by induction on $\dim R$, that $R/L$ is diffeomorphic to
some~$\real^{d_1}$. Therefore 
 \begin{enumerate} \renewcommand{\theenumi}{\alph{enumi}}
 \item \label{R=RdPf-R=prod}
 $R$ is diffeomorphic to $(R/L) \times L$ and 
 \item \label{R=RdPf-L=R}
 $L$ is homotopy equivalent to~$R$ 
 \end{enumerate}
 \fullsee{trivialbundle}{M=Rn}. Because $R$ is $1$-connected,
\pref{R=RdPf-L=R}~implies that $L$ is $1$-connected; hence, $L$ is a
1-connected, solvable Lie group, so we may assume, by induction on $\dim R$,
that $L$ is diffeomorphic to some~$\real^{d_2}$. Thus,
\pref{R=RdPf-R=prod}~implies that $R$ is diffeomorphic to $\real^{d_1} \times
\real^{d_2} \diffeo \real^{d_1+d_2}$, as desired.
 \end{proof}

\begin{cor}[(of proof)] \label{Rnormal}
 If $R$ is a $1$-connected, solvable Lie group, then every connected, closed,
normal subgroup of~$R$ is $1$-connected.
 \end{cor}

The following proposition is a nearly complete generalization of
Lemma~\ref{abelian} to the class of solvable groups. There are two
exceptions:
 \begin{enumerate}
 \item[\protect\ref{solvable-H=Rn})] Of course, subgroups of a
solvable group may not be abelian, so the conclusion in
\fullref{abelian}{H=Rn} that $H$ is isomorphic to some~$\real^k$ must be
weakened to the conclusion that $H$ is diffeomorphic to some~$\real^k$.
 \item[\protect\ref{solvable-HcapL})] The intersection of connected subgroups
is not always connected \see{HcapLdisconn}, so we add the restriction that $L$
is normal to~\fullref{abelian}{HcapL}. (We remark that no such restriction is
necessary if $R \subset AN$, because the exponential map is a diffeomorphism
from~$\Lie R$ onto~$R$ in this case \cite{Dixmier-exp, Saito2}.)
 \end{enumerate}

\begin{eg} \label{HcapLdisconn}
 Let 
 $$ R = \bigset{
 \begin{pmatrix}
 e^{2\pi i t} & x + iy & 0 \\
 0            & 1      & 0 \\
 0            & 0      & e^t
 \end{pmatrix}
 }{
 t, x, y \in \real
 } ,
 \qquad
 g^t =
 \begin{pmatrix}
 e^{2\pi i t} & 0      & 0 \\
 0            & 1      & 0 \\
 0            & 0      & e^t
 \end{pmatrix}
 , $$
 $$
 u =
 \begin{pmatrix}
 1            & 1      & 0 \\
 0            & 1      & 0 \\
 0            & 0      & 0
 \end{pmatrix}
 , \qquad
 h^t 
 = u^{-1} g^t u
 =
 \begin{pmatrix}
 e^{2\pi i t} & e^{2\pi i t} - 1 & 0 \\
 0            & 1                & 0 \\
 0            & 0                & e^t
 \end{pmatrix}
 .
 $$
  Then $R$, being
diffeomorphic to~$\real^3$, is 1-connected; and $\{g^t\}$ and~$\{h^t\}$ are
connected subgroups. But
 $$ \{g^t\} \cap \{h^t\} =
 \bigset{ 
 \begin{pmatrix}
 1            & 0      & 0 \\
 0            & 1      & 0 \\
 0            & 0      & e^n
 \end{pmatrix}
 }{
 n \in \integer
 } $$
 is not connected.
 \end{eg}

\begin{prop}[{\cite[Thms.~12.2.2 and 12.2.3, pp.~137--138]{Hochschild-Lie}}]
\label{solvable}
 Let $R$ be a $1$-connected, solvable Lie group. 
 \begin{enumerate}
 \item \label{solvable-H=Rn}
 If $H$ is a connected subgroup of~$R$, then $H$ is closed, simply connected,
and diffeomorphic to some~$\real^d$.
 \item \label{solvable-HcapL}
 If $H$ and $L$ are connected subgroups of~$R$, and $L$ is normal, then
$H \cap L$ is connected.
 \item \label{solvable-nocpct}
 If $C$ is a compact subgroup of~$R$, then $C$ is trivial.
 \end{enumerate}
 \end{prop}

\begin{proof}
 \pref{solvable-HcapL}
 We may assume $L$ is nontrivial, so $\dim (R/L) < \dim R$. Thus, by
induction on $\dim R$, using \pref{solvable-H=Rn}, we may assume that $HL/L$
is a closed, simply connected subgroup of $R/L$. Then, since $H/(H \cap L)$
is homeomorphic to $HL/L$, we see that
 \begin{equation} \label{pi1(H/HcapL)}
 \mbox{$H/(H \cap L)$ is simply connected} ,
 \end{equation}
 so Lemma~\ref{R/Hsc} implies that $H \cap L$ is connected. 

\pref{solvable-H=Rn} Because $R$ is solvable, there is a connected, closed,
proper, normal subgroup~$L$ of~$R$, such that $R/L$ is abelian.
 We know that $L$ is 1-connected \see{Rnormal}, so, by induction on $\dim R$,
we may assume that every connected subgroup of~$L$ is closed and simply
connected. From~\pref{solvable-HcapL}, we know that $H \cap L$ is connected,
so we conclude that $H \cap L$ is closed, and 
 \begin{equation} \label{pi1(HcapL)}
 \pi_1(H \cap L) = 0 .
 \end{equation}
 From~\pref{HtpyExactFibration} (with $H$ in the place of~$R$, and $L$~in the
place of~$H$), together with \pref{pi1(H/HcapL)} and~\pref{pi1(HcapL)}, we
conclude that $\pi_1(H) = 0$; that is, $H$ is simply connected. So
\pref{R=Rd} implies $H$ is diffeomorphic to some~$\real^d$.

Because both $HL/L$ and $H \cap L$ are closed, it is not difficult to see
that $H$ is closed.

\pref{solvable-nocpct} Because $R$ is solvable, there is a connected, closed,
proper, normal subgroup~$L$ of~$R$, such that $R/L$ is abelian. We know that
$R/L$ is 1-connected \see{R/Hsc}, so $R/L$ has no nontrivial, compact
subgroups \fullsee{abelian}{nocpct}; thus, we must have $C \subset L$.
Therefore, $C$ is a compact subgroup of~$L$. Then, since $L$ is 1-connected
\see{Rnormal}, we may conclude, by induction on $\dim R$, that $C$ is
trivial.
 \end{proof}

\begin{cor} \label{ANsc}
 $AN$ is a $1$-connected, solvable Lie group.
 \end{cor}

\begin{proof}
 Because $G$ is linear, it is a subgroup of some $\GL(n,\real)$.
 Replacing $G$ by a conjugate, we may assume that $AN$ is contained in the
group~$B$ of upper triangular matrices with positive diagonal entries
\cf{simultaneous}. The matrix entries provide an obvious diffeomorphism
from~$B$ onto $(\real^+)^n \times \real^{n(n-1)/2} \diffeo \real^{n(n+1)/2}$,
so $B$ is 1-connected. Thus, Proposition~\fullref{solvable}{H=Rn} implies
that $AN$ is simply connected.
 \end{proof}

The following observation will be used in Sections~\ref{dimHSect}
and~\ref{ExistenceSect}.

\begin{cor} \label{fiberbundle}
 Let $F$ be a connected subgroup of~$AN$, and suppose we have a proper,
$C^\infty$ action of~$F$ on a manifold~$M$. Then $M$ and $M/F$ have the same
homology.
 \end{cor}

\begin{proof}
 Because the action is proper, we know that the stabilizer of each point
of~$M$ is compact. However, $F$ has no nontrivial compact subgroups
\fullsee{solvable}{nocpct}. Thus, the action is free.

 Because the action is free, proper, and $C^\infty$, it is easy to see that
the manifold $M$ is a principal fiber bundle over the quotient $M/F$
\cite[Thm.~1.1.3]{Palais-Slice}. Furthermore, the fiber~$F$ of the bundle is
contractible \fullsee{solvable}{H=Rn}, so
Lemma~\fullref{trivialbundle}{H=Rn} implies that $M$ homotopy equivalent
to~$M/F$. Therefore, the spaces $M$ and $M/F$ have the same homology.
 \end{proof}

For the special case where $M/F$ is a homogeneous space of a solvable group,
the following more detailed result describes the topology of $M/F$, not just
its homology.

\begin{prop}[{(Mostow \cite[Prop.~11.2]{Mostow-FSS})}] \label{AN/H=Rd} 
 If $H$ is any connected subgroup of a $1$-connected, solvable Lie group~$R$,
then $R/H$ is diffeomorphic to the Euclidean space $\real^d$, for some~$d$.
 \end{prop} 

\begin{proof}
 Because $R$ is solvable, it has a nontrivial, connected, closed, abelian,
normal subgroup~$L$. Since $L$ is abelian and $H \cap L$ is connected
\fullsee{solvable}{HcapL}, we know that $L/(H \cap L)$ is a 1-connected
abelian group \see{R/Hsc}, so it is isomorphic to some
$\real^{d_1}$ \fullsee{abelian}{H=Rn}. 

We know $H$ is closed \fullsee{solvable}{H=Rn}. Also, since $L$ is
nontrivial, we have $\dim (R/L) < \dim R$, so we may assume, by induction
on $\dim R$, that
 $$ R/(HL) \diffeo (R/L) / (HL/L) $$
 is diffeomorphic to some~$\real^{d_2}$.

Now $R$ is a principal $HL$-bundle over $R/(HL)$. Because $R/(HL) \diffeo
\real^{d_2}$, this bundle is trivial \fullsee{trivialbundle}{M=Rn}: $R$ is
$HL$-equivariantly diffeomorphic to $R/(HL) \times HL$. Then
 $$ R/H
 \diffeo R/(HL) \times HL/H
 \diffeo R/(HL) \times L/(H \cap L)
 \diffeo \real^{d_2} \times \real^{d_1}
 = \real^{d_1+d_2} ,$$
 as desired.
 \end{proof}

\subsection{$T$-invariant subspaces of~$\Lie A + \Lie N$}

The following well-known observation puts an important restriction on the
subspaces of $\Lie A + \Lie N$ that are normalized by a torus. It is an
ingredient in our case-by-case analysis of all possible subgroups of $AN$ in
Sections~\ref{SUFlargeSect} and~\ref{ProofSect}.

\begin{lem} \label{rootdecomp}
 Let
 \begin{itemize}
 \item $\Phi^+$ be the set of weights of~$A$ on~$\Lie N$ {\upshape(}in
other words, the set of all positive real roots of~$G${\upshape)};
 \item $T$ be a subgroup of~$A$;
 \item $\omega \in \Phi^+ \cup \{0\}$;
 \item $\Lie N^{=\omega} = \bigoplus_{\sigma|_T = \omega|_T} \rsp_\sigma$,
where the sum is over all $\sigma \in \Phi^+ \cup \{0\}$, such that the
restriction of~$\sigma$ to~$T$ is the same as the restriction of~$\omega$
to~$T$;
 \item $\Lie N^{\neq\omega} = \bigoplus_{\sigma|_T \neq \omega|_T}
\rsp_\sigma$, where the sum is over all $\sigma \in \Phi^+ \cup \{0\}$, such
that the restriction of~$\sigma$ to~$T$ is not the same as the restriction
of~$\omega$ to~$T$.
 \end{itemize}
 If $\Lie U$ is any $\real$-subspace of~$\Lie A + \Lie N$ normalized by~$T$,
then $\Lie U = (\Lie U \cap \Lie N^{=\omega}) \oplus (\Lie U \cap \Lie
N^{\neq\omega})$.
 \end{lem}

\begin{proof}
 Since $T \subset A$, we know that the elements of $\Ad_G T$ are
simultaneously diagonalizable (over~$\real$), so their restrictions to the
invariant subspace $\Lie U$  are also simultaneously diagonalizable
(cf.\ \cite[Thms.~26 and~27 in \S3.12, pp.~167--168]{ZariskiSamuel1}). Thus,
$\Lie U$ is a direct sum of weight spaces:
 $$ \Lie U = \bigoplus_{\psi \in \Psi} \Lie U_\psi .$$
 For each weight~$\psi$ of~$T$ on~$\Lie U$, we have
 $$ \Lie U_\psi
 = \Lie U \cap \rsp_\psi
 = \Lie U \cap \rsp^{= \psi} ,$$
 so
 $$ \Lie U_{\omega|_T} = \Lie U \cap \Lie N^{=\omega}$$
 and
 $$ \bigoplus_{\psi \neq \omega|_T} \Lie U_\psi
 = \Lie U \cap \Lie N^{\neq\omega} . $$
 The conclusion follows.
 \end{proof}

\section{Lower bound on the dimension of~$H$}  \label{dimHSect}
 
In this section, we prove  Corollary~\ref{tess->dim>1,2}, an
\emph{a~priori} lower bound on $\dim H$. On the way, we recall a result of
T.~Kobayashi that will also be used several times in later sections, and we
establish that crystallographic groups have only one end.

\subsection{T.~Kobayashi's Dimension Theorem}

The following theorem is essentially due to T.~Kobayashi. (Kobayashi assumed
that $H$ is reductive, but H.~Oh and D.~Witte \cite[Thm.~3.4]{OhWitte-CK}
pointed out that, by using Lemma~\ref{HcanbeAN}, this restriction can be
eliminated.) The proof here is based on Kobayashi's original argument and the
modifications of Oh-Witte, but uses less sophisticated topology. Namely,
instead of group cohomology and the spectral sequence of a covering space, we
use only some basic properties of homology groups of manifolds (including
Lemma~\ref{fiberbundle}). These comments also apply to
Theorem~\ref{construct-tess}.

\begin{thm}[{(Kobayashi, cf.\ \cite[Thm.~1.5]{Kobayashi-necessary},
\cite[Thm.~4.7]{Kobayashi-properaction})}] \label{noncpctdim}
 Let $H$ and~$H_1$ be closed, connected subgroups of~$G$, and assume there is
a crystallographis group $\Gamma$ for $G/H$, such that $\Gamma$ acts properly
discontinuously on $G/H_1$.
 Then:
 \begin{enumerate}
 \item \label{noncpctdim-notess}
 We have $d(H_1) \le d(H)$.
 \item \label{noncpctdim-tess}
 If $d(H_1) \ge d(H)$, then  $\Gamma \backslash G/H_1$ is compact, so $G/H_1$
has a tessellation.
 \end{enumerate}
 \end{thm}

\begin{proof}
 By Lemma~\ref{HcanbeAN}, we may assume $H, H_1 \subset AN$. (So $d(H)
= \dim H$ and $d(H_1) = \dim H_1$ \see{d(H)=dimH}.) 

  From Lemma~\ref{fiberbundle}, we know that $\Gamma \backslash G$ and $\Gamma
\backslash G/H_1$ have the same homology. Therefore
 $$ \max\{\, k \mid \homology_k(\Gamma \backslash G) \neq 0 \,\}
 = \max\{\, k \mid \homology_k(\Gamma \backslash G/H_1) \neq 0 \,\}
 \le \dim G/H_1 ,$$
 with equality if and only if $\Gamma \backslash G / H_1$ is compact
\cite[Cor.~8.3.4, p.~260]{Dold}.
 Similarly, we have
 $$ \max\{\, k \mid \homology_k(\Gamma \backslash G) \neq 0 \,\}
 = \dim G/H . $$
 Combining these two statements, we conclude \pref{noncpctdim-notess} that
$\dim G/H \le \dim G/H_1$ and, furthermore, \pref{noncpctdim-tess} that
equality holds if and only if $\Gamma \backslash G/H_1$ is compact.
 \end{proof}

\begin{cor}[(Kobayashi)] \label{noncpct-dim-notess}
 Let $H$ and~$H_1$ be closed, connected subgroups of~$G$, such that $d(H_1) >
d(H)$. If there is a compact subset~$C$ of~$A$, such that $\mu(H_1) \subset
\mu(H) C$, then $G/H$ does not have a tessellation.
 \end{cor}

\begin{proof}
 Suppose $\Gamma$ is a crystallographic group for $G/H$. (This will lead to a
contradiction.) Because $\Gamma$ acts properly discontinuously on $G/H$, the
assumption on $\mu(H_1)$ implies that $\Gamma$ also acts properly
discontinuously on $G/H_1$ \cf{proper<>mu(L)}. So
Theorem~\fullref{noncpctdim}{notess} yields a contradiction.
 \end{proof}

\subsection{Crystallographic groups have only one end}

It is easy to see that crystallographic groups are finitely generated; we now
show that they have only one end \see{tess->1end}.

\begin{defn}[{(cf.\ \cite[$0.2.A'_2$, p.~4]{Gromov-asymptotic})}]
\label{1endDefn}
 Let $F$ be a finite generating set for an (infinite) group~$\Gamma$.
We say that $\Gamma$ has \emph{only one end} if, for every partition $\Gamma
= A_1 \cup A_2 \cup C$ of~$\Gamma$ into three disjoint sets $A_1$,
$A_2$, and~$C$, such that $A_1$ and~$A_2$ are infinite, but $C$ is finite,
there exists $\gamma \in A_1$ and $f \in F \cup F^{-1}$, such that $\gamma f
\in A_2$. (This does not depend on the choice of the generating set~$F$.)
 \end{defn}

The following observation is a straightforward reformulation of
Definition~\ref{1endDefn} (obtained by letting $A_2 = \Gamma \smallsetminus
(A_1 \cup C')$ and $C = C' \smallsetminus A_1$). 

\begin{lem} \label{1endcomplement}
 Let $F$ be a finite generating set for an infinite group~$\Gamma$. If $A_1$
and~$C'$ are subsets of~$\Gamma$, such that
 \begin{itemize}
 \item $A_1$ is infinite,
 \item $C'$ is finite, and
 \item $A_1 f \in A_1 \cup C'$, for every $f \in F \cup F^{-1}$,
 \end{itemize}
 then the complement $\Gamma \smallsetminus A_1$ is finite.
 \end{lem}

\begin{rem}[{(cf.\ \cite[pp.~25--26, p.~32, and Prop.~2.14]{Cohen-coho1})}]
 Definition~\ref{1endDefn} is often stated in the language of Cayley graphs:
 The \emph{Cayley graph} of~$\Gamma$, with respect to the generating set~$F$,
is the graph $\operatorname{Cay}(\Gamma;F)$ whose vertex set~$V$ and edge set~$E$ are given
by:
 \begin{align*}
 V &= \Gamma; \\
 E &= \{\, (\gamma, \gamma f) \mid \gamma \in \Gamma, f \in F \cup F^{-1} \,\}
.
 \end{align*}
 The group~$\Gamma$ has only one end if and only if, for every finite
subset~$C$ of~$\Gamma$, the graph $\operatorname{Cay}(\Gamma; F) \smallsetminus C$ has only
one infinite component.
 \end{rem}

The following lemma is not difficult, but, unfortunately, we do not have a
proof that is both short and elementary.

\begin{lem}[{(see proof of Lemma~\fullref{dimT}{A})}] \label{HN=AN}
 If $HN=AN$, then, for some $x \in N$, the conjugate $x^{-1} H x$ is
normalized by~$A$.
 \end{lem}

\begin{cor} \label{codim>1}
 If $d(G) - d(H) \le 1$, and $G/H$ is not compact, then $G/H$ does not have a
tessellation.
 \end{cor}

\begin{proof}
 It suffices to show that $H$ is a Cartan-decomposition subgroup of~$G$
\see{CDS->notess}. 

 We may assume, without loss of generality, that $H \subset AN$
\see{HcanbeAN}; then
 $$\dim H + 1 = d(H) + 1 \ge d(G) = \dim(AN) $$
 \seeand{d(H)=dimH}{d(G)}.
 A theorem of B.~Kostant \cite[Thm.~5.1]{Kostant} implies that $N$ is a
Cartan-decomposition subgroup, so we may assume $N \not\subset H$; then
$\dim(H \cap N) \le \dim N - 1$. Therefore
 $$ \dim(HN) = \dim H + \dim N - \dim(H \cap N)
 \ge \dim H + 1
 \ge \dim(AN) .$$
 Hence $HN = AN$, so, from Lemma~\ref{HN=AN}, we see that, after replacing
$H$ by a conjugate subgroup, we may assume that $H$ is normalized by~$A$.
Then, letting $\omega = 0$ and $T = A$ in \pref{rootdecomp}, we see that
$\Lie H = (\Lie H \cap \Lie A) + (\Lie H \cap \Lie N)$. Since $HN = AN$, we
have $\Lie H + \Lie N = \Lie A + \Lie N$, so this implies that $\Lie A
\subset \Lie H$; therefore $H$ contains~$A$. Since $A$ is a
Cartan-decomposition subgroup \see{AisCDS}, this implies $H$ is a
Cartan-decomposition subgroup, as desired.
 \end{proof}

\begin{defn}[{(cf.\ \cite[$0.2.A'_2$, p.~4]{Gromov-asymptotic})}]
 A topological space~$M$ is \emph{connected at~$\infty$} if every compact
subset~$\mathcal{C}$ is contained in a compact subset~$\mathcal{C}'$, such
that the complement $M \smallsetminus \mathcal{C}'$ is connected.
 \end{defn}

\begin{prop} \label{tess->1end}
 If $\Gamma$ is a crystallographic group for $G/H$, then $\Gamma$
is finitely generated and has only one end.
 \end{prop}

\begin{proof} Assume, without loss of generality, that $H \subset AN$
\see{HcanbeAN}. Then $H$ is torsion free, so $\Gamma$ must act freely on
$G/H$; therefore $\Gamma \backslash G/H$ is a compact manifold (rather than an
orbifold).
 Because $\Gamma$ is essentially the fundamental group of $\Gamma \backslash
G/H$ (specifically, $\Gamma \iso \pi_1(\Gamma \backslash G/H)/\pi_1(G/H)$),
and the fundamental group of any compact manifold is finitely generated
\cite[Thm.~6.16, p.~95]{Raghunathan}, we know that $\Gamma$ is finitely
generated.

 From the Iwasawa decomposition $G = KAN$, we see that $G/H$ is  homeomorphic
to $K \times (AN/H)$, and Proposition~\ref{AN/H=Rd} asserts that $AN/H$ is
homeomorphic to~$\real^d$, for some~$d$. Obviously, we must have $d = \dim
(AN) - \dim H$, and we may assume $G/H$ is not compact (otherwise, $\Gamma$
is finite, so the desired conclusion is obvious), so Corollary~\ref{codim>1}
implies that $d > 1$. Thus, we conclude that $G/H$ is connected at~$\infty$.

To complete the proof, we use a standard argument (cf.\ \cite[$0.2.C_1$,
p.~5]{Gromov-asymptotic}) to show that, because $G/H$ is connected
at~$\infty$ and $\Gamma \backslash G/H$ is compact, the group~$\Gamma$
has only one end. To begin, note that there is a compact subset~$\mathcal{C}$
of~$G/H$, such that $\Gamma \mathcal{C} = G/H$. 
 Let 
 $$F_0 = \{\, f \in \Gamma \mid \mathcal{C} \cap f \mathcal{C} \neq \emptyset
\,\} $$
 (cf.\ \cite[(ii), p.~195]{PlatonovRapinchuk}).
 Because $\Gamma$ acts properly discontinuously on $G/H$, we know that $F_0$
is finite; let $F$ be a finite generating set for~$\Gamma$, such that $F_0
\subset F$.

Suppose $\Gamma = A_1 \cup A_2 \cup C$, with $|A_1| = |A_2| = \infty$
and $|C| < |\infty|$. (We wish to show there exist $\gamma \in A_1$ and $f
\in F$, such that $\gamma f \in A_2$; this establishes that $\Gamma$ has only
one end.)
 Because $G/H$ is connected at~$\infty$, there is a compact subset
$\mathcal{C}'$ of~$G/H$, containing $C \mathcal{C}$, such that $(G/H)
\smallsetminus \mathcal{C}'$ is connected. Because $C \mathcal{C} \subset
\mathcal{C}'$, we have
 $$ (G/H) \smallsetminus \mathcal{C}' = (\Gamma \mathcal{C}) \smallsetminus
\mathcal{C}' \subset A_1 \mathcal{C} \cup A_2 \mathcal{C} .$$
 Because $\Gamma$ acts properly discontinuously on $G/H$, we know $A_1
\mathcal{C}$ and $A_2 \mathcal{C}$ are closed (and neither is contained
in~$\mathcal{C}'$), so connectivity implies that 
 $A_1 \mathcal{C} \cap A_2 \mathcal{C} \neq \emptyset$:
 there exist $\gamma \in A_1$ and $\gamma' \in A_2$, such that 
 $\gamma \mathcal{C} \cap \gamma' \mathcal{C} \neq \emptyset$.
 Let $f = \gamma^{-1} \gamma'$; then
 $\gamma \in A_1$, $\gamma f = \gamma' \in A_2$, and
 $$ \mathcal{C} \cap f \mathcal{C} = \gamma^{-1} (\gamma \mathcal{C} \cap
\gamma' \mathcal{C}) \neq \emptyset ,$$
 so $f \in F_0 \subset F$, as desired.
 \end{proof}

\subsection{Walls of~$A^+$ and a lower bound on $d(H)$}

\begin{prop} \label{tess->misswall}
 Assume $\Rrank G = 2$. Let
 \begin{itemize}
 \item $L_1$ and $L_2$ be the two walls of~$A^+$, and
 \item $\Gamma$ be a crystallographic group for $G/H$.
 \end{itemize}
 If $H$ is not compact, then there exists $k \in \{1,2\}$, such that, for
every compact subset~$C$ of~$A$, the intersection $\mu(\Gamma) \cap L_k C$ is
finite.
 \end{prop}

\begin{center}
 \begin{figure}
 \includegraphics[scale=0.5]{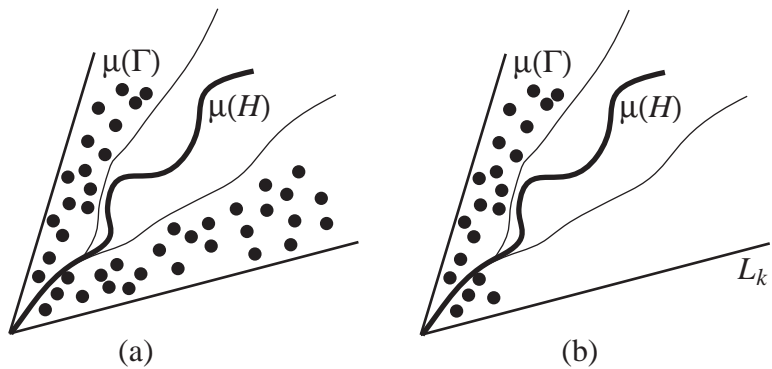}
 \caption{Proposition~\ref{tess->misswall}: (a)~$\mu(\Gamma)$ cannot be on
both sides of~$\mu(H)$, because $\Gamma$ has only one end. (b)~Therefore,
$\mu(\Gamma)$ stays away from $L_k$.}
 \label{misswallfig}
 \end{figure}
 \end{center}

\begin{proof}[Proof {\rm (cf.\ Figure~\ref{misswallfig})}]
 Suppose there is a compact subset~$C$ of~$A$, such that each of $\mu(\Gamma)
\cap L_1 C$ and $\mu(\Gamma) \cap L_2 C$ is infinite. (This will lead to a
contradiction.)
 Let $F$ be a (symmetric) finite generating set for~$\Gamma$ \see{tess->1end}.
We may assume $C$ is so large that $\mu(\gamma F) \subset \mu(\gamma) C$ for
every $\gamma \in \Gamma$ \see{bddchange}. We may also assume that $C$ is
convex and symmetric.

 Because $\Gamma$ acts properly on $G/H$, there is a compact
subset~$\mathcal{C}$ of~$A$, such that $\mu(H) C \cap \mu(\Gamma) \subset
\mathcal{C}$ \see{proper<>mu(L)}. Furthermore, we may assume that $\mu(L_1)C
\cap \mu(L_2) C \subset \mathcal{C}$.

Let 
 \begin{itemize}
 \item $M = \cup_{\gamma \in \Gamma} \mu(\gamma) C \smallsetminus \mathcal{C}$, 
 \item $M_1$ be the union of all the connected components of~$M$ that contain
a point of~$L_1$, and
 \item $A_1 = \Gamma \cap \mu^{-1}(M_1)$.
 \end{itemize}
 Then $A_1$ is infinite (because $\mu(\Gamma) \cap L_1 C$ is infinite). 
 Also, for any $\gamma\in A_1$ and $f \in F$, we have $\mu(\gamma f) \in
\mu(\gamma) C$, so $\gamma f \in A_1 \cup \mu^{-1}(\mathcal{C})$. Since
$\Gamma$ has only one end \see{tess->1end}, this implies $\Gamma
\smallsetminus A_1$ is finite \see{1endcomplement}. Because $\mu(\Gamma) \cap
L_2 C$ is infinite, we conclude that $M_1 \cap L_2 \neq \emptyset$.

Because $\mu(H)$ separates~$L_1$ from~$L_2$, and every connected component
of~$M_1$ contains a point of~$L_1$, we conclude that $\mu(H) \cap
M_1 \neq \emptyset$. This contradicts the fact that $\mu(H) C \cap \mu(\Gamma)
\subset \mathcal{C}$.
 \end{proof}

\begin{cor} \label{tess->missHk}
 Assume $\Rrank G = 2$.
 Let
 \begin{enumerate}
 \item $L_1$ and $L_2$ be the two walls of~$A^+$; and
 \item $H_1$ and $H_2$ be closed, connected, nontrivial 
subgroups of $G$,
 \end{enumerate}
 such that 
 $$\mu(H_k) \approx L_k$$
 for $k = 1,2$.
 If $H$ is not compact, then any crystallographic group for $G/H$ acts
properly discontinuously on either $G/H_1$ or $G/H_2$.
 \end{cor}

\begin{proof}
 Suppose $\Gamma$ acts properly discontinuously on \textbf{neither} $G/H_1$
\textbf{nor} $G/H_2$. (This will lead to a contradiction.) From
Proposition~\ref{proper<>mu(L)}, we know there is a compact subset~$C$
of~$A$, such that each of $\mu(\Gamma) \cap \mu(H_1) C$ and $\mu(\Gamma) \cap
\mu(H_2) C$ is infinite. Then, since $\mu(H_k) \approx L_k$, we may
assume (by enlarging~$C$) that each of $\mu(\Gamma) \cap L_1 C$ and
$\mu(\Gamma) \cap L_2 C$ is infinite. This contradicts the conclusion of
Proposition~\ref{tess->misswall}.
 \end{proof}

\begin{cor} \label{tess->dim>1,2}
 Assume $\Rrank G = 2$.
 Let
 \begin{enumerate}
 \item $L_1$ and $L_2$ be the two walls of~$A^+$; and
 \item $H_1$ and $H_2$ be closed, connected, nontrivial subgroups of $G$;
 \end{enumerate}
 such that 
 $$\mu(H_k) \approx L_k$$
 for $k = 1,2$.
 If $G/H$ has a tessellation, and $H$ is not compact, then 
 $$d(H) \ge \min\{ d(H_1), d(H_2) \} .$$
 \end{cor}

\begin{proof}
 The desired conclusion is obtained by combining Corollary~\ref{tess->missHk}
with Theorem~\fullref{noncpctdim}{notess}.
 \end{proof}

\begin{rem}
 For $G = \SL(3,\real)$, there does not exist a connected subgroup~$H_k$,
such that $\mu(H_k) \approx L_k$ \cf{SL3-B+}. Thus,
Corollary~\ref{tess->dim>1,2} does not provide a lower bound on $d(H)$ in
this case.
 \end{rem}

\section{One-dimensional subgroups} \label{1DSect}

Although the following conjecture does not seem to have been stated
previously in the literature, it is perhaps implicit in \cite{OhWitte-CK}.

\begin{conj} \label{d(H)=1->notess?}
 If $d(H) = 1$, then $G/H$ does not have a tessellation.
 \end{conj}

In this section, we establish that the conjecture is valid in two cases: if
either $\Rrank G \le 2$ \see{1Drank2} or $G$ is almost simple \see{1Dsimple}.
Each of these illustrates a general theorem: for groups of real rank two, the
conjecture follows from a theorem of Y.~Benoist and F.~Labourie that is based
on differential geometry; H.~Oh and D.~Witte observed that, for simple groups,
the conjecture follows from a theorem of G.~A.~Margulis that is based on
unitary representation theory.

The following example is the only case of Conjecture~\ref{d(H)=1->notess?}
that is needed in later sections. (It is used in the proof of
Theorem~\ref{SL3->notess}.) Because $\Rrank \bigl( \SL(3,\F) \bigr) = 2$ and
$\SL(3,\F)$ is almost simple, this example is covered both by the theorem of
Benoist-Labourie and by the theorem of Margulis, but it would be interesting
to have an easy proof.

\begin{prop}[{(see \ref{1Drank2} or \ref{1Dsimple})}] \label{wallinSL3}
 Assume $G = \SL(3,\F)$, for $\F = \real$, $\complex$, or~$\quaternion$, and
let
 \begin{equation} \label{SL3wallSubgrp}
 H_1 =
 \bigset{
 \begin{pmatrix}
 e^t & 0 & 0 \\
 0 & e^t & 0 \\
 0 & 0 & e^{-2t}
 \end{pmatrix}
 }{
 t \in \real
 }
 \subset G
 .
 \end{equation}
 Then $G/H_1$ does not have a tessellation.
 \end{prop}

Let us begin with an easy observation.

\begin{lem} \label{1Drank1}
 If $d(H) = 1$ and $\Rrank G = 1$, then $G/H$ does not have a tessellation.
 \end{lem}

\begin{proof}
 We may assume $H \subset AN$ \see{HcanbeAN}. From \pref{Rrank1-CDS}, we know
that $H$ is a Cartan-decomposition subgroup, so Lemma~\ref{CDS->notess}
implies that $G/H$ must be compact; thus, the trivial group~$e$ is a
crystallographic group for $G/H$. However, since
 $$d(G) = \dim A + \dim N \ge 1 + 1 > 1 = d(H) $$
 \see{d(G)}, and $e$~acts properly discontinuously on $G/G$, this contradicts
Theorem~\fullref{noncpctdim}{notess}.
 \end{proof}

\begin{rem}
 The proof of Lemma~\ref{1Drank1} shows that the dimension of every connected,
cocompact subgroup of~$G$ is at least $d(G)$. This is a result of M.~Goto and
H.--C.~Wang \cite[(1.2), p.~263]{GotoWang}.
 \end{rem}

\subsection{The geometric method of Y.~Benoist and F.~Labourie}

\begin{thm}[{(Benoist-Labourie \cite[Cor.~3]{BenoistLabourie})}]
\label{center->notess}
 If $H$ is reductive and contains an element of~$A$ in its center, then $G/H$
does not have a tessellation.
 \end{thm}

To illustrate the idea behind Theorem~\ref{center->notess}, we give a direct
proof of the following special case (under an additional technical assumption
\see{1DA->notessPf-H=Z}), which is sufficient for our needs. (Note that
Condition~\pref{1DA->notessPf-H=Z} is satisfied for the subgroup~$H_1$ of
Proposition~\ref{wallinSL3}.) Benoist and Labourie prove the general case by
using a slightly different 1-form in place of the form~$\omega$ that we
define in Step~\ref{1DA->notessPf-Omega}.

\begin{cor}[{(Benoist-Labourie)}] \label{1DA->notess}
 If $H$ is a one-dimensional subgroup of~$A$, then $G/H$ does not have a
tessellation.
 \end{cor}

\begin{proof}[{Proof \rm \cite{BenoistLabourie}}]
 Suppose $\Gamma$ is a crystallographic group for $G/H$. (This will lead to a
contradiction.) By passing to finite-index subgroups, we may assume that $G$
and $H$ are connected.

We will construct a volume form~$\nu$ on $\Gamma \backslash G/H$ that is
exact: $\nu = d \xi$. The integral of~$\nu$ over $\Gamma \backslash G/H$
is the volume of $\Gamma \backslash G/H$, which is obviously not~$0$, but
Stokes' Theorem implies that the integral of any exact form over a closed
manifold is~$0$. This is a contradiction.

\setcounter{step}{0}

\begin{step} \label{1DA->notessPf-Omega}
 Construction of a certain $2$-form~$\breve{\Omega}$ on $\Gamma \backslash
G/H$.
 \end{step}
 Let $\Lie M = \Lie H^\perp$ be the orthogonal complement to~$\Lie H$, under
the Killing form. Then $\Lie M$ is an $(\Ad_G H)$-invariant subspace of~$\Lie
G$, such that $\Lie G = \Lie H + \Lie M$ and $\Lie H \cap \Lie M = 0$. Let
$\omega_e \colon \Lie G \to \Lie H$ be the projection with kernel~$\Lie M$,
and let $\omega$ be the corresponding left-invariant $\Lie H$-valued 1-form
on~$G$.

The space $G$ is a homogeneous principal $H$-bundle over $G/H$. It is well
known \cite[Thm.~II.11.1, p.~103]{KobayashiNomizu1} that $\omega$ is
the connection form of a $G$-invariant connection on this bundle, and that
the curvature form~$\Omega$ of this connection is given by
 \begin{equation} \label{Omega(M,M)}
  \mbox{$\Omega(X, Y) = - \frac{1}{2} \omega \bigl( [X,Y] \bigr)$
 \qquad
 for $X,Y \in \Lie M$.}
 \end{equation}
 Also, because $H$ is abelian, the Structure Equation \cite[Thm.~II.5.2,
p.~77]{KobayashiNomizu1} implies
 \begin{equation} \label{structeq}
  \mbox{$\Omega(X, Y) = d \omega(X,Y)$
 \qquad
 for $X,Y \in \Lie G$.}
 \end{equation}
 By identifying the $1$-dimensional Lie algebra~$\Lie H$ with~$\real$,
we may think of $\omega$ and~$\Omega$ as ordinary (that is, $\real$-valued)
differential forms.

 Because $\omega$ and (hence)~$\Omega$ are left-invariant, they determine
well-defined forms $\overline{\omega}$ and~$\overline{\Omega}$ on $\Gamma
\backslash G$. (Note that $\overline{\omega}$ is a connection form on the
principal $H$-bundle $\Gamma \backslash G$ over $\Gamma \backslash G/H$, and
the curvature form of this connection is~$\overline{\Omega}$.) Because $H$ is
abelian, we have $\overline{\Omega}_{gh} = \overline{\Omega}_g$ for all $h
\in H$ (cf.\ \cite[Prop.~II.5.1(c), p.~76]{KobayashiNomizu1}), so the
horizontal form $\overline{\Omega}$ determines a well-defined form
$\breve{\Omega}$ on the base space $\Gamma \backslash G/H$.

\begin{step}
 Construction, for a certain~$s$, of a certain $s$-form~$\breve{\mu}$ on
$\Gamma \backslash G/H$.
 \end{step}
 Identifying $\Lie H$ with~$\real$ provides an ordering on the weights
of~$\Lie H$. Let 
 \begin{itemize}
 \item $\Lie M_0 = \Lie C_{\Lie G}(\Lie H) \cap \Lie M$ be the $0$~weight
space of $\ad_{\Lie G} \Lie H$ on~$\Lie M$;
 \item $\Lie M_+$ be the sum of the positive weight spaces of $\ad_{\Lie G}
\Lie H$ on~$\Lie M$;
 \item $\Lie M_-$ be the sum of the negative weight spaces of $\ad_{\Lie G}
\Lie H$ on~$\Lie M$;
 \item $s = \dim \Lie M_0$; and
 \item $\mu$ be a nontrivial left-invariant $s$-form on~$G$, such that
 \begin{equation} \label{mu(H+M)}
 \mu( \Lie H \oplus \Lie M_+ \oplus \Lie M_-, *, \ldots, *) = 0 . 
 \end{equation}
 \end{itemize}
 Because $\mu$ is left-invariant, and 
 $\mu( \Lie H, *, \ldots, *) = 0$,
 the form~$\mu$ determines a well-defined differential form~$\breve{\mu}$ on
$\Gamma \backslash G/H$.
 (We remark that, because $\Lie G = \Lie H \oplus \Lie M_+ \oplus \Lie M_0
\oplus \Lie M_-$, Condition~\pref{mu(H+M)} implies that the form~$\mu$ is
unique up to a scalar multiple.) 

\begin{step} \label{1DA->notessPf-volume}
 For a certain~$r$, the wedge product $\breve{\mu} \wedge
\breve{\Omega}^{\wedge r}$ is a volume form on $\Gamma \backslash G/H$.
 \end{step}
 Let $r = \dim \Lie M_+$. It suffices to show that the restriction of
$\Omega$ to $\Lie M_+ \oplus \Lie M_-$ is a symplectic form. It is obviously
skew, so we need only show that it is nondegenerate. Thus, letting
 \begin{align*}
 \Lie Z
 &= \{\, X \in \Lie M_+ \oplus \Lie M_- \mid \Omega (X, \Lie M_+ \oplus \Lie M_-) = 0 \bigr) \\
 &= \{\, X \in \Lie M_+ \oplus \Lie M_- \mid \omega \bigl( [X, \Lie M_+ \oplus \Lie M_-] \bigr)
= 0 \,\}
 & \mbox{\pref{Omega(M,M)}} \\
 &= \{\, X \in \Lie M_+ \oplus \Lie M_- \mid
 [X, \Lie M_+ \oplus \Lie M_-] \subset \Lie M \,\}
 & \mbox{(definition of~$\omega$)}
\\
 &= \{\, X \in \Lie M_+ \oplus \Lie M_- \mid
 \langle \Lie H \mid [X, \Lie M_+ \oplus \Lie M_-]
\rangle_{\text{Killing}} = 0 \,\}
 & \mbox{($\Lie M = \Lie H^\perp$)}
 ,
 \end{align*}
 we wish to show $\Lie Z = 0$.

There is no harm in passing to the complexification~$\Lie G^{\complex}$
of~$\Lie G$. Let $\Lie T^{\complex}$ be a Cartan subalgebra of~$\Lie
G^{\complex}$ that contains~$\Lie H$. Because $\Lie T^{\complex}$ preserves
the Killing form, centralizes~$\Lie H$, and normalizes $\Lie M_+ \oplus \Lie
M_-$, we know that $\Lie Z^{\complex}$ is $(\ad_{\Lie G^{\complex}} \Lie
T^{\complex})$-invariant; thus, $\Lie Z^{\complex}$ is a sum of root spaces. 

Suppose there exists a nonzero element~$X$ of~$\Lie Z^{\complex}$, such that
$X$ belongs to some root space~$\Lie G^{\complex}_{\alpha}$. (This will
lead to a contradiction.) There exists $Y \in \Lie G^{\complex}_{-\alpha}$,
such that
 $$ \langle t \mid [X,Y] \rangle_{\text{Killing}} = \alpha(t) $$
 for all $t \in \Lie T^{\complex}$ \cite[Prop.~8.3(c), p.~37]{Humphreys-Lie}.
 Because $\Lie G^{\complex}_{\alpha} \subset \Lie M_+ \oplus \Lie M_-$, we
have $\alpha(\Lie H) \neq 0$. Then $(-\alpha)(\Lie H) \neq 0$, so 
 $\Lie G^{\complex}_{-\alpha} \subset \Lie M_+ \oplus \Lie M_-$. We now know
that $X \in \Lie Z^{\complex}$ and $Y \in \Lie M_+ \oplus \Lie M_-$, so
 $\langle \Lie H \mid [X, Y] \rangle_{\text{Killing}} = 0$.
 We therefore conclude, from the definition of~$Y$, that $\alpha(\Lie H) =
0$. This is a contradiction.

\begin{step}
 The form $\breve{\Omega}$ is exact: we may write $\breve{\Omega} = d \breve{\phi}$.
 \end{step}
 Let $\breve{\omega}_0$ be the connection form of a flat connection on $\Gamma
\backslash G$ over $\Gamma \backslash G/H$. (Since the principal bundle is
trivial \see{trivialbundle}, it is obvious that there is a flat connection.)
For any vector field~$X$ on $\Gamma \backslash G/H$, let $\widetilde X$ be the
lift of~$X$ to a vector field on $\Gamma \backslash G$ that is horizontal
with respect to the flat connection~$\breve{\omega}_0$. 

Since $H$ is abelian, there is a well-defined 1-form~$\breve{\phi}$ on $\Gamma
\backslash G/H$ given by
 $$ \breve{\phi}(X) = \overline{\omega}(\widetilde X) .$$
 Then
 \begin{align*}
 d \breve{\phi}(X,Y)
 &= \frac{1}{2} \left(
 X \bigl( \breve{\phi}(Y) \bigr) - Y \bigl(
\breve{\phi}(X) \bigr) -
 \breve{\phi} \bigl( [X,Y] \bigr) \right)
 & \mbox{(definition of~$d$)} \\
 &= \frac{1}{2} \left( \widetilde X \bigl( \overline{\omega}(\widetilde Y)
\bigr) - \widetilde Y \bigl( \overline{\omega}(\widetilde X) \bigr) -
 \overline{\omega} \bigl( \widetilde{[X,Y]} \bigr) \right)
 & \mbox{(definition of~$\breve{\phi}$)} \\
 &= \frac{1}{2} \left( \widetilde X \bigl( \overline{\omega}(\widetilde Y)
\bigr) - \widetilde Y \bigl( \overline{\omega}(\widetilde X) \bigr) -
 \overline{\omega} \bigl( [\widetilde X,\widetilde Y] \bigr) \right)
 & \mbox{($\breve{\omega}_0$ is flat)} \\
 &= d \overline{\omega}(\widetilde X, \widetilde Y)
 & \mbox{(definition of~$d$)} \\
 &= \overline{\Omega}(\widetilde X, \widetilde Y)
 & \mbox{\pref{structeq}} \\
 &= \breve{\Omega}(X,Y)
 & \mbox{(definition of~$\breve{\Omega}$)} .
 \end{align*}

\begin{pfassump} 
 For simplicity, assume that 
 \begin{equation} \label{1DA->notessPf-H=Z}
 \mbox{every hyperbolic element of the center of ${\Lie C}_{\Lie G}(\Lie H)$
is contained in~$\Lie H$.}
 \end{equation}
 \end{pfassump}%

\begin{step} \label{1DA->notessPf-closed}
 We have $d\breve{\mu} \wedge \breve{\Omega}^{\wedge(r-1)} = 0$. 
 \end{step}
 Let
 \begin{itemize}
 \item $X_1,\ldots,X_r$ be a basis of~$\Lie M_+$, and
 \item $Y_1,\ldots,Y_r$ be the dual basis of~$\Lie M_-$, with respect to the
symplectic form~$\Omega$ on $\Lie M_+ \oplus \Lie M_-$.
 \end{itemize}
 Thus, $\Omega_e(X_j,Y_k) = \delta_{j,k}$.
 
Let $Z_1,\ldots,Z_s$ be a basis of~$\Lie M_0$,
 write 
 $$[X_j,Y_k] = \sum_\ell a_{j,k}^\ell Z_\ell \pmod {\Lie H + \Lie
M_+ + \Lie M_-} ,$$
 and define
 $$ W = \sum_j [X_j,Y_j] = \sum_{j,\ell} a_{j,j}^\ell Z_\ell .$$

\begin{substep} \label{1DA->notessPf-closed-Windep}
 $W$ is independent of the choice of the basis $\{X_j\}$ of $\Lie M_+$
{\upshape(}with the understanding that $\{Y_j\}$ must be the dual basis of
$\Lie M_-${\upshape)}.
 \end{substep}
 Let
 \begin{equation} \label{elemop}
 X'_j = 
 \begin{cases}
 a X_1 + b X_2 & \mbox{if $j = 1$} \\
 X_j & \mbox{if $j \ge 2$}
 \end{cases} 
 \end{equation}
 for some $a,b \in \real$ with $a \neq 0$. Then
 $$Y'_j = 
 \begin{cases}
 (1/a) Y_1 & \mbox{if $j = 1$} \\
 Y_2 - (b/a) Y_1 & \mbox{if $j = 2$} \\
 Y_j & \mbox{if $j \ge 3$}
 ,
 \end{cases}
 $$
 so 
 \begin{align*}
 W'
 &= [X'_1,Y'_1] + [X'_2,Y'_2] + \sum_{j \ge 3} [X'_j,Y'_j] \\
 &= \bigl( [X_1,Y_1] + (b/a)[X_2,Y_1] \bigr)
  + \bigl( [X_2,Y_2] - (b/a) [X_2,Y_1] \bigr) + \sum_{j \ge 3} [X_j,Y_j] \\
 &=  \sum_j [X_j,Y_j] \\
 &= W .
 \end{align*}
 Since $X_1,\ldots,X_r$ can be transformed into any other basis of~$\Lie M_+$
by a sequence of elementary operations as in \pref{elemop}, we conclude that
$W$ is independent of the choice of basis, as desired.

\begin{substep}
 We have $W \in \Lie H$.
 \end{substep}
 Substep~\ref{1DA->notessPf-closed-Windep} implies that $W$ is centralized by
$C_G(H)$, so $W$ is in the center of $\Lie C_{\Lie G}(\Lie H)$. 

 Let $\sigma$ be a Cartan involution of~$G$ with $\sigma(h) = -h$ for $h \in
\Lie H$. Substep~\ref{1DA->notessPf-closed-Windep} implies $\sigma(W) = -W$.
(From Substep~\ref{1DA->notessPf-closed-Windep}, we see that $W$ depends only
on~$\Lie H$ and the chosen identification of~$\Lie H$ with~$\real$; $\sigma$
reverses the choice of identification.) 

Thus, $W$ is a hyperbolic element of the center of $\Lie C_{\Lie G}(\Lie H)$.
By Assumption~\ref{1DA->notessPf-H=Z}, this implies $W \in \Lie H$, as
desired.

\begin{substep}  
 Completion of Step~\ref{1DA->notessPf-closed}.
 \end{substep}
 Let
 $$ X_1^*, \ldots, X_r^*, Y_1^*, \ldots, Y_r^*, Z_1^*, \ldots, Z_s^* $$
 be the basis of $(\Lie G/\Lie H)^*$ dual to
 $$ X_1, \ldots, X_r, Y_1, \ldots, Y_r, Z_1, \ldots, Z_s . $$
 We may assume $\mu = Z_1^* \wedge \cdots \wedge Z_s^*$. Then, because
 $[\Lie M_0 + \Lie M_-, \Lie M_-] \subset \Lie M_-$ and 
 $[\Lie M_0 + \Lie M_+, \Lie M_+] \subset \Lie M_+$, we have
 $$ dZ_\ell^*( \Lie M_0 + \Lie M_-, \Lie M_- ) = 0 
 = dZ_\ell^*( \Lie M_0 + \Lie M_+, \Lie M_+ ) ,$$
 so 
 $$ dZ_\ell^* = - \sum_{j,k} a_{j,k}^\ell X_j^* \wedge Y_k^* \pmod{\Lie M_0^*
\wedge \Lie M_0^*} .$$
 Therefore
 \begin{align*}
 (s+1) \, d \breve{\mu}
 &= \sum_{\ell} (-1)^{\ell-1} dZ_\ell^* \wedge Z_1^* \wedge \cdots \wedge
\widehat{Z_\ell^*} \wedge \cdots \wedge  Z_s^* \\
 &= \sum_{j,k,\ell} (-1)^\ell a_{j,k}^\ell
 X_j^* \wedge Y_k^* \wedge Z_1^* \wedge \cdots \wedge \widehat{Z_\ell^*}
\wedge \cdots \wedge  Z_s^* .
 \end{align*}
 From the choice of $Y_1, \ldots, Y_r$, we have
 $\Omega = 2 \sum_{j=1}^r X_j^* \wedge Y_j^*$,
 so
 $$ \Omega^{\wedge(r-1)} 
 = 2^{r-1} (r-1)! \, \sum_{j=1}^r X_1^* \wedge Y_1^* \wedge \cdots \wedge
\widehat{X_j^*} \wedge \widehat{Y_j^*} \wedge \cdots \wedge X_r^* \wedge Y_r^*
 $$
 and
 $$\Omega^{\wedge r} 
 = 2^{r} r! \, X_1^* \wedge Y_1^* \wedge \cdots \wedge X_r^* \wedge Y_r^*
 .$$
 Hence
 \begin{align*}
 (s+1) \, d \breve{\mu} \wedge \breve{\Omega}^{\wedge(r-1)}
 &= \frac{1}{2r} \, \sum_{j,\ell} 
 (-1)^{\ell} a_{j,j}^\ell
 Z_1^* \wedge \cdots \wedge \widehat{Z_\ell^*} \wedge \cdots \wedge  Z_s^*
 \wedge \breve{\Omega}^{\wedge r} \\
 &= - \frac{s+2r}{2r} \, \iota_{W} (\breve{\mu} \wedge
\breve{\Omega}^{\wedge r}) .
 \end{align*}
 Since $W \in \Lie H$, we have 
 $ \iota_{W} (\mu_e \wedge \Omega_e^{\wedge r}) = 0$, so the desired
conclusion follows.

\begin{step} \label{1DA->notessPf-volexact}
 $\breve{\mu} \wedge \breve \Omega^{\wedge r}$ is exact.
 \end{step}
 We have
 \begin{align*}
  d \bigl( \breve{\mu} \wedge \breve{\phi} \wedge \breve{\Omega}^{\wedge(r-1)} \bigr)
 &= \pm d \breve{\mu} \wedge \breve{\Omega}^{\wedge(r-1)} \wedge \breve{\phi}
 \pm \breve{\mu} \wedge d \breve{\phi}  \wedge \breve{\Omega}^{\wedge(r-1)}
 \pm \breve{\mu} \wedge \breve{\phi} \wedge d \breve{\Omega}^{\wedge(r-1)} \\
 &= 0 \wedge \breve{\phi}
 \pm \breve{\mu} \wedge \breve{\Omega}  \wedge \breve{\Omega}^{\wedge(r-1)}
 \pm \breve{\mu} \wedge \breve{\phi} \wedge 0 \\
 &= \pm \breve{\mu} \wedge \breve{\Omega}^{\wedge r} .
 \end{align*}

\begin{step} 
 A contradiction.
 \end{step}
 From Step~\ref{1DA->notessPf-volume}, we know that
 $$ \int_{\Gamma \backslash G/H} \breve{\mu} \wedge \breve{\Omega}^{\wedge r}
\neq 0 .$$
 On the other hand, Step~\ref{1DA->notessPf-volexact} implies that this
integral is zero.
 This is a contradiction.
 \end{proof}

\begin{notation}[{\cite[p.~320]{Benoist}}] \label{oppinv}
 Let $\tau$ be the opposition involution in~$A^+$; that is, for $a\in A^+$,
$\tau(a) = \mu(a^{-1})$ is the unique element of~$A^+$ that is conjugate
(under an element of the Weyl group) to $a^{-1}$. Thus, for all $h \in G$, we
have
 $$ \mu(h^{-1}) = \tau \bigl( \mu(h) \bigr) .$$

See \pref{SL3-oppinv} for an explicit description of the opposition
involution in $G = \SL(3,\real)$. For some groups, such as $G = \SO(2,n)$, we
have $\mu(h^{-1}) = \mu(h)$ for all $h \in G$ \see{mu(h-1)}; in such a case,
the opposition involution is simply the identity map on~$A^+$.
 \end{notation}

\begin{cor} \label{1Drank2}
 If $d(H) = 1$ and $\Rrank G \le 2$, then $G/H$ does not have a tessellation.
 \end{cor}

\begin{proof}
 Suppose $\Gamma$ is a crystallographic group for $G/H$. (This will lead to a
contradiction.)

From~\pref{1Drank1}, we know $\Rrank G = 2$. Let $L_1$ and $L_2$ be the two
walls of~$A^+$ and, for $k \in \{1,2\}$, let $H_k = L_k \cup L_k^{-1}$.
Because $L_k$ is a ray (that is, a one-parameter semigroup), it is clear that
$H_k$ is a subgroup of~$A$. 

From Proposition~\ref{tess->misswall}, we know that there is some $k \in
\{1,2\}$, such that 
 $$ \mbox{$\mu(\Gamma) \cap L_k C$ is finite,} $$
 for every compact subset~$C$ of~$A$. Since $\Gamma = \Gamma^{-1}$, we have
$\tau \bigl( \mu(\Gamma) \bigr) = \mu(\Gamma)$, so this implies that
 $$ \mbox{$\mu(\Gamma) \cap \tau(L_k) C$ is finite,} $$
 for every compact subset~$C$ of~$A$. Also, because $L_k \subset A^+$, we
have $\mu(L_k) = L_k$, so
 $$\mu(H_k)
 = \mu( L_k \cup L_k^{-1})
 = \mu(L_k) \cup \tau \bigl( \mu(L_k) \bigr) 
 = L_k \cup \tau(L_k) .$$
 Therefore
 $$ \mbox{$
 \mu(\Gamma) \cap \mu(H_k) C
 = \bigl( \mu(\Gamma) \cap L_k C \bigl) \cup \bigl( \mu(\Gamma) \cap
\tau(L_k) C \bigl)
 $ is finite,} $$
 for every compact subset~$C$ of~$A$. Hence, Corollary~\ref{proper<>mu(L)}
implies that $\Gamma$ acts properly discontinuously on $G/H_k$. Then, because
$d(H) = 1 = d(H_k)$ \see{d(H)=dimH}, Theorem~\fullref{noncpctdim}{tess}
implies that $G/H_k$ has a tessellation. This contradicts
Corollary~\ref{1DA->notess}.
 \end{proof}

\subsection{The representation-theoretic method of G.~A.~Margulis}

\begin{defn}[{\cite[Defn.~2.2, Rmk.~2.2]{Margulis-CK}}]
\label{tempered-def}
 The subgroup~$H$ is \emph{tempered} in~$G$ if there exists a
(positive) function $f \in L^1(H)$ (with respect to a left-invariant Haar
measure on~$H$), such that, for every unitary representation~$\pi$ of~$G$,
either
 \begin{itemize}
 \item $|\langle \pi (h) \phi \mid \psi \rangle  | \leq f(h) \,
\lVert\phi\rVert \lVert\psi\rVert$ for all $h \in H$ and all $K$-fixed
vectors $\phi$ and~$\psi$; or
 \item some nonzero vector is fixed by every element of~$\pi(G)$.
 \end{itemize}
 \end{defn}

For many examples of tempered subgroups of simple Lie groups, see
\cite{Oh-tempered}.

\begin{thm}[{(Margulis \cite[Thm.~3.1]{Margulis-CK})}]
\label{MargulisTempered}
 If $H$ is noncompact and tempered, then $G/H$ does not have a tessellation.
 \end{thm}

\begin{proof}[{Proof \cite{Margulis-CK}}]
 Suppose $\Gamma$ is a crystallographic group for $G/H$. (This will lead to a
contradiction.) To simplify the notation somewhat (and because this is the
only case we need), let us assume $H = \{h^t\}$ is a one-parameter subgroup
of~$G$.

Because $G$ and~$\Gamma$ are unimodular (recall that $G$ is semisimple
\see{standing} and $\Gamma$ is discrete), there is a $G$-invariant measure
(in fact, a $G$-invariant volume form) on the homogeneous space $\Gamma
\backslash G$ \cite[Lem~1.4, p.~18]{Raghunathan}. Thus, the natural
representation~$\pi$ of~$G$ on $L^2(\Gamma \backslash G)$, defined by
 $$ \mbox{
 $\bigl( \pi(g) \phi \bigr)(x) = \phi(x g^{-1})$,
 \qquad
 for $\phi \in L^2(\Gamma \backslash G)$, $g \in G$, $x \in \Gamma
\backslash G$,
 } $$
 is unitary.

Because $H$ is noncompact, and acts properly on $\Gamma \backslash G$, we
know that any compact subset of $\Gamma \backslash G$ has infinitely many
pairwise-disjoint translates (all of the same measure), so we see that
 \begin{equation} \label{MargulisTemperedPf-notcpct}
 \mbox{$\Gamma \backslash G$ is not compact, and has infinite volume} .
 \end{equation}
 Therefore, $\pi$ has no nonzero $G$-invariant vectors, so, because $H$ is
tempered, we know that there is some 
 $f \in L^1(\real)$,
 such that
 \begin{equation} \label{MargulisTemperedPf-<f}
 f(t) \, \lVert\phi\rVert_2 \lVert\psi\rVert_2
 \ge \left| \int_{\Gamma \backslash G} \phi( x h^t) \, \psi(x) \, dx \right|
 \end{equation}
 for all $t \in \real$ and all $K$-invariant $\phi,\psi \in L^2(\Gamma
\backslash G)$.

Because $\Gamma \backslash G/H$ is compact, there is a compact subset~$C$
of~$G$, such that 
 $\Gamma C H = G$;
 let $\overline C$ be the image of~$C$ in $\Gamma \backslash G$.
 From the choice of~$C$, we know, for each $x \in \Gamma \backslash G$,
that there is some $T_x \in \real$, such that
 \begin{equation} \label{MargulisTemperedPf-Tx}
 x h^{T_x} \in \overline{C} .
 \end{equation}
 Because 
 $\bigcup_{t=0}^1 \overline{C} h^t K$
 is a compact subset of $\Gamma \backslash G$, there is a
positive, continuous function~$\phi$ on $\Gamma \backslash G$ with compact
support, such that
 \begin{equation} \label{MargulisTemperedPf-phi(xht)}
 \mbox{$\phi(x h^t) \ge 1$ for all $x \in \overline{C}$ and all $t \in
[0,1]$,} 
 \end{equation}
 and, by averaging over~$K$, we may assume that $\phi$ is $K$-invariant.

Fix some large $T \in \real^+$. Because $\bigcup_{t = -(T+1)}^{T+1}
\overline{C} h^t K$ is compact, and $\Gamma \backslash G$ has infinite volume
\see{MargulisTemperedPf-notcpct}, there is some $K$-invariant
continuous function $\psi_T$ on $\Gamma \backslash G$, such that
$\lVert\psi_T\rVert_2 = 1$, 
 \begin{equation} \label{MargulisTemperedPf-psi<1}
  \mbox{$0 \le \psi_T(x) \le 1$ for all $x \in \Gamma \backslash G$} ,
 \end{equation}
 and
 \begin{equation} \label{MargulisTemperedPf-Tx>T+1}
 \mbox{$|T_x| > T + 1$ for all $x$ in the support of~$\psi_T$} .
 \end{equation}
 We have
 \begin{align*}
 \lVert\phi\rVert_2 \int_{|t|>T} f(t) \, dt
 &\ge \int_{\Gamma \backslash G} \int_{|t|>T} \phi( x h^t) \, \psi_T(x) \, dt
\, dx 
 & \mbox{(\pref{MargulisTemperedPf-<f} and $\lVert\psi_T\rVert_2 = 1$)} \\
 &\ge \int_{\Gamma \backslash G} \int_0^1 \phi( x h^{T_x + t}) \, \psi_T(x) \,
dt \, dx 
 & \mbox{\pref{MargulisTemperedPf-Tx>T+1}} \\
 &\ge \int_{\Gamma \backslash G} \psi_T(x) \,
dx 
 & \mbox{(\pref{MargulisTemperedPf-Tx} and
\pref{MargulisTemperedPf-phi(xht)})} \\
 &\ge 1
 & \mbox{(\pref{MargulisTemperedPf-psi<1} and $\lVert\psi_T\rVert_2 = 1$)}
 .
 \end{align*}
 However, because $f \in L^1(\real)$, we know that 
 $\lim_{T \to \infty} \int_{|t|>T} f(t) \, dt = 0$.
 This is a contradiction.
 \end{proof}

We state the following well-known result of representation theory without
proof. As is explained in \cite[\S3, p.~140]{KatokSpatzier}, it can be
obtained by combining work of R.~Howe \cite[Cor.~7.2 and \S7]{Howe} and
M.~Cowling \cite[Thm.~2.4.2]{Cowling}. (The assumption that $\Rrank G \ge 2$
can be relaxed: it suffices to assume that $G$ is not locally isomorphic to
$\SO(1,n)$ or $\SU(1,n)$.)

Fix any matrix norm $\lVert \cdot \rVert$ on~$G$; for example, we may let
 $\lVert g \rVert = \max_{j,k} |g_{j,k}|$.

\begin{thm}[{(Cowling, Howe)}] \label{expdecay}
 If $G$ is almost simple, and $\Rrank G \ge 2$, then there are constants $C>0$
and $p>0$ such that, for every unitary representation~$\pi$ of~$G$, either
 \begin{enumerate}
 \item $ | \langle \pi(g) \phi \mid \psi \rangle|
 \leq C \lVert \phi \rVert \lVert \psi \rVert \lVert g\rVert^{-p}$
 for all $g \in G$ and all $\pi(K)$-fixed vectors $\phi$ and~$\psi$; or
 \item some nonzero vector is fixed by every element of~$\pi(G)$.
 \end{enumerate}
 \end{thm}

Although we cannot prove Theorem~\ref{expdecay} here, we present an
elementary proof of the following related result, which, unfortunately, is
qualitative, rather than quantitative.
 On the other hand, this simple result applies to all vectors, not only the
$K$-fixed vectors, and it applies to all semisimple groups, including
$\SO(1,n)$ and $\SU(1,n)$.
 It was first proved by R.~Howe and C.~Moore \cite[Thm.~5.1]{HoweMoore}
and (independently) R.~Zimmer \cite[Thm.~5.2]{Zimmer-orbitspace}.

\begin{thm}[{\cite[Thm.~2.2.20, p.~23]{ZimmerBook}}] \label{matcoeffs->0}
 If
 \begin{itemize}
 \item $G$ is connected and almost simple;
 \item $\pi$ is a unitary representation of~$G$ on a Hilbert
space~$\hilbert$, such that no nonzero vector is fixed by $\pi(G)$; and
 \item $\{g_j\}$ is a sequence of elements of~$G$, such that $\lVert g_j
\rVert \to \infty$,
 \end{itemize}
 then $\langle \pi(g_j) \phi \mid \psi \rangle \to 0$, for every $\phi,\psi
\in \hilbert$.
 \end{thm}

\begin{proof}[{Proof \rm (Ellis-Nerurkar \cite{EllisNerurkar})}]

 \setcounter{case}{0}

\begin{case} \label{matcoeffs->0-A}
 Assume $\{g_j\} \subset A$.
 \end{case} 
 By passing to a subsequence, we may assume $\pi(g_j)$ converges weakly, to
some operator~$E$; that is,
 $$ \langle \pi(g_j) \phi \mid \psi \rangle
 \to \langle E \phi \mid \psi \rangle
 \mbox{ \ for every $\phi,\psi \in \hilbert$} .$$
 Let 
 \begin{equation} \label{horodefn}
 U = \{\, v \in G \mid g_j^{-1} v g_j \to e \,\} 
 \mbox{ \qquad and \qquad}
 U^- = \{\, u \in G \mid g_j u g_j^{-1} \to e \,\} .
 \end{equation}
 For $u \in U^-$, we have
 $$ \langle E\pi(u) \phi \mid \psi \rangle
 = \lim \langle \pi(g_j u) \phi \mid \psi \rangle
 = \lim \langle \pi(g_j u g_j^{-1}) \pi(g_j) \phi \mid \psi \rangle
 = \lim \langle \pi(g_j) \phi \mid \psi \rangle
 = \langle E \phi \mid \psi \rangle ,$$
 so $E \pi(u) = E$. Therefore $E \bigl( (\hilbert^{U^-})^\perp \bigr) = 0$.

We have
 $$ \langle E^* \phi \mid \psi \rangle
 =  \langle \phi \mid E \psi \rangle
 = \lim \langle \phi \mid \pi(g_j) \psi \rangle
 = \lim \langle \pi(g_j^{-1}) \phi \mid  \psi \rangle ,$$
 so the same argument, with $E^*$ in the place of~$E$ and $g_j^{-1}$ in the
place of~$g_j$, shows that $E^* \bigl( (\hilbert^U)^\perp \bigr) = 0$.

 Because $\pi$ is unitary, we know that $\pi(g_j)$ is
normal (that is, commutes with its adjoint) for every~$j$;
thus, the limit~$E$ is also normal: we have $E^* E = E
E^*$. Therefore
 $$ \lVert E\phi\rVert^2
 = \langle E\phi \mid E \phi \rangle
 = \langle (E^* E)\phi \mid \phi \rangle
 = \langle (E E^*)\phi \mid \phi \rangle
 = \langle E^*\phi \mid E^* \phi \rangle
 = \lVert E^* \phi\rVert^2 ,$$
 so $\ker E = \ker E^* $. 

Thus, 
 $$\ker E
 = \ker E + \ker E^*
 \supset (\hilbert^{U^-})^\perp + (\hilbert^U)^\perp
 = (\hilbert^{U^-} \cap \hilbert^U)^\perp
 = (\hilbert^{\langle U, U^- \rangle})^\perp
 .$$
 By passing to a subsequence of $\{g_j\}$, we may assume $\langle U, U^-
\rangle = G$ \see{horosubseq}. Then $\hilbert^{\langle U, U^- \rangle} =
\hilbert^G = 0$,
 so $\ker E
 \supset 0^\perp
 = \hilbert$.
 Hence, for all $\phi,\psi \in \hilbert$, we have
 $$ \lim \langle \pi(g_j) \phi \mid \psi \rangle
 = \langle E \phi \mid \psi \rangle
 = \langle 0 \mid \psi \rangle
 = 0 ,$$
 as desired.

\begin{case}
 The general case.
 \end{case}
 From the Cartan Decomposition $G = KAK$, we may write $g_j = c_j a_j c'_j$,
with $c_j,c'_j \in K$ and $a_j \in A$. Because $K$ is compact, we may assume,
by passing to a subsequence, that $\{c_j\}$ and $\{c'_j\}$ converge: say, $c_j
\to c$ and $c'_j \to c'$. Then
 \begin{align*}
 \lim \langle \pi(g_j) \phi \mid \psi \rangle
 &= \lim \langle \pi(c_j a_j c'_j) \phi \mid \psi \rangle \\
 &= \lim \langle \pi(a_j)  \pi(c'_j) \phi \mid \pi(c_j)^{-1}
\psi \rangle \\
 &= \lim \bigl\langle \pi(a_j) \bigl( \pi(c') \phi \bigr) 
 \mathrel{\big|} \pi(c)^{-1}
\psi \bigr\rangle \\
 &= 0 ,
 \end{align*}
 by Case~\ref{matcoeffs->0-A}.
 \end{proof}

The following example illustrates Lemma~\ref{horosubseq}.

\begin{eg} 
 Let $G = \SL(3,\real)$, define $H_1$ as in
Proposition~\ref{wallinSL3}, and suppose $\{g_j\}$ is some sequence of
elements of~$H_1$, such that $\lVert g_j \rVert \to \infty$. We may write
 $$ g_j = 
 \begin{pmatrix}
 e^{t_j} & 0 & 0 \\
 0 & e^{t_j} & 0 \\
 0 & 0 & e^{-2t_j}
 \end{pmatrix}
 , $$
 where $t_j \in \real$.
 By passing to a subsequence, we may assume that either $t_j \to \infty$ or
$t_j \to -\infty$. If $t_j \to \infty$, then, in the notation
of~\pref{horodefn}, we have
 $$
 U = 
 \left\{
 \begin{pmatrix}
 1 & 0 & * \\
 0 & 1 & * \\
 0 & 0 & 1
 \end{pmatrix}
 \right\}
 \mbox{ \qquad and \qquad}
 U^- = 
 \left\{
 \begin{pmatrix}
 1 & 0 & 0 \\
 0 & 1 & 0 \\
 * & * & 1
 \end{pmatrix}
 \right\}
 ;$$
 if $t_j \to -\infty$, then $U$ and~$U^-$ are interchanged.
 Thus, in either case, $\Lie U$ is the sum of two root spaces of~$\Lie G$, and
$\Lie U^-$ is the sum of the two opposite root spaces. It is not difficult to
see that $[\Lie U, \Lie U^-]$ is the sum of~$\Lie A$ and the remaining two
root spaces. Therefore, we have $\langle \Lie U, \Lie U^- \rangle = \Lie G$,
so $\langle U, U^- \rangle = G$.
 \end{eg}

\begin{lem} \label{horosubseq}
 If $G$ and~$\{g_j\}$ are as in Theorem~\ref{matcoeffs->0}, and $\{g_j\}
\subset A$, then, after replacing $\{g_j\}$ by a subsequence, we have $\langle
U,U^-\rangle = G$, where $U$ and~$U^-$ are defined in \pref{horodefn}.
 \end{lem}

\begin{proof}
 By passing to a subsequence, we may assume $\{g_j\}$ is contained in a single
Weyl chamber, which we may take to be $A^+$. Then, by passing to a
subsequence yet again, we may assume, for every positive real root~$\alpha$,
that either $\alpha(g_j) \to \infty$ or $\alpha(g_j)$ is bounded. Let
 \begin{itemize}
 \item $\Phi^+$ be the set of positive real roots;
 \item $\Delta$ be the set of positive simple real roots;
 \item $ \Psi = \{\, \alpha \in \Phi^+ \mid 
 \mbox{$\alpha(g_j)$ is bounded} \,\}$;
 and
 \item $T = \cap_{\psi\in \Psi} \ker \psi = \cap_{\psi\in \Psi \cap \Delta}
\ker \psi$. 
 \end{itemize}
 There is a compact subset~$C$ of~$A$, such that $\{g_j\} \subset C T$, so,
because $\lVert g_j \rVert \to \infty$, we know that $T$ is not trivial.

For each real root~$\alpha$, let $\rsp_\alpha$ be the corresponding root
subspace of~$\Lie G$. Then
 $$ \Lie U = \bigoplus_{\alpha \in \Phi^+ \smallsetminus \Psi} \rsp_\alpha
 \mbox{ \qquad and \qquad}
 \Lie U^- = \bigoplus_{\alpha \in \Phi^+ \smallsetminus \Psi} \rsp_{-\alpha}
 .$$
 Now, for $\alpha \in \Phi^+$, we have $\alpha \in \Psi$ if and only if
$\alpha$ is in the linear span of $\Psi \cap \Delta$. Thus, we see that $\Lie
U$ is precisely the unipotent radical of the standard parabolic subalgebra
$\Lie P = \Lie C_{\Lie G}(T) + \Lie A + \Lie N$ corresponding to the set
$\Psi \cap \Delta$ of simple roots \cite[4.2,
pp.~85--86]{BorelTits-Reductive}. 

Similarly, $\Lie U^-$ is the unipotent radical of the opposite parabolic
algebra $\Lie P^- = \Lie C_{\Lie G}(T) + \Lie A + \Lie N^-$. Because $G$ is
simple, the unipotent radicals of opposite parabolics generate~$\Lie G$
\cite[Prop.~4.11, p.~89]{BorelTits-Reductive}, so $\langle U, U^- \rangle =
G$, as desired.
 \end{proof}

\begin{cor}[{(of Theorem~\ref{expdecay})}] \label{1D->tempered}
 Assume $G$ is simple, and $\Rrank G \ge 2$. If $H$ is a one-parameter
subgroup of $AN$, then either
 \begin{enumerate}
 \item $H$ is tempered; or
 \item $H \subset N$.
 \end{enumerate}
 \end{cor}

\begin{proof}[{Proof \rm {\cite[Prop.~3.7]{OhWitte-CK}}}]
 Write $H = \{h^t\}$. From the Real Jordan Decomposition \pref{JordanDecomp},
we may assume, after replacing~$H$ by a conjugate subgroup, that $h^t = a^t
u^t$, where $a^t\in A$ is a hyperbolic one-parameter subgroup, and $u^t\in N$
is a unipotent one-parameter subgroup, such that $a^t$ and~$u^t$ commute with
each other.

We may assume $H \not\subset N$, so $a^t$ is nontrivial. Since the growth of
the hyperbolic one-parameter subgroup~$a^t$ is exponential, while that of the
unipotent one-parameter subgroup~$u^t$ is polynomial, there is some $\epsilon
> 0$, such that 
 $$ \lVert h^t\rVert = \lVert a^t u^t\rVert > \sqrt{\lVert a^t\rVert} >
e^{\epsilon |t|}$$
 for large $t \in \real$. 
 Since the function $C/e^{p \epsilon |t|}$ is in $L^1(\real)$, it follows
from Theorem~\ref{expdecay} that $H$ is tempered, as desired.
 \end{proof}

\begin{lem} \label{1DinN->notess}
  If $d(H) = 1$ and $H \subset N$, then $G/H$ does not have a
tessellation.
 \end{lem}

\begin{proof}[Proof {\rm \cite[Prop.~3.7]{OhWitte-CK}}]
 We have $\dim H = d(H) = 1$ \see{d(H)=dimH}, so $H$ is a connected,
one-dimensional, unipotent subgroup. Hence, the Jacobson-Morosov Lemma
\cite[Thm.~9.7.4, p.~432]{HelgasonBook} implies that there exists a connected,
closed subgroup $H_1$ of~$G$, such that $H_1$ contains~$H$, and $H_1$ is
locally isomorphic to $\SL(2,\real)$. Then $H$ is a Cartan-decomposition
subgroup of~$H_1$ \see{Rrank1-CDS}, so there is a compact subset~$C$ of~$A$,
such that $\mu(H_1) \subset \mu(H) C$ \see{bddchange}. Also, we have $d(H_1)
= 2 > 1 = d(H)$. Therefore, Theorem~\ref{noncpct-dim-notess} applies.
 \end{proof}

\begin{cor}[{(Oh-Witte \cite[Prop.~3.7]{OhWitte-CK})}] \label{1Dsimple}
 If $d(H) = 1$ and $G$ is simple, then $G/H$ does not have a tessellation.
 \end{cor}

\begin{proof}
 We may assume $H \subset AN$ \see{HcanbeAN}, so $\dim H = d(H) = 1$
\see{d(H)=dimH}. 
 \begin{itemize}
 \item If $\Rrank G < 2$, then Lemma~\ref{1Drank1} applies.
 \item If $H \subset N$, then Lemma~\ref{1DinN->notess} applies.
 \item If $\Rrank G \ge 2$ and $H \not\subset N$, then
Corollary~\ref{1D->tempered} implies that $H$ is tempered, so
Theorem~\ref{MargulisTempered} applies.
 \end{itemize}
 \end{proof}

\section{Homogeneous spaces of $\SL(3,\real)$,  $\SL(3,\complex)$, and
$\SL(3,\quaternion)$} \label{SL3notessSect}

Y.~Benoist \cite[Cor.~1]{Benoist} and G.A.~Margulis (unpublished) proved
(independently) that $\SL(3,\real)/\SL(2,\real)$ does not have a tessellation.
Using Benoist's method, H.~Oh and D.~Witte \cite[Prop.~1.10]{OhWitte-CK}
generalized this result by replacing $\SL(2,\real)$ with any closed, connected
subgroup~$H$, such that neither $H$ nor $\SL(3,\real)/H$ is compact. The same
argument applies even if $\real$ replaced with either $\complex$
or~$\quaternion$. However, the proof of Benoist (which applies in a more
general context) relies on a somewhat lengthy argument to establish one
particular lemma. Here, we adapt Benoist's method to obtain a short proof of
Theorem~\ref{SL3->notess} that avoids any appeal to the lemma.

\begin{notation} \label{SL3-oppinv}
  Assume $G = \SL(3,\F)$, for $\F = \real$, $\complex$,
or~$\quaternion$.
 \begin{itemize}
 \item Let $\tau$ be the opposition involution in~$A^+$ \see{oppinv};
 \item Let $B^+=\{\,a\in A^+ \mid \tau(a)=a \,\}$.
 \end{itemize}
 \end{notation}

More concretely, we have
  \begin{gather*}
   A^+ =
  \bigset{
  \begin{pmatrix}
  a_1 & 0 & 0 \\
  0 & a_2 & 0 \\
  0 & 0 & a_3
  \end{pmatrix}
  }{
  \begin{matrix}
  a_1,a_2,a_3 \in \real^+ , \\
  a_1 a_2 a_3 = 1 , \\
  a_1 \ge a_2 \ge a_3
  \end{matrix}
  }; \\
 \tau
  \begin{pmatrix}
  a_1 & 0 & 0 \\
  0 & a_2 & 0 \\
  0 & 0 & a_3
  \end{pmatrix}
  =
  \begin{pmatrix}
  a_3^{-1} & 0 & 0 \\
  0 & a_2^{-1} & 0 \\
  0 & 0 & a_1^{-1}
  \end{pmatrix}
  ; \\
  B^+ =
  \bigset{
  \begin{pmatrix}
  a & 0 & 0 \\
  0 & 1 & 0 \\
  0 & 0 & a^{-1}
  \end{pmatrix}
  }{
  \begin{matrix}
 \\
  a \ge 1  \\
  \\
 \end{matrix}
  }.
  \end{gather*}

\begin{lem} \label{1DinSL3}
 If $G = \SL(3,\F)$ and $d(H) = 1$, then $G/H$ does not have a tessellation.
 \end{lem}

\begin{proof}
 Since $\Rrank G = 2$, the desired conclusion follows from
Corollary~\ref{1Drank2}; since $G$ is simple, it also follows from
Corollary~\ref{1Dsimple}.
 However, we give a proof that requires only the special case described in
Proposition~\ref{wallinSL3}, rather than the full strength of \pref{1Drank2}
or~\pref{1Dsimple}.

Suppose $\Gamma$ is a crystallographic group for $G/H$. (This will lead to a
contradiction.) Let $L_1$ and~$L_2$ be the two walls of~$A^+$. From
Proposition~\ref{tess->misswall}, we know that there exists $k \in \{1,2\}$,
such that $\mu(\Gamma) C \cap L_k$ is finite, for every compact subset~$C$
of~$A$. 

Because $\Gamma^{-1} = \Gamma$, we have $\tau \bigl( \mu(\Gamma) \bigr) =
\mu(\Gamma)$. On the other hand, $\tau$~interchanges $L_1$ and~$L_2$. Thus,
the preceding paragraph implies that $\mu(\Gamma) C \cap (L_1 \cup L_2)$ is
finite, for every compact subset~$C$ of~$A$. 

For $H_1$ as in \pref{SL3wallSubgrp}, we have $\mu(H_1) = L_1 \cup L_2$, so
the conclusion of the preceding paragraph implies that $\Gamma$ acts properly
discontinuously on $G/H_1$ \see{proper<>mu(L)}. Now
Theorem~\fullref{noncpctdim}{tess} implies $\Gamma \backslash G/H_1$ is
compact; thus, $G/H_1$ has a tessellation. This contradicts
Proposition~\ref{wallinSL3}.
 \end{proof}

For completeness, we include the proof of the following simple proposition.

\begin{prop}[{\cite[Prop.~7.3]{OhWitte-CK}}] \label{SL3-B+}
 Assume $G = \SL(3,\F)$.
 If $H$ is a closed, connected subgroup of~$AN$ with $\dim H \ge  2$, then
$B^+ \subset \mu(H)$.
 \end{prop}

\begin{center}
 \begin{figure}
 \includegraphics[scale=0.5]{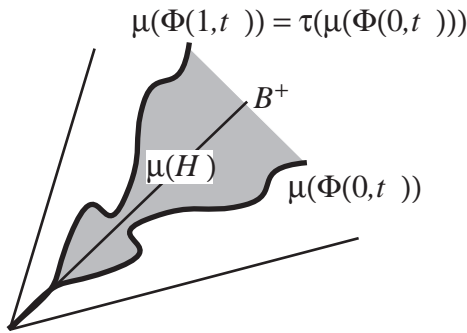}
 \caption{Proposition~\ref{SL3-B+}: $\mu(H)$ is on both sides of~$B^+$, so it
must contain $B^+$.}
 \label{SL3B+Fig}
 \end{figure}
 \end{center}

\begin{proof}[{Proof \rm (cf.\ Figure~\ref{SL3B+Fig} and proof of
\ref{CDS<>h_m})}]
 Since $H \subset AN$ and $\dim H \ge 2$, it is easy to construct a
continuous, proper map $\Phi\colon [0,1] \times \real ^+ \to H$ such that
$\Phi(1,t) = \Phi(0,t)^{-1}$, for all $t \in \real^+$ (cf.\
Figure~\ref{ConstructPhiFig}). For example, choose two linearly independent
elements $u$ and~$v$ of~$\Lie H$, and define
 $$ \Phi(s,t) = \exp \bigl( t \cos(\pi s) u + t \sin(\pi s) v \bigr) .$$

\begin{center}
 \begin{figure}
 \includegraphics[scale=0.5]{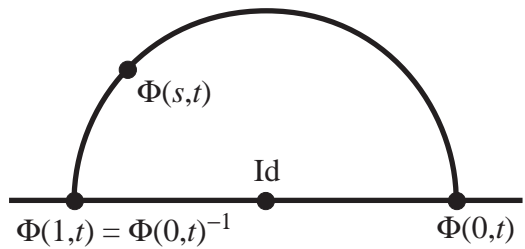}
 \caption{Construction of $\Phi(s,t)$.}
 \label{ConstructPhiFig}
 \end{figure}
 \end{center}

 If we identify~$A$ with its Lie algebra~$\Lie a$, then~$A^+$ is a convex cone
in~$\Lie a$ and the opposition involution~$\tau$ is the reflection in~$A^+$
across the ray~$B^+$. Thus, for any $a\in A^+$, the points $a$ and~$\tau(a)$
are on opposite sides of~$B^+$, so any continuous curve in~$A^+$ from~$a$
to~$\tau(a)$ must intersect~$B^+$. In particular, for each $t \in \real^+$,
the curve 
 $$ \{\, \mu \bigl( \Phi(s,t) \bigr) \mid 0 \le s \le 1\,\} $$ from $\mu
\bigl( \Phi(0,t) \bigr)$ to $\mu \bigl( \Phi(1,t) \bigr)$ must
intersect~$B^+$. Thus, we see, from an elementary continuity argument, that
$\mu \bigl[ \Phi\bigl( [0,1] \times \real ^+ \bigr) \bigr]$ contains~$B^+$.
Therefore, $B^+$ is contained in $\mu(H)$.
 \end{proof}

\begin{center}
 \begin{figure}
 \includegraphics[scale=0.5]{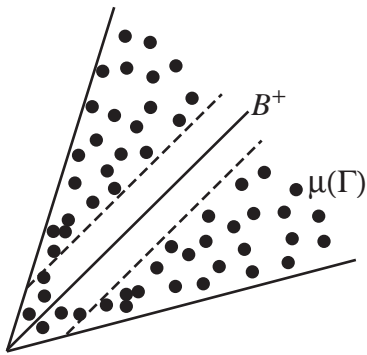}
 \caption{Proof of Theorem~\ref{SL3->notess}: $\mu(\Gamma)$ stays away
from~$B^+$, because $B^+ \subset \mu(H)$. Also, half of $\mu(\Gamma)$ is on
each side of~$B^+$, because $\tau \bigl( \mu(\Gamma) \bigr) = \mu(\Gamma)$.
This contradicts the fact that $\Gamma$ has only one end.}
 \label{SL3notessFig}
 \end{figure}
 \end{center}

\begin{proof}[{\bf Proof of Theorem~\ref{SL3->notess} \rm (cf.\
Figure~\ref{SL3notessFig} and proof of \ref{tess->dim>1,2})}]
  Suppose $\Gamma$ is a crystallographic group for $G/H$. (This will lead to
a contradiction.) We may assume $H \subset AN$ \see{HcanbeAN}.

 Let $F$ be a (symmetric) finite generating set for~$\Gamma$, and choose a
compact, convex, symmetric subset~$C$ of~$A$ so large that 
 $$\mu(\gamma F) \subset \mu(\gamma) C$$
 for every $\gamma \in \Gamma$ \see{bddchange}.

 From Lemma~\ref{1DinSL3}, we know that $\dim H \ge 2$, so
Proposition~\ref{SL3-B+} implies that $B^+ \subset \mu(H)$. Then, because
$\Gamma$ acts properly on $G/H$, we conclude that $\mu(\Gamma) \cap B^+ C$ is
finite \see{proper<>mu(L)}. Since $\mu$ is a proper map, this implies that
$\Gamma \cap \mu^{-1}(B^+ C)$ is finite.

Let $A_1$ and $A_2$ be the two components of $A^+ \smallsetminus B^+$. Because
$\Gamma^{-1} = \Gamma$, we know that $\tau \bigl( \mu(\Gamma) \bigr) =
\mu(\Gamma)$. Then, because $\tau$ interchanges $A_1$ and $A_2$, we conclude
that $\tau \bigl( \mu(\Gamma) \cap A_1 \bigr) = \mu(\Gamma) \cap A_2$.
Therefore, $\mu(\Gamma) \cap A_1$ and $\mu(\Gamma) \cap A_2$ have the same
cardinality, so they must both be infinite. So 
 $$ \mbox{each of $\Gamma \cap \mu^{-1}(A_1)$ and $\Gamma \cap \mu^{-1}(A_2)$
is infinite.} $$
 Because $\Gamma$ has only one end \see{tess->1end}, this implies there exist 
 $$\gamma \in \bigl( \Gamma \cap \mu^{-1}(A_1) \bigr) \smallsetminus \mu^{-1}(B^+
C) ,$$
 such that 
 \begin{equation} \label{SL3->notessPf-gammaf}
 \gamma f \in \bigl( \Gamma \cap \mu^{-1}(A_2) \bigr)
\smallsetminus \mu^{-1}(B^+ C) ,
 \end{equation}
 for some $f \in F$.
 Then $\mu(\gamma) \in A_1$, $\mu(\gamma f) \in A_2$, and
 $$ \mu(\gamma f) \in \mu(\gamma F) \subset \mu(\gamma) C .$$
 Using the fact that $C$ is symmetric and the fact that $C$ contains the
identity element~$e$, we conclude that
 $$ \mu(\gamma) \in \bigl( \mu(\gamma f) C \bigr) \cap A_1
 \mbox{\qquad and\qquad}
 \mu(\gamma f) \in \bigl( \mu(\gamma f) C \bigr) \cap A_2
 ; $$
 therefore $\mu(\gamma f) C$ intersects both~$A_1$ and~$A_2$.
 Since $B^+$ separates $A_1$ from~$A_2$, and $C$ is connected, this implies
that $\mu(\gamma f) C$ intersects $B^+$; hence
 $\mu(\gamma f) \in B^+ C$. This contradicts the fact that $\gamma f \notin
\mu^{-1}(B^+ C)$ \see{SL3->notessPf-gammaf}.
 \end{proof}

\section{Explicit coordinates on $\so(2,n)$ and $\su(2,n)$} \label{coordsSect}

From this point on, we focus almost entirely on $\SO(2,n)$ and $\SU(2,n)$.
(The only exception is that some of the examples constructed in
Section~\ref{ExistenceSect} are for other groups.) In this section, we define
the group $\SU(2,n;\F)$, which allows us to provide a fairly unified treatment
of $\SO(2,n)$ and $\SU(2,n)$ in later sections.

\subsection{The group $\SU(2,n;\F)$}

\begin{notation} \ \hbox{ }
 \begin{itemize}
 \item We use $\F$ to denote either $\real$ or~$\complex$.

 \item Let $\df = \dimR \F$, so $\df \in \{1,2\}$.

 \item We use $\Fim$ to denote the purely imaginary elements of~$\F$, so
 $$ \Fim =
 \begin{cases}
 \hfil 0 & \mbox{if $\F = \real$} \\
 i \real & \mbox{if $\F = \complex$} .
 \end{cases}
 $$

\item For $\phi \in \F$, there exist unique $\Re \phi \in \real$ and $\Im
\phi \in \Fim$, such that $\phi = \Re \phi + \Im \phi$. (\textbf{Warning:}
in our notation, the imaginary part of $a+bi$ is $bi$, \emph{not}~$b$.)

 \item For $\phi \in \F$, we use $\cjg{\phi}$ to denote the conjugate $\Re
\phi - \Im \phi$ of~$\phi$. (If $\F = \real$, then $\cjg{\phi} = \phi$.)

 \item For a row vector $x \in \F^{n-2}$, or, more generally, for any
matrix~$x$ with entries in~$\F$, we use $x^\dagger$ to denote the
conjugate-transpose of~$x$.
 \end{itemize}
 \end{notation}

\begin{notation} \label{SUFDefn}
 For
 $$ J =
 \begin{pmatrix}
 0      &      0 & 0 & \cdots & 0 & 0      & 1      \\
 0      &      0 & 0 & \cdots & 0 & 1      & 0      \\
 0      &      0 &   &        &   & 0      & 0      \\
 \vdots & \vdots &   & \Id    &   & \vdots & \vdots \\
 0      &      0 &   &        &   & 0      & 0      \\
 0      &      1 & 0 & \cdots & 0 & 0      & 0      \\
 1      &      0 & 0 & \cdots & 0 & 0      & 0
 \end{pmatrix}
 \in \SL(n+2, \F) ,$$
 we define 
 $$\SU(2,n;\F)
 = \{\, g \in \SL(n+2,\F) \mid g J g^\dagger = J \,\} $$
and
 $$\su(2,n;\F)
 = \{\, u \in \Liesl(n+2,\F) \mid u J + J u^\dagger = 0 \,\}. $$
 Then:
 \begin{itemize}
 \item $\SU(2,n;\real)$ is a realization of $\SO(2,n)$, 
 \item $\SU(2,n;\complex)$ is a realization of $\SU(2,n)$, and
 \item $\su(2,n;\F)$ is the Lie algebra of $\SU(2,n;\F)$.
 \end{itemize}
 We choose
 \begin{itemize}
 \item $A$ to
consist of the diagonal matrices in $\SU(2,n;\F)$ that have nonnegative real
entries,
 \item $N$ to consist of the upper-triangular matrices in
$\SU(2,n;\F)$ with only $1$'s on the diagonal,
 and
 \item $K = \SU(2,n;\F) \cap \SU(n+2)$.
 \end{itemize}
 A straightforward matrix calculation shows that the Lie algebra of $AN$ is
 \begin{equation}
 \label{SUF-AN}
 \Lie a + \Lie N =
 \bigset{
 \begin{pmatrix}
 t_1 & \phi & x & \eta & \xx \\
 0   & t_2  & y & \yy & -\cjg{\eta} \\
 0   & 0    & 0 & -y^{\dagger} & -x^{\dagger} \\
 0   & 0   & 0 & -t_2 & -\cjg{\phi} \\
 0   & 0   & 0 & 0 & -t_1 \\
 \end{pmatrix}
  }{
 \begin{matrix}
 t_1,t_2 \in \real, \\
 \phi,\eta \in \F, \\
 x,y \in \F^{n-2}, \\
 \xx, \yy \in \Fim
 \end{matrix}}
 .
 \end{equation}
 \end{notation}

\begin{rem} \label{tworows}
 From \eqref{SUF-AN}, we see that the first two rows of any element
of $\Lie a + \Lie N$ are sufficient to determine the entire matrix. In fact,
it is also not necessary to specify the last entry of the second row of the
matrix.
 \end{rem}

\begin{rem} \label{d(SU2)}
 From \pref{d(G)} and \pref{SUF-AN}, we see that $d \bigl( \SU(2,n;\F) \bigr)
=  \dim( \Lie A +
\Lie N) = 2 \df n$.
 \end{rem}

\begin{notation}
 Because $N$ is simply connected and nilpotent, the exponential map is a
diffeomorphism from~$\Lie N$ to~$N$ (indeed, its inverse, the
logarithm map, is a polynomial \cite[Thm.~8.1.1, p.~107]{Hochschild-Algic}, so
each element of~$N$ has a unique representation in the form $\exp u$ with $u
\in \Lie N$. Thus, each element~$h$ of~$N$ determines corresponding values of
$\phi$, $x$, $y$, $\eta$, $\xx$ and~$\yy$ (with $t_1 = t_2 = 0$). We write
 $$ \phi_h, x_h, y_h, \eta_h, \xx_h, \yy_h $$
 for these values.
 \end{notation}

\begin{center}
 \begin{figure}
 \includegraphics[scale=0.5]{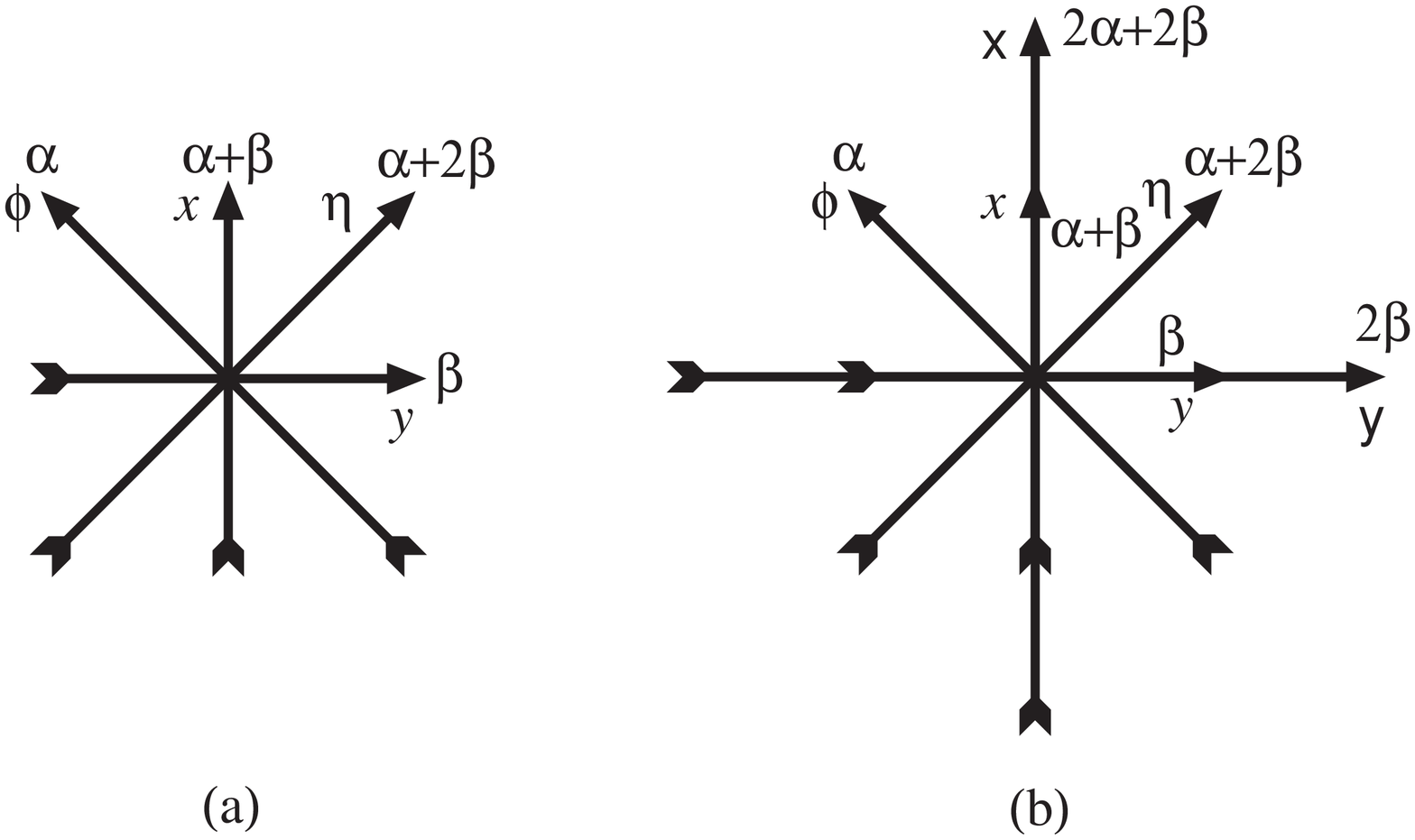}
 \caption{The real root systems of (a)~$\SU(2,n;\real) = \SO(2,n)$ and
(b)~$\SU(2,n;\complex) = \SU(2,n)$.}
 \label{rootspict}
 \end{figure}
 \end{center}

\begin{notation} \label{simpleroots}
 We let $\alpha$ and~$\beta$ be the simple real roots of $\SU(2,n;\F)$,
defined by 
 $$ \mbox{$\alpha(a) = a_{1,1}/a_{2,2}$ and $\beta(a) = a_{2,2}$} ,$$
 for a (diagonal) element~$a$ of~$A$.
 Thus, the positive real roots (see Figure~\ref{rootspict}) are
 $$ \begin{cases}
 \alpha, \beta, \alpha+\beta, \alpha+2\beta,
 & \mbox{if $\F = \real$} \\
 \alpha, \beta, \alpha+\beta, \alpha+2\beta, 2\beta, 2\alpha+2\beta  &
\mbox{if $\F = \complex$} . \\
 \end{cases} $$
 Concretely:
 \begin{itemize}
 \item the root space $\rsp_\alpha$ is the $\phi$-subspace in~$\Lie N$,
 \item the root space $\rsp_\beta$ is the $y$-subspace in~$\Lie N$,
 \item the root space $\rsp_{\alpha+\beta}$ is the $x$-subspace in~$\Lie N$,
 \item  the root space $\rsp_{\alpha+2\beta}$ is the $\eta$-subspace
in~$\Lie N$,
 \item  the root space $\rsp_{2\beta}$ is the $\yy$-subspace
in~$\Lie N$ (this is~0 if $\F = \real$), and
 \item  the root space $\rsp_{2\alpha+2\beta}$ is the $\xx$-subspace
in~$\Lie N$ (this is~0 if $\F = \real$).
 \end{itemize}
 \end{notation}

\begin{defn}
 Let
 \begin{align*}
 \Dn &= \rsp_{\alpha+2\beta} + \rsp_{2\alpha+2\beta} +
\rsp_{2\beta} \\
 &= \{\, z \in \Lie N \mid \phi_z = 0, x_z = y_z = 0 \,\} \\
 &= 
 \bigset{
  \begin{pmatrix}
 0 & 0    & 0 & \eta & \xx \\
   & 0    & 0 & \yy & -\cjg{\eta} \\
 & & \dots
 \end{pmatrix}
 }{
  \begin{matrix}
  \eta \in \F, \\
  \xx, \yy \in \Fim
  \end{matrix}
  }
 .
 \end{align*}
 and, for a given Lie algebra~$\Lie H \subset \Lie N$,
 $$ \Dh = \Dn \cap \Lie H.$$
 Note that if $\phi_u = 0$ for every $u \in \Lie H$, then $[\Lie H, \Lie H]
\subset \Dh$ and $\Dh$ is contained in the center of~$\Lie H$
\cf{[u,v]}.
 \end{defn}

\begin{rem}
 By definition \see{A+Defn}, we have 
 $$ A^+ = \{\, a \in A \mid \alpha(a) \ge 1, \beta(a) \ge 1 \,\} .$$
 Therefore, from the definition of~$\alpha$ and~$\beta$ \see{simpleroots}, we
see that
 \begin{equation} \label{A+}
 A^+ = \{\, a \in A \mid a_{1,1} \ge a_{2,2} \ge 1
\,\} .
 \end{equation}
 \end{rem}

\begin{rem}
  For $\F = \quaternion$, the division algebra of real quaternions, the group
$\SU(2,n;\quaternion)$ is a realization of $\Symp(2,n)$. Most of the work
in this paper carries over, but the upper bound on $\dim H$ given in
Theorem~\ref{maxnolinear} is not sharp in this case (and it does not seem to
be easy to improve this result to obtain a sharp bound). Thus, we have not
obtained any interesting conclusions about the nonexistence of tessellations
of homogeneous spaces of $\Symp(2,n)$.
 \end{rem}

\subsection{The subgroups $\SU(1,n;\F)$ and $\Symp(1,m;\F)$}

We now describe how the four important families of homogeneous spaces of
Example~\ref{KulkarniEg} are realized in terms of $\SU(2,n;\F)$.

\begin{defn} \label{SU1nDefn}
 Let
 \begin{itemize}
 \item $\SU(1,n;\real) = \SO(1,n)$;
 \item  $\Symp(1,n;\real) = \SU(1,n)$;
 \item $\SU(1,n;\complex) = \SU(1,n)$; and
 \item $\Symp(1,n;\complex) = \Symp(1,n)$.
 \end{itemize}
 Then, for an appropriate choice of the embeddings in
Example~\ref{KulkarniEg}, we have
 \begin{equation} \label{SU1nAN}
 \su(1,n;\F) \cap (\Lie A + \Lie N) =
 \bigset{
  \begin{pmatrix}
 t & \phi & x & \phi & \xx \\
 0 & 0    & 0 & 0 & -\cjg{\phi} \\
 & & \dots
 \end{pmatrix}
 }{
  \begin{matrix}
  t \in \real, \\
  \phi \in \F, \\
  x \in \F^{n-2}, \\
  \xx \in \Fim
  \end{matrix}
  }
 \end{equation}
 and (if $2m \le n$) we have
  \begin{multline} \label{Sp1mAN}
  \Liesp(1,m;\F) \cap (\Lie A + \Lie N)
 =  \\
 \setcounter{MaxMatrixCols}{15}  
 \bigset{
  \begin{pmatrix}
 t & 0    & x_1 & x_2 & x_3 & x_4 & \dots & x_{2m-3} & x_{2m-2} & 0 
& \dots & 0
& \eta & \xx \\
 0 & t    & -\cjg{x_2} & \cjg{x_1} & -\cjg{x_4} & \cjg{x_3} & \dots &
-\cjg{x_{2m-2}} & \cjg{x_{2m-3}} & 0 & \dots & 0 & -\xx & -\cjg{\eta} \\
 & & & & & & & \dots
 \end{pmatrix}
 }{
  \begin{matrix}
  t \in \real, \\
  x_j \in \F, \\
  \eta \in \F, \\
  \xx \in \Fim
  \end{matrix}
  }
 .
 \end{multline}
 \end{defn}

\begin{rem} \label{d(Sp)}
 From \pref{d(G)}, \pref{SU1nAN} and \pref{Sp1mAN}, we see that
 \begin{itemize}
 \item $ d \bigl( \SU(1,n;\F) \bigr)
 = \dim \bigl( \su(1,n;\F) \cap (\Lie A + \Lie N) \bigr)
 = \df n$ and
 \item $ d \bigl( \Symp(1,m;\F) \bigr)
 = \dim \bigl( \Liesp(1,m;\F) \cap (\Lie A + \Lie N) \bigr)
 = 2 \df m$.
 \end{itemize}
 \end{rem}

\subsection{Formulas for exponentials and brackets}

The arguments in later sections often require the calculation of $\exp
u$, for some $u \in \Lie N$, or of $[u,v]$, for some $u,v \in \Lie N$. We now
provide these calculations for the reader's convenience.

\begin{rem}
 For
 $$ u =
 \begin{pmatrix}
 0 & \phi & x & \eta & \xx \\
 0 & 0    & y & \yy & -\cjg{\eta} \\
 0 & 0    & 0 & -y^{\dagger} & -x^{\dagger} \\
 0 & 0    & 0 & 0 & -\cjg{\phi} \\
 0 & 0    & 0 & 0 & 0 \\
 \end{pmatrix}
 \in \Lie N
 ,$$
 we have
 \begin{equation} \label{exp}
 \exp(u)=
 \begin{pmatrix}
 \s 1 & \phi & x+\frac{1}{2} \phi y
      &  \eta -\frac{1}{2} x y^{\dagger} + \frac{1}{2} \phi \yy -
\frac{1}{6} \phi |y|^2
      & \vphantom{\vrule height 20pt depth 20pt}
 \genfrac{}{}{0pt}{0}
 {
 - \frac{1}{2} |x|^2
 - \Re(\phi \cjg{\eta})
 +\frac{1}{24} |\phi|^2 |y|^2
   }{
 + \left(
 \xx
 - \frac{1}{6} \phi \yy \cjg{\phi}
 + \frac{1}{3} \Im(x y^{\dagger} \cjg{\phi} )
 \right)
 } \\
 \s 0 & 1 & y
    & \yy -\frac{1}{2} |y|^2
    & -\cjg{\eta} -\frac{1}{2} y x^{\dagger}-\frac{1}{2} \yy \cjg{\phi}
+\frac{1}{6} |y|^2 \cjg{\phi}  \\
 \s 0 & 0 &\Id& -y^{\dagger} & -x^{\dagger}+\frac{1}{2}
y^{\dagger} \cjg{\phi} \\
 \s 0 & 0 & 0 & 1 & -\phi^{\dagger} \\
 \s 0 & 0 & 0 & 0 & 1 \\
 \end{pmatrix}
 .
 \end{equation}

 When $\phi = 0$, this simplifies to:
 \begin{equation} \label{exp(phi=0)}
 \exp(u)=
 \begin{pmatrix}
 \s 1 & 0 & x
      &  \eta -\frac{1}{2} x y^{\dagger}
      & \vphantom{\vrule height 20pt depth 20pt}
 \xx - \frac{1}{2} |x|^2 \\
 \s 0 & 1 & y
    & \yy -\frac{1}{2} |y|^2
    & -\cjg{\eta} -\frac{1}{2} y x^{\dagger}  \\
 \s 0 & 0 &\Id& -y^{\dagger} & -x^{\dagger} \\
 \s 0 & 0 & 0 & 1 & 0 \\
 \s 0 & 0 & 0 & 0 & 1 \\
 \end{pmatrix}
 . 
 \end{equation}
 Similarly, when $y = 0$, we have
 \begin{equation} \label{exp(y=0)}
 \exp(u)=
 \begin{pmatrix}
 \s 1 & \phi & x
      &  \eta  + \frac{1}{2} \phi \yy
      & \vphantom{\vrule height 20pt depth 20pt}
 \genfrac{}{}{0pt}{0}
 {
 - \frac{1}{2} |x|^2
 - \Re(\phi \cjg{\eta})
   }{
 + \left(
 \xx
 - \frac{1}{6} \phi \yy \cjg{\phi}
 \right)
 } \\
 \s 0 & 1 & 0
    & \yy
    & -\cjg{\eta} -\frac{1}{2} \yy \cjg{\phi} \\
 \s 0 & 0 &\Id& 0 & -x^{\dagger} \\
 \s 0 & 0 & 0 & 1 & -\phi^{\dagger} \\
 \s 0 & 0 & 0 & 0 & 1 \\
 \end{pmatrix}
 .
 \end{equation}
 \end{rem}

\begin{rem} 
 For
 \begin{equation} \label{umatrix}
 u =
 \begin{pmatrix}
 0 & \phi & x & \eta & \xx \\
   & 0 & y & \yy & -\cjg{\eta} \\
   &   & \cdots \\
 \end{pmatrix}
 \mbox{\qquad and\qquad}
\tilde u =
 \begin{pmatrix}
 0 & \tilde \phi & \tilde x & \tilde \eta & \tilde \xx \\
   & 0 & \tilde y & \tilde \yy & - \cjg{\tilde {\eta}} \\
   &   & \cdots \\
 \end{pmatrix}
 ,
 \end{equation}
 we have
 \begin{equation} \label{[u,v]}
 [u, \tilde u] =
 \begin{pmatrix}
 0 & 0 & \phi \tilde y - \tilde \phi y
 & - x \tilde y^{\dagger} + \tilde x y^{\dagger} + \phi \tilde \yy -
\tilde \phi \yy
 & -2 \Im(x \tilde x^{\dagger} + \phi \cjg{\tilde \eta}
 - \tilde \phi \cjg{\eta}) \\
   & 0 & 0 & -2 \Im (y \tilde y^{\dagger})
 &  \tilde y x^{\dagger} - y \tilde x^{\dagger} + \tilde \yy \cjg{\phi} -
 \yy \cjg{\tilde \phi} \\
  &  &  & \cdots \\
 \end{pmatrix}
 .
 \end{equation}
 \end{rem}

\begin{rem}[{\cite[Eq.~(2.13.8), p.~104]{Varadarajan}}] \label{conjugation}
 For $u,v \in \Lie H$, we have
 $$ \exp(-v) u \exp(v)
 = u
 + [u,v]
 + \frac{1}{2} \bigl[ [u,v] , v \bigr]
 + \frac{1}{3!} \Bigl[ \bigl[ [u,v] , v \bigr] , v \Bigr]
 + \cdots
 .$$
 Combining this with \eqref{[u,v]} allows us to calculate the effect of
conjugating by an element of~$N$. 

For example, suppose $u \in \Lie N$, with $\phi_u = 0$ and $y_u = 0$, and
suppose $v \in \rsp_{\alpha+\beta}$. We see, from~\pref{[u,v]}, that
$\phi_{[u,v]} = 0$ and that $x_{[u,v]} = y_{[u,v]} = 0$, so  
 $\bigl[ [u,v] , v \bigr] = 0$ \see{[u,v]}. Therefore
 $$  \exp(-v) u \exp(v) = u + [u,v] .$$
 \end{rem}

\section{Calculating the Cartan projection} \label{CalcSect}

Y.~Benoist \cite[Lem.~2.4]{Benoist} showed that calculating values of the
Cartan projection~$\mu$ is no more difficult than calculating the norm of a
matrix \see{mucalc}. In this section, we describe this elegant method and
some of its consequences, in the special case $G = \SU(2,n;\F)$.

\begin{standing}
 Throughout this section, we assume $G = \SU(2,n;\F)$.
 \end{standing}

\subsection{The basic definitions}

\begin{notation}
 We employ the usual Big Oh and little oh notation: for functions $f_1,f_2$
on a subset~$X$ of~$G$, we say
  $$ \mbox{\emph{$f_1 = O(f_2)$ for $h \in X$}} $$
  if there is a constant~$C$, such that, for all $h \in X$ with $\lVert
h\rVert$ large, we have $\lVert f_1(h)\rVert \le C \lVert f_2(h)\rVert$. (The
values of each~$f_j$ are assumed to belong to some finite-dimensional normed
vector space, typically either~$\complex$ or a space of complex matrices.
Which particular norm is used does not matter, because all norms are
equivalent up to a bounded factor.) We say
  $$ \mbox{\emph{$f_1 = o(f_2)$ for $h \in X$}} $$
  if $\lVert f_1(h)\rVert/\lVert f_2(h)\rVert \to 0$ as $h \to \infty$. (We
use $h \to \infty$ to mean $\lVert h\rVert \to \infty$.) Also, we write
  $$f_1 \asymp f_2$$
  if $f_1 = O(f_2)$ and $f_2 = O(f_1)$.
 \end{notation}

We use the following norm on $\SU(2,n;\F)$, because it is easy to calculate.
The reader is free to make a different choice, at the expense of changing $=$
to~$\asymp$ in a few of the calculations.

\begin{defn}
 For $h \in \SU(2,n;\F)$, we define $\lVert h\rVert$ to be the maximum
absolute value among the matrix entries of~$h$. That is,
 $$ \lVert h\rVert = \max_{1 \le j,k \le n+2} |h_{j,k}| .$$
 \end{defn}

\begin{defn} \label{rhoDefn}
 Define $\rho \colon \SU(2,n;\F) \to \GL(\F^{n+2} \wedge \F^{n+2})$
by $\rho(h) = h \wedge h$, so $\rho$ is the second exterior power of the
standard representation of $\SU(2,n;\F)$. Thus, we may define $\lVert \rho(h)
\rVert$ to be the maximum absolute value among the determinants of all the $2
\times 2$ submatrices of the matrix~$h$. That is,
 $$ \lVert\rho(h)\rVert = \max_{1 \le j,k,\ell,m \le n+2}
 \left| \det
 \begin{pmatrix}
 h_{j,k} & h_{j,\ell} \\
 h_{m,k} & h_{m,\ell}
 \end{pmatrix}
 \right|
 .$$
 \end{defn}

From \pref{exp}, \pref{exp(phi=0)}, and~\pref{exp(y=0)}, it is clear that the
$2 \times 2$ minor in the top right corner is often larger than the other $2
\times 2$ minors, so we give it a special name.

\begin{defn} \label{DeltaDefn}
 For $h \in \Mat_{n+2}(\F)$, define
 $$\Delta(h) =
 \det
 \begin{pmatrix}
 h_{1,n+1} & h_{1,n+2} \\
 h_{2,n+1} & h_{2,n+2}
 \end{pmatrix}
 . $$
 \end{defn}

\subsection{Y.~Benoist's method for using matrix norms to calculate~$\mu$}
\label{CalcSect-Benoist}

\begin{lem} \label{calcrho}
 For $a \in A^+$, we have
 $\lVert a\rVert = a_{1,1}$ and 
 $\lVert\rho(a)\rVert = a_{1,1} a_{2,2}$.
 \end{lem}

\begin{proof}
 From \pref{SUF-AN}, we see that
 \begin{equation} \label{ajj}
  a_{j,j} =
 \begin{cases}
 1 & \mbox{if $3 \le j \le n$} \\
 1/a_{2,2} & \mbox{if $j = n+1$} \\
 1/a_{1,1} & \mbox{if $j = n+2$} .
 \end{cases}
 \end{equation}
 Thus, from~\pref{A+}, we see that 
 $$a_{1,1} \ge a_{2,2} \ge a_{j,j}$$
 for $j \ge 3$ (and, since $a$ is diagonal, we have $a_{j,k} = 0$ for $j \neq
k$).
 Therefore, the desired conclusions follow from the definitions of $\lVert
a\rVert$ and $\lVert\rho(a)\rVert$.
 \end{proof}

\begin{prop}[{(Benoist, cf.\ \cite[Lem.~2.4]{Benoist})}] \label{mu=rho}
 We have 
 \begin{gather}
 \mu(h) \asymp h, \label{h=mu(h)} \\
 \rho \bigl( \mu(h) \bigr) \asymp \rho(h), \label{rho(h)=rho(mu)} \\
 \mbox{$\mu(h)_{1,1} \asymp \lVert h\rVert$,
 and
 $\mu(h)_{2,2} \asymp \lVert \rho(h) \rVert/\lVert h\rVert$,} \label{mucalc}
 \end{gather}
 for $h \in \SU(2,n;\F)$.
 \end{prop}

\begin{proof}
 Choose $k_1,k_2 \in K$, such that $\mu(h) = k_1 h k_2$. Because
 $\lVert xy\rVert = O \bigl( \lVert x\rVert \lVert y\rVert \bigr)$
 for $x,y \in \SU(2,n;\F)$, and 
 $\max_{k \in K} \lVert k\rVert < \infty$
 (since $K$ is compact), we have
 $$ \lVert \mu(h) \rVert = \lVert k_1 h k_2 \rVert = O \bigl( \lVert h\rVert
\bigr) $$
 and
 $$ \lVert h\rVert = \lVert k_1^{-1} \mu(h) k_2^{-1} \rVert = O \bigl(
\lVert\mu(h)\rVert \bigr), $$
 so \pref{h=mu(h)} holds.
 Similarly, we have
 $$ \lVert \rho \bigl( \mu(h) \bigr) \rVert
 = \lVert \rho(k_1) \rho(h) \rho(k_2) \rVert
 \asymp \lVert\rho(h)\rVert ,$$
 so \pref{rho(h)=rho(mu)} holds.

 For $a \in A^+$, we know, from \pref{calcrho}, that
 $a_{1,1} = \lVert a\rVert$ and
 $ a_{2,2} = \lVert\rho(a)\rVert / a_{1,1}$.
 Thus, letting $a = \mu(h)$, and using \pref{h=mu(h)}
and~\pref{rho(h)=rho(mu)}, we see that
 $$ \mu(h)_{1,1} = \lVert\mu(h)\rVert \asymp \lVert h\rVert $$
 and
 $$ \mu(h)_{2,2} =
 \frac{\lVert \rho \bigl( \mu(h) \bigr) \rVert} {\mu(h)_{1,1}}
 \asymp \frac{\lVert\rho(h)\rVert}{\lVert h\rVert} ,$$
 as desired.
 \end{proof}

\begin{rem}
 Proposition~\ref{mu=rho} generalizes to any reductive group~$G$
\cite[Lem.~2.3]{Benoist}. However, one may need to use a different
representation in the place of~$\rho$. In fact, if $\Rrank G = r$, then
$r$~representations of~$G$ are needed; for $G = \SU(2,n;\F)$, we have $\Rrank
G = 2$, and the two representations we use are $\rho$ and the identity
representation $I(h) = h$.
 \end{rem}

\begin{cor} \label{gn=hn}
 Let $g_n \to \infty$ and $h_n \to \infty$ be two sequences of elements of
$\SU(2,n;\F)$.
 We have
 $$ \mbox{$g_n \asymp h_n$ and $\rho(g_n) \asymp \rho(h_n)$} $$
 if and only if
 $$ \mbox{there is a compact subset~$C$ of~$A$, such that, for all $n \in
\integer^+$, we have $\mu(g_n) \in \mu(h_n) C$} .$$
 \end{cor}

\begin{proof}
 ($\Rightarrow$) Let $a = \mu(h_n)^{-1} \mu(g_n)$.
 From~\pref{mucalc}, we see that 
 $ \mu(g_n)_{j,j} \asymp \mu(h_n)_{j,j}$
 for $j \in \{1,2\}$, so, using \pref{ajj}, we have
 $$ a_{j,j}
 = \frac{\mu(g_n)_{j,j}}{\mu(h_n)_{j,j}}
 = \begin{cases}
 O(1) & \mbox{if $1 \le j \le 2$;} \\
 1/1 = 1 & \mbox{if $3 \le j \le n$;} \\
 \mu(h_n)_{2,2}/\mu(g_n)_{2,2} =  O(1) & \mbox{if $j = n+1$;} \\
 \mu(h_n)_{1,1}/\mu(g_n)_{1,1} =  O(1) & \mbox{if $j = n+2$} . 
 \end{cases}
 $$
 Therefore $a = O(1)$, as desired.

($\Leftarrow$) Because $C$ is compact, we have
 $$ \mbox{$\mu(g_n) \asymp \mu(h_n)$ and
 $\rho \bigl( \mu(g_n) \bigr) \asymp \rho \bigl( \mu(h_n) \bigr)$} $$
  (cf.\ proof of \pref{h=mu(h)} and \pref{rho(h)=rho(mu)}). Then the desired
conclusions follow from \pref{h=mu(h)} and~\pref{rho(h)=rho(mu)}.
 \end{proof}

\begin{proof}[Proof of Proposition~\ref{bddchange} for $G = \SU(2,n;\F)$]
 Because $C$ is compact, we have
 $g' \asymp g$ and $\rho(g') \asymp \rho(g)$ for any $g' \in CgC$ (cf.\ proof
of \pref{h=mu(h)} and \pref{rho(h)=rho(mu)}).
 Thus, the desired conclusion follows from Corollary~\ref{gn=hn}.
 \end{proof}

Because of Proposition~\ref{mu=rho}, we will often need to calculate $\lVert
h\rVert$ and $\lVert\rho(h)\rVert$. The following observation and its
corollary sometimes simplifies the work, by allowing us to replace $h$
with~$h^{-1}$.

\begin{lem} \label{mu(h-1)}
 We have
 $\mu(h^{-1}) = \mu(h)$
 for $h \in \SU(2,n;\F)$.
 \end{lem}

\begin{proof}
 Define $J$ as in \pref{SUFDefn}, and choose $k_1,k_2 \in K$, such that
$\mu(h) = k_1 h k_2$. For any $a \in A^+$, we see, using \pref{SUF-AN}
or~\pref{ajj}, that $J a^{-1} J = a$, so
 $$ (J k_2^{-1}) h^{-1} (k_1^{-1} J) = J \, \mu(h)^{-1} \, J = \mu(h) .$$
 Note that $\det J = 1$. Also, we have $J^2 = \Id$ and $J^\dagger = J$, so it
is obvious that $J J J^\dagger = J$ and $J J^\dagger = \Id$. Therefore 
 $$ J \in \SU(2,n;\F) \cap \SU(n+2) = K .$$
 Thus, from the definition of~$\mu$, we conclude that 
 $\mu(h^{-1}) = \mu(h)$, as desired.
 \end{proof}

The following corollary is obtained by combining Lemma~\ref{mu(h-1)} with
Corollary~\ref{gn=hn}.

\begin{cor} \label{rho(h-1)}
 We have
 $h^{-1} \asymp h$ and
 $\rho(h^{-1}) \asymp \rho(h)$
 for $h \in \SU(2,n;\F)$.
 \end{cor}

\subsection{The walls of~$A^+$}

\begin{notation}
 For $k \in \{1,2\}$, set
  \begin{equation} \label{Lk-defn}
  L_k = \{\, a \in A^+ \mid a_{2,2} = a_{1,1}^{k-1} \,\} .
  \end{equation}
  From~\pref{A+}, we see that $L_1$ and $L_2$ are the two walls
of~$A^+$. From \pref{calcrho}, we have
 \begin{equation} \label{rho(L)}
  \mbox{$\rho(a) \asymp \lVert a\rVert^k$ for $a \in L_k$} .
  \end{equation}
 \end{notation}

We reproduce the proof of the following result, because it is both
short and instructive. (Although we have no need for it here, let us point
out that the converse of this proposition also holds, and that there is no
need to assume $H \subset AN$.) Because of this proposition (and
Corollary~\ref{CDS->notess}), Section~\ref{SUFlargeSect} will study the
existence of curves $h^t$, such that $h^t \asymp \lVert h^t\rVert^k$, for $k
\in \{1,2\}$. 

\begin{prop}[{(Oh-Witte \cite[Prop.~3.24]{OhWitte-CDS})}]
\label{CDS<>h_m}
 Let $H$ be a closed, connected subgroup of $AN$ in $\SU(2,n;\F)$. If, for
each $k \in \{1,2\}$, there is a continuous curve $h^t$ in~$H$, such that
 $\rho(h^t) \asymp \lVert h^t\rVert^k \to \infty$ as $t \to \infty$, then $H$
is a Cartan-decomposition subgroup.
 \end{prop}

\begin{center}
 \begin{figure}
 \includegraphics[scale=0.5]{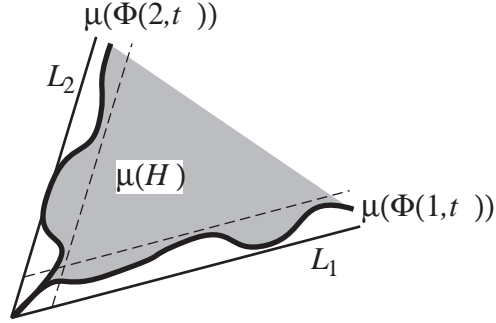}
 \caption{Proposition~\ref{CDS<>h_m}: if $\mu(H)$ contains a curve near each
wall of~$A^+$, then it also contains the interior.}
 \label{mu(H)=wall}
 \end{figure}
 \end{center}

\begin{proof}[Proof {\upshape(cf.\ Figure~\ref{mu(H)=wall} and proof of
Prop.~\ref{SL3-B+})}]
  By hypothesis, there is a continuous, proper map $\Phi\colon \{1,2\}
\times \real ^+ \to H$, such that $\rho \bigl( \Phi(k,t) \bigr) \asymp
\lVert\Phi(k,t)\rVert^k$. Because $H \subset AN$, we know that $H$ is
homeomorphic to some Euclidean space~$\real^m$ \fullsee{solvable}{H=Rn}.

Suppose, for the moment, that $\dim H = 1$. (This will lead to a
contradiction.) We know that $\rho(h) \asymp h$ for $h \in \Phi(1,\real^+)$.
Because $h^{-1} \asymp h$ and $\rho(h^{-1}) \asymp \rho(h)$ \see{rho(h-1)}, we
must also have $\rho(h) \asymp h$ for $h \in \Phi(1,\real^+)^{-1}$. There is
no harm in assuming $\phi(1,0) = \Id$; then
 $\Phi(1,\real^+) \cup \Phi(1,\real^+)^{-1} = H$ (because $\dim H = 1$), so
we conclude that $\rho(h) \asymp h$ for all $h \in H$. This contradicts the
fact that $\rho(h) \asymp \lVert h\rVert^2$ for $h \in \Phi(2,\real^+)$.

We may now assume $\dim H \ge 2$. Then, because $H$ is homeomorphic to
$\real^m$, it is easy to extend $\Phi$ to a continuous and proper map
$\Phi'\colon [1,2] \times \real ^+ \to H$. From \eqref{rho(L)}
and~\pref{gn=hn}, we know that the curve $\mu\bigl(\Phi'(k,t) \bigr)$ stays
within a bounded distance from the wall~$L_k$; say $\operatorname{dist}\bigl[
\bigl(\Phi' (k,t) \bigr), L_k \bigr] < C$ for all~$t$. We may assume $C$ is
large enough that $\operatorname{dist}\bigl( \Phi'(s,1), e \bigr) < C$ for
all $s \in [1,2]$. Then an elementary homotopy argument shows that $\mu \bigl[
\Phi'\bigl( [1,2] \times \real ^+ \bigr) \bigr]$ contains
  $$\{\, a \in A^+ \mid \operatorname{dist}(a, L_1 \cup L_2) > C \, \} ,$$
  so $\mu \bigl[ \Phi'\bigl( [1,2] \times \real ^+ \bigr) \bigr] \approx A^+$.
Because $\mu(H) \supset \mu \bigl[ \Phi'\bigl( [1,2] \times \real ^+ \bigr)
\bigr]$, we conclude from Theorem~\ref{CDSvsmu} that $H$ is a
Cartan-decomposition subgroup.
 \end{proof}

\begin{rem}
 When $\Rrank G = 1$, the Weyl chamber~$A^+$ has only one point at infinity.
Thus, if $H$ is any noncompact subgroup, then the closure of $\mu(H)$ must
contain this point at infinity. This is why it is easy to prove that any
noncompact subgroup of~$G$ is a Cartan-decomposition subgroup
\see{Rrank1-CDS}.

 The idea of Proposition~\ref{CDS<>h_m} is that if $\Rrank G = 2$, then the
points at $\infty$ of the Weyl chamber~$A^+$ form a closed interval. If the
closure of $\mu(H)$ contains the two endpoints of this interval, then, by
continuity, it must also contain all the points in between.

Unfortunately, we have no good substitute for this proposition when $\Rrank G
> 2$. The points at $\infty$ of~$A^+$ form a closed disk (topologically
speaking). It is easy to define a map~$f$ from one disk to another, such that
the image of~$f$ contains the entire boundary sphere, but does not contain the
interior of the disk. Thus, it does not suffice to show only that the
closure of $\mu(H)$ contains the boundary of the disk at $\infty$; rather,
one needs additional homotopical information to guarantee that no interior
points are missed.
 \end{rem}

\begin{lem} \label{mu(SUorSp)}
  Let $G = \SU(2,n; \F)$, and fix some $m \le n/2$.
  Then $\mu \bigl( \SU(1,n;\F) \bigr)$ and $\mu \bigl( \Symp(1,m;\F) \bigr)$ are
the two walls of $A^+$.

 We have
 \begin{enumerate}
 \item \label{mu(SUorSp)-SU}
 $\rho(h) \asymp h$ for $h \in \SU(1,n;\F)$; and
 \item \label{mu(SUorSp)-Sp}
 $\rho(h) \asymp \lVert h\rVert^2$ for $h \in \Symp(1,m;\F)$.
 \end{enumerate} 
 \end{lem}

\begin{proof}
 Let $H = \SU(1,n; \F)$ or $\Symp(1,m;\F)$. Then $H \cap K$ is a maximal
compact subgroup of~$H$. From the Cartan decomposition 
 $$H = (K \cap H) (A \cap H) (K \cap H) ,$$
 and the definition of~$\mu$,
 we conclude that
  $\mu(H) = \mu(A \cap H)$. In the notation of \eqref{Lk-defn}, we see (from
Definition~\ref{SU1nDefn}) that $H \cap A = L_k \cup L_k^{-1}$, where
  $$ k =
  \begin{cases}
  1 & \mbox{if $H = \SU(1,n; \F)$}; \\
  2 & \mbox{if $H = \Symp(1,m; \F)$}.
  \end{cases}
  $$
 Then, since $\mu(a^{-1}) = \mu(a)$ \see{mu(h-1)} and $\mu(a) = a$ for $a \in
A^+$, we conclude that
 $\mu(H) = L_k$ is a wall of~$A^+$. Furthermore, we have
  $\rho(a) \asymp \lVert a\rVert^k$ for $a \in \mu(H)$ \see{rho(L)},
 so
  $\rho(h) \asymp \lVert h\rVert^k$ for $h \in H$ \see{mu=rho}.
 \end{proof}

\begin{cor} \label{SU1inH}
 If there is a continuous curve $h^t \to \infty$ in~$H$, such that $\rho(h^t)
\asymp h^t$, then there is a compact subset~$C$ of~$G$, such that $\SU(1,n;\F)
\subset CHC$.
 \end{cor}

\begin{proof}
 For any (large) $g \in \SU(1,n;\F)$, we see from continuity (more precisely,
from the Intermediate Value Theorem) that there exists $t \in \real^+$, such
that
 $$ \lVert h^t\rVert = \lVert g\rVert .$$
 Then, by assumption and from \fullref{mu(SUorSp)}{SU}, we have
 $$ \rho(h^t) \asymp h^t \asymp g \asymp \rho(g) ,$$
 so there is a compact subset~$C'$ of~$A$,
such that
  $$ \mu \bigl( \SU(1,n;\F) \bigr)  \subset  \{\, \mu(h^t) \mid t \in
\real^+\,\} C' 
 \subset \mu(H) C' $$
 \see{gn=hn}. Therefore
  $$ \SU(1,n;\F) \subset K \, \mu \bigl( \SU(1,n;\F) \bigr) \, K
  \subset K \, \mu(H) C' \, K
  \subset K (KHK) C'K ,$$
 as desired.
 \end{proof}

The following corollary can be proved by a similar argument. (Recall that the
equivalence relation~$\sim$ is defined in \pref{simdefn}.)

\begin{cor} \label{HsimSUorSp}
 Assume $H$ is not compact.
 \begin{enumerate}
 \item \label{HsimSUorSp-SU}
 We have $H \sim \SU(1,n;\F)$ if and only if $\rho(h) \asymp h$ for $h \in H$.
 \item \label{HsimSUorSp-Sp}
 We have $H \sim \Symp(1,m;\F)$ if and only if $\rho(h) \asymp \lVert
h\rVert^2$ for $h \in H$.
 \end{enumerate}
 \end{cor}

Because of Proposition~\ref{CDS<>h_m}, we will often want to show that a
curve $h^t$ satisfies $\rho(h^t) \asymp \lVert h\rVert^k$, for some $k \in
\{1,2\}$. The following lemma does half of the work.

\begin{lem} \label{Owalls}
 Let $X$ be a subset of $\SU(2,n;\F)$.
 \begin{enumerate}
 \item \label{Owalls-linear}
 If $\rho(h) = O(h)$ for $h \in X$, then $\rho(h) \asymp h$ for $h \in X$.
 \item \label{Owalls-square}
 If $\lVert h\rVert^2 = O \bigl( \rho(h) \bigr)$ for $h \in X$, then $\rho(h)
\asymp \lVert h\rVert^2$ for $h \in X$.
 \end{enumerate}
 \end{lem}

\begin{proof}
 From \pref{calcrho} and~\pref{A+}, we have
 $$ \lVert a\rVert = a_{1,1} \le a_{1,1} \, a_{2,2} = \lVert \rho(a) \rVert$$
 and 
 $$ \lVert \rho(a) \rVert = a_{1,1} \, a_{2,2} \le a_{1,1}^2 =
\lVert a\rVert^2 $$
 for $a \in A^+$. Thus, letting $a = \mu(h)$, and using \pref{h=mu(h)}
and~\pref{rho(h)=rho(mu)}, we have:
 $$ \lVert h\rVert \asymp \lVert \mu(h) \rVert \le \lVert \rho \bigl( \mu(h)
\bigr) \rVert \asymp \lVert\rho(h)\rVert $$
 and
 $$ \lVert \rho(h) \rVert \asymp \lVert \rho \bigl( \mu(h) \bigr) \rVert
 \le \lVert \mu(h) \rVert^2 \asymp \lVert h\rVert^2 ,$$
 so $h = O \bigl( \rho(h) \bigr)$ and 
 $\rho(h) = O \bigl( \lVert h\rVert^2 \bigr)$.
 The desired conclusions follow.
 \end{proof}

For convenience, we record the following simple observation. (For the proof,
cf.\ the proof of \pref{h=mu(h)} and \pref{rho(h)=rho(mu)}.)

\begin{lem} \label{conjcurve}
 Let 
 \begin{itemize}
 \item $k \in \{1,2\}$,
 \item $g \in G$, and
 \item $h^t \to \infty$ be a continuous curve in~$H$.
 \end{itemize}
 If $\rho(h^t) \asymp \lVert h^t\rVert^k$, then 
 $\rho(g^{-1} h^t g) \asymp \lVert g^{-1} h^t g \rVert^k$.
 \end{lem}

\subsection{Homogeneous functions of the same degree}

The following well-known, elementary observation is used frequently in the
later sections.

\begin{lem} \label{O(linear)}
 Let $V'$ be a subspace of a finite-dimensional real vector space~$V$, and let
$f_1 \colon V \to W_1$ and $f_2 \colon V \to W_2$ be linear transformations.
 \begin{enumerate}
 \item \label{O(linear)-O}
 If $f_1^{-1}(0) \cap V' = \{0\}$ {\upshape(}or, more generally, if
$f_1^{-1}(0) \cap V' \subset f_2^{-1}(0)${\upshape)}, then there is a linear
transformation $f \colon W_1 \to W_2$, such that $f_2(v) = f \bigl( f_1(v)
\bigr)$ for all $v \in V'$. Therefore $f_2 = O(f_1)$ on~$V'$.
 \item \label{O(linear)-=}
 If $f_1^{-1}(0) \cap V' = f_2^{-1}(0) \cap V'$, then $f_1 \asymp f_2$
on~$V'$.
 \end{enumerate}
 \end{lem}

\begin{proof}
 \pref{O(linear)-O} By passing to a subspace, we may assume $V' = V$. Then,
by modding out $f_1^{-1}(0)$, we may assume $f_1$ is an isomorphism onto its
image. Define $f' \colon f_1(V) \to W_2$ by $f'(w) = f_2 \bigl( f_1^{-1}(w)
\bigr)$, and let $f \colon W_1 \to W_2$ be any extension of~$f'$.

For $v \in V'$, we have
 $$ \lVert f_2(v) \rVert = \lVert f \bigl( f_1(v) \bigr) \rVert
 \le \lVert f\rVert \, \lVert f_1(v) \rVert ,$$
 so $f_2 = O(f_1)$.

\pref{O(linear)-=} From~\pref{O(linear)-O}, we have $f_2 = O(f_1)$ and $f_1 =
O(f_2)$, so $f_1 \asymp f_2$.
 \end{proof}

\begin{eg}
 Let $\Lie H$ be a real Lie subalgebra of $\Lie A + \Lie N$, and assume there
does not exist a nonzero element~$u$ of~$\Lie H$, such that $x_u = 0$ and
$y_u = 0$. Then there exist $\real$-linear transformations $R,S \colon
\F^{n-2} \to \F$, such that $\eta_u = R(x_u) + S(y_u)$ for all $u \in \Lie
H$. (Similarly, $\phi_u$, $\xx_u$, and~$\yy_u$ are also functions of
$(x_u,y_u)$.) Furthermore, we have $u \asymp |x_u| + |y_u|$.
 \end{eg}

The following well-known result is a generalization of the fact that all
norms on a finite-dimensional vector space are equivalent up to a bounded
factor.

\begin{lem}\label{QuadrForms}
 If $V$ is any finite-dimensional real vector space, and $f_1,f_2 \colon V
\to\real$ are two continuous, homogeneous functions of the same degree, such
that $f_1^{-1}(0) = f_2^{-1}(0) = \{0\}$, then $f_1 \asymp f_2$.
 \end{lem}

\begin{proof}
 By continuity, the function $f_1/f_2$ attains a non-zero minimum and a finite
maximum on the unit sphere. Because $f_1/f_2$ is homogeneous of degree zero,
these values bound $f_1/f_2$ on all of $V \smallsetminus \{0\}$.
 \end{proof}

\section{Existence of tessellations} \label{ExistenceSect}

In this section, we show how to construct several families of homogeneous
spaces that have tessellations. All of these examples are based on a method
of T.~Kobayashi \see{construct-tess} that generalizes
Example~\ref{KulkarniEg}. 

\subsection{The general Kulkarni-Kobayashi construction}

As explained in the comments before Theorem~\ref{noncpctdim}, the following
theorem is essentially due to T.~Kobayashi.

\begin{thm}[{(Kobayashi, cf.~\cite[Thm.~4.7]{Kobayashi-properaction})}]
\label{construct-tess}
 If
 \begin{itemize}
 \item $H$ and~$L$ are closed subgroups of~$G$, with only finitely many
connected components;
 \item $L$ acts properly on $G/H$;
 \item $d(L) + d(H) = d(G)$; and
 \item there is a cocompact lattice~$\Gamma$ in~$L$,
 \end{itemize}
 then $G/H$ has a tessellation.
 {\upshape(}Namely, $\Gamma$ is a crystallographic group for $G/H$.{\upshape)}
 \end{thm}

\begin{proof}
 Because $\Gamma$ is a closed subgroup of~$L$, we know that it acts properly
on $G/H$ \see{CHCproper}. Thus, it suffices to show that $\Gamma \backslash
G/H$ is compact.

 From Lemma~\ref{HcanbeAN}, we see that there is no harm in assuming $H
\subset AN$, and that there is a closed, connected subgroup~$L'$ of~$G$, such
that
 \begin{itemize}
 \item $L'$ is conjugate to a subgroup of $AN$,
 \item $d(L') + d(H) = d(G)$, and
 \item $L'C = LC$, for some compact subset~$C$ of~$G$. 
 \end{itemize}
 (Unfortunately, we cannot assume $L \subset AN$: we may not be able to
replace $L$ with~$L'$, because there may not be a cocompact lattice in~$L'$.
For example, there is not lattice in $AN$, because any group with a lattice
must be unimodular \cite[Rem.~1.9, p.~21]{Raghunathan}.)

It suffices to show that $L' \backslash G/H$ is compact. (Because $L' \subset
LC$ is compact, and $\Gamma \backslash L$ is compact, this implies that
$\Gamma \backslash G/H$ is compact, as desired.)

We know that  $L'$ acts properly on $G/H$ \see{CHCproper}, so $L' \times H$
acts properly on $G$, with quotient $L' \backslash G/H$.
Therefore, Lemma~\ref{fiberbundle} implies that $L' \backslash G/H$ has
the same homology as~$G$; in particular,
 $$ \homology_{\dim K}(L' \backslash G/H) \iso \homology_{\dim K}(G) .$$
 From the Iwasawa decomposition $G = KAN$, and because $AN$ is homeomorphic
to $\real^{d(G)}$ \seeand{ANsc}{R=Rd}, we know that $G$ is homeomorphic
to $K \times \real^{d(G)}$. Since $\real^{d(G)}$ is contractible, this implies
that $G$ is homotopy equivalent to~$K$, so $G$ and $K$ have the same
homology; in particular,
 $$ \homology_{\dim K}(G) = \homology_{\dim K}(K) \neq 0 .$$
Since
 \begin{align*}
 \dim(L' \backslash G/H) 
 &= \dim G - \dim L' - \dim H \\
 &= \dim G - \bigl( d(L') + d(H) \bigr) \\
 &= \dim G - d(G) \\
 &= \dim(KAN) - \dim(AN) \\
 &= \dim K ,
 \end{align*}
 this implies that the top-dimensional homology of the manifold $L' \backslash
G/H$ is nontrivial. Therefore $L' \backslash G/H$ is compact
\cite[Cor.~8.3.4]{Dold}, as desired.
 \end{proof}

Our results for $G = \SU(2,2m;\F)$ are based on the following special case of
the theorem.
 The converse of this corollary is proved in Section~\ref{ProofSect}
\see{SUF-known}.

Recall the equivalence relation~$\sim$, introduced in Notation~\ref{simdefn}.

\begin{cor}[{(Kobayashi \cite[Prop.~4.9]{Kobayashi-properaction})}]
\label{SUevenTessExists}
 Let $H$ be a closed, connected subgroup of $G = \SU(2,2m;\F)$.
 If
 \begin{itemize}
 \item $d(H) = 2\df m$; and
 \item either $H \sim \SU(1,2m;\F)$ or $H \sim \Symp(1,m;\F)$,
 \end{itemize}
 then $G/H$ has a tessellation.
 \end{cor}

\begin{proof}
 Let $L_+ = \SU(1,2m;\F)$ and $L_- = \Symp(1,m;\F)$. By assumption, we have $H
\sim L_\varepsilon$, for some $\varepsilon \in \{+,-\}$; let $L =
L_{-\varepsilon}$. Because $\mu(L_+)$ and $\mu(L_-)$ are the two walls of~$A^+$
\see{mu(SUorSp)}, we know that $L = L_{-\varepsilon}$ acts properly
on~$G/L_{\varepsilon}$ \see{proper<>mu(L)}; since $H \sim L_\varepsilon$, this
implies that $L$ acts properly on~$G/H$ \see{CHCproper}.
 Also, we have
 $$d(L) + d(H) = 2\df m + 2\df m = d(G) ,$$
 \seeand{d(Sp)}{d(G)},
 and there is a cocompact lattice in~$L$ \fullcf{classical}{Borel}.
 Thus, the desired conclusion follows from Theorem~\ref{construct-tess}.
 \end{proof}

\subsection{Deformations of $\SO(2,2m)/\SU(1,m)$ and $\SU(2,2m)/\Symp(1,m)$}

The homogeneous spaces described here were found by H.~Oh and D.~Witte
 \cite[Thms.~4.1 and~4.6]{OhWitte-announce}, \cite[Thm.~1.5]{OhWitte-CK}.

\begin{notation} \label{HB-defn}
 For any $\real$-linear $B \colon \F^{n-2} \to \F^{n-2}$, we define
 $$ \Lie H_B =
 \bigset{\begin{pmatrix}
 t & 0 & x  & \eta & \xx \\
 & t & x B & -\xx & -\cjg{\eta} \\
 & & \dots \\
 \end{pmatrix}
 }{
 \begin{matrix}
 t \in \real, \\
 x\in
 \F^{n-2}, \\
 \eta \in \F, \\
 \xx \in \Fim
 \end{matrix}
 }
 \subset \Lie A + \Lie N
 .
 $$
 We write $xB$, rather than $Bx$, because $x$ is a row vector.

 It is easy to see, using, for instance, the formula for the bracket in
\eqref{[u,v]}, that if
 \begin{equation} \label{vBwB}
 \mbox{$\Im \bigl( (vB)(wB)^\dagger \bigr) = -\Im (v w^\dagger)$  for every
$v,w \in \F^{n-2}$} ,
 \end{equation}
 then $\Lie H_B$ is a real Lie subalgebra of
$\Lie A+\Lie N$; we let $H_B$
denote the corresponding connected Lie subgroup of~$AN$.

From \pref{d(H)=dimH}, we have
 \begin{equation} \label{d(HB)}
 \begin{split}
 d(H_B) &= \dim \Lie H_B \\
 &= \dim \real + \dim \F^{n-2} + \dim \F + \dim \Fim \\
 &= 1 + \df(n-2) + \df + (\df -1) \\
 &= \df n .
 \end{split}
 \end{equation}
 \end{notation}

\begin{rem} \label{HB=Sp1m}
 Assume $n = 2m$. By comparing \pref{Sp1mAN} with~\pref{HB-defn}, we see that
there is a $\real$-linear map $B_0 \colon \F^{2m-2} \to \F^{2m-2}$, such that
 $ \Liesp(1,m;\F) \cap (\Lie A + \Lie N) = H_{B_0} $
 (and $B_0$ satisfies~\pref{vBwB}).
 Thus, in general, $H_B$ is a deformation of $\Liesp(1,m;\F) \cap (\Lie A +
\Lie N)$.
 \end{rem}

\begin{thm}[{(Oh-Witte \cite[Thms.~4.1 and~4.6]{OhWitte-announce})}] \label{HBthm}
 Let $B \colon \F^{2m-2} \to \F^{2m-2}$ be $\real$-linear. If
 \begin{itemize}
 \item Condition~\pref{vBwB} holds, and
 \item \label{HBthm-xBnotinFx}
 $xB \notin \F x$, for every nonzero $x \in \F^{2m-2}$,
 \end{itemize}
 then
 \begin{enumerate}
 \item \label{HBthm-mu}
 $\rho(h) \asymp \lVert h\rVert^2$ for $h \in H_B$; and
 \item \label{HBthm-tess}
 $\SU(2,2m;\F)/ H_B$ has a tessellation.
 \end{enumerate}
 \end{thm}

\begin{proof}
  \pref{HBthm-mu} Given $h \in H_B$, write $h = au$, with $a \in A$ and $u
\in N$. We may assume that $a_{1,1} \ge 1$ (by replacing $h$ with~$h^{-1}$ if
necessary \see{rho(h-1)}). It suffices to show $\lVert h\rVert^2 = O \bigl(
a_{1,1} a_{2,2} + |\Delta(h)| \bigr)$ (for then $\lVert h\rVert^2 = O \bigl(
\rho(h) \bigr)$, so Lemma~\fullref{Owalls}{square} applies).

 \setcounter{case}{0}

\begin{case} \label{HB-hinN}
 Assume $a$ is trivial.
 \end{case}
 From~\eqref{exp(phi=0)} and~\pref{HB-defn}, we see that
 $$ h = O \bigl( |x_h|^2 + |\eta_h| + |\xx_h| \bigr) ,$$
 so
 $$ \lVert h\rVert^2 = O \bigl( |x_h|^4 + |\eta_h|^2 + |\xx_h|^2 \bigr) .$$

From \pref{exp(phi=0)} and~\pref{DeltaDefn}, we have
 $$ - \Re \Delta(h)
 = \frac{1}{4} \bigl(|x_h|^2 |y_h|^2 - |x y^\dagger|^2\bigr) 
 + \bigl( |\eta_h|^2 + \xx_h \yy_h \bigr) .$$
 From~\pref{xBnotinFx}, we see that 
 $|x|^2 |xB|^2 - |x (xB)^\dagger|^2 > 0$
 for every nonzero $x \in \F^{2m-2}$, so Lemma~\ref{QuadrForms} implies
 $$ |x_h|^4  \asymp |x_h|^2 |y_h|^2 - |x_h y_h^\dagger|^2 .$$
 Also, because $\yy_h = - \xx_h$ (and $\xx_h \in \Fim$), we have 
 $$|\eta_h|^2 + \xx_h \yy_h = |\eta_h|^2 + |\xx_h|^2 \ge 0 .$$
 Thus,
 $$ - \Re \Delta(h)
 \asymp |x_h|^4 + \bigl( |\eta_h|^2 + |\xx_h|^2 \bigr) ,$$
 so 
 $\lVert h\rVert^2 = O \bigl( \Re \Delta(h) \bigr) = O \bigl( \Delta(h)
\bigr)$, as desired.

\begin{case}
 The general case.
 \end{case}
 From Case~\ref{HB-hinN}, we know
 $\lVert u\rVert^2 = O \bigl( 1 + |\Delta(u)| \bigr)$.
 Then, because
 $\lVert h\rVert \le \lVert a\rVert \lVert u\rVert = a_{1,1} \lVert u\rVert$,
we have
 $$ \lVert h\rVert^2
 \le a_{1,1}^2 \lVert u\rVert^2 
 = O \bigl( a_{1,1}^2 ( 1 + |\Delta(u)|) \bigr)
 = O \bigl( a_{1,1}^2 + a_{1,1}^2 |\Delta(u)| \bigr).$$
 Then, since $a_{1,1} = a_{2,2}$ and $\Delta(h) = a_{1,1} a_{2,2} \Delta(u)$,
we conclude that
 $\lVert h\rVert^2 = O \bigl( a_{1,1}
a_{2,2} + |\Delta(h)| \bigr)$,
 as desired.

\pref{HBthm-tess} From \pref{HBthm-mu} and~\fullref{HsimSUorSp}{Sp},
we see that $H_B \sim \Symp(1,m; \F)$. Then, because $d(H_B) = \df (2m)$
\see{d(HB)}, Theorem~\ref{SUevenTessExists} implies that $\SU(2,2m;\F)/H_B$
has a tessellation.
 \end{proof}

\begin{lem} \label{Bsymplectic}
 Let $B \colon \F^{n-2} \to \F^{n-2}$ be $\real$-linear.
Condition~\pref{vBwB} holds if and only if either
 \begin{enumerate}
 \item $\F = \real$; or
 \item $\F = \complex$ and $B' \in \Symp(2n-4;\real)$, where $x B' = \cjg{xB}$
and we use the natural identification of~$\complex^{n-2}$ with~$\real^{2n-4}$.
 \end{enumerate}
 \end{lem}

\begin{proof}

\setcounter{case}{0}

\begin{case} 
 Assume $\F = \real$.
 \end{case}
 Because $\Im z = 0$ for every $z = 0$, it is obvious that
\pref{vBwB} holds.

\begin{case} 
 Assume $\F = \complex$.
 \end{case}
 If \pref{vBwB} holds, then
 \begin{align*}
 \Im \bigl( (vB')(wB')^\dagger \bigr)
 &= \Im \bigl( (\overline{vB})(\overline{wB})^\dagger \bigr) \\
 &= \Im \overline{\bigl( ({vB})({wB})^\dagger \bigr)} \\
 &= -\Im \bigl( ({vB})({wB})^\dagger \bigr) \\
 &= - \bigl( -\Im (vw^\dagger) \bigr) \\
 &= \Im (vw^\dagger) ,
 \end{align*}
 so $B'$ is symplectic.
 The argument is reversible.
 \end{proof}

\begin{rem} \label{xBnotinFx}
  \begin{itemize}
  \item For $\F = \real$, the assumption that $xB \notin \F x$ simply requires
that $B$ have no real eigenvalues.

  \item For $\F = \complex$, we do not know a good description of the linear
transformations~$B$ that satisfy $xB \notin \F x$, although it is easy to see
that this is an open set (and not dense). A family of examples was
constructed by H.~Oh and D.~Witte (see~\ref{HBeg} below).

 \item  If $n$ is odd, then there does not exist $B \colon \F^{n-2} \to
\F^{n-2}$ satisfying the assumption that $xB \notin \F x$. For $\F = \real$,
this is simply the elementary fact that a linear transformation on an
odd-dimensional real vector space must have a real eigenvalue. For $\F =
\complex$, see Step~\ref{dim(U)<2npf-U/Z} of the proof of
Proposition~\ref{dim(U)<2n}.

 \item If $n$ is even, then, by varying~$B$, one can obtain uncountably  many
pairwise non-conjugate subgroups~$H_B$, such that $\SU(2,n;\F)/H_B$ has a
tessellation. For $\F = \real$, this is proved in
\cite[Thm.~1.5]{OhWitte-CK}). For $\F = \complex$, a similar argument can be
applied to the examples constructed in~\pref{HBeg} below.
  \end{itemize}
 \end{rem}

\begin{eg}[{(Oh-Witte \cite[Thm.~4.6(1)]{OhWitte-announce})}] \label{HBeg}
 Assume $n$~is even, let $B' \in \SO(n-2; \real)$, such that $B'$ has no real
eigenvalue, and define an $\real$-linear map $B \colon \complex^{n-2} \to
\complex^{n-2}$ by $xB = \cjg{x} B'$. Let us verify that $B$ satisfies the
conditions of Theorem~\ref{HBthm} (for $\F = \complex$).

Let $x_1,x_2,y_1,y_2 \in \real^{n-2}$. From the definition of~$B$, and because
$B' \in \SO(n-2; \real)$, we have
 \begin{align*}
 \Im \left( \bigl( (x_1 + i x_2)B \bigr) \bigr( (y_1 + i y_2)B \bigr)^\dagger
\right)
 &= \Im \left( \bigl( (x_1 - i x_2)B' \bigr) \bigr( (y_1 - i y_2) B'
\bigr)^\dagger \right) \\
 &= i \left( (x_1 B') (y_2 B')^\dagger - (x_2 B') (y_1 B')^\dagger \right)
\\
 &= i ( x_1 y_2^\dagger - x_2 y_1^\dagger ) \\
 &= - \Im \bigl( (x_1 +i x_2) (y_1 + i y_2)^\dagger \bigr)
 .
 \end{align*}

Suppose $Bx = \lambda x$, for some $\lambda \in \complex$. Because $B \in
\SO(2n-4;\real)$, we must have $|\lambda| = 1$. Then
 $$ B'(x + \cjg{\lambda x})
 = B' x + \cjg{\lambda} B' \cjg{x}
 = \cjg{B' \cjg{x}} + \cjg{\lambda} Bx
 = \cjg{Bx} + \cjg{\lambda} (\lambda x)
 = \cjg{\lambda x} + x .$$
 Because $B'$ has no real eigenvalues, we know that $1$ is not an eigenvalue
of~$B'$, so we conclude that $x + \cjg{\lambda x} = 0$. Similarly,
because $-1$ is not an eigenvalue of~$B'$, we see that $x - \cjg{\lambda x} =
0$. Therefore
 $$x = \frac{1}{2} \bigl( (x + \cjg{\lambda x}) + (x - \cjg{\lambda x}) \bigr)
= \frac{1}{2} ( 0 + 0 ) = 0 .$$
 \end{eg}

\subsection{Deformations of $\SU(2,2m)/\SU(1,2m)$}

These examples are new for $\F = \complex$, but provide nothing interesting
for $\F = \real$ \fullsee{Hc}{real}.

\begin{notation} \label{SUegsDefn}
 For $c \in (0,1]$, we define
 $$ \LieHc =
 \bigset{
 \begin{pmatrix}
 t & \phi & x & \Re \phi + c \Im \phi & \xx \\
  & 0 & 0 & 0 & * \\
 & & \dots \\
 \end{pmatrix}
 }{
 \begin{matrix}
 t \in \real, \\
 \phi \in \F, \\
 x \in \F^{n-2}, \\
 \xx \in\Fim
 \end{matrix}
 }
 .$$

It is easy to see, using, for instance, the formula for the bracket in
\eqref{[u,v]}, that $\LieHc$ is a real Lie subalgebra
of $\Lie A+\Lie N$ (even without the assumption that $0 < c \le 1$); we
let $\Hc$ be the corresponding connected Lie subgroup of $AN$.

From \pref{d(H)=dimH}, we have
 \begin{equation} \label{d(Hc)}
 \begin{split}
 d(\Hc) 
 &= \dim \LieHc \\
 &= \dim \real + \dim \F + \dim \F^{n-2} + \dim \Fim \\
 &= 1 + \df + \df(n-2) + (\df -1) \\
 &= \df n .
 \end{split}
 \end{equation}
 \end{notation}

\begin{rem} \label{Hc}
 Let $\su(1,n;\F)$ be embedded into $\su(2,n;\F)$ as in~\ref{SU1nAN}. 
  \begin{enumerate}
  \item \label{Hc-real}
 If $\F=\real$, then $c$ is irrelevant in the definition of $\LieHc$
(because $\Im \phi = 0$); therefore
 $\LieHc=\su(1,n;\real)\cap (\Lie A + \Lie N)$.
  \item If $\F=\complex$, then
 $\LieHcc1 = \su(1,n;\complex)\cap (\Lie A + \Lie N)$.
  \end{enumerate}
 Thus, in general, $\LieHc$ is either $\su(1,n;\F)\cap (\Lie A + \Lie N)$
or a deformation of it.
 \end{rem}

\begin{thm} \label{SUegs}
 Assume $\F = \complex$, and $n = 2m$ is even.
 If $c \in (0,1]$, then
 \begin{enumerate}
 \item \label{SUegs-linear}
 $\rho(h) \asymp h$ for $h \in \Hc$; and
 \item \label{SUegs-tess}
 $\SU(2,2m;\F)/\Hc$ has a tessellation.
 \end{enumerate}
 \end{thm}

\begin{proof}
 \pref{SUegs-linear}  Given $h \in \Hc$, it suffices to
show that $\rho(h) = O(h)$ \fullsee{Owalls}{linear}. Write $h = au$,
with $a \in A$ and $u \in N$. We may assume that $a_{1,1} \ge 1$ (by replacing
$h$ with~$h^{-1}$ if necessary \see{rho(h-1)}). 

  Let $Q \colon \complex \oplus \complex^{n-2} \oplus \complex \to \real$ be
the real quadratic form 
 $$Q(\phi,x,\eta)=|x|^2+2\Re(\phi \cjg{\eta}) ,$$
 and let $V$ be the $\real$-subspace of $\complex \oplus \complex^{n-2} 
\oplus \complex$ defined by
 $$
 V =
 \bigset{
  ( \phi,x,\eta)
  }{
  \begin{matrix}
 \phi \in \complex , \\
  x \in \complex^{n-2}, \\
  \eta = \Re\phi + c \Im\phi
  \end{matrix}
  } . $$

\setcounter{step}{0}

\begin{step} \label{SUegspf-Q=x2+phi2}
 For $v \in V$, we have $Q(v) \asymp |\phi|^2 + |x|^2$.
 \end{step}
 For $(\phi,x,\eta) \in V \smallsetminus \{0\}$, we have
 \begin{align*}
 Q(\phi,x,\eta)
 &= |x|^2+2\Re(\phi \cjg{\eta}) \\
 &= |x|^2 + 2 \Re \bigl( \phi (\cjg{\Re\phi + c \Im\phi}) \bigr)  \\
 &= |x|^2 + 2 (\Re\phi)^2 - 2 c (\Im\phi)^2 \\
 &> 0
 \end{align*}
 (because $c>0$ and $\Im \phi$ is purely imaginary). Thus, the restriction
of~$Q$ to~$V$ is positive definite, so the desired conclusion follows
from Lemma~\ref{QuadrForms}.

\begin{step} \label{SUegspf-u=x2+phi2+xx}
 We have
 $u_{1,n+2} \asymp \bigl( |\phi_u|^2 + |x_u|^2 \bigr) + |\xx_u|$.
 \end{step}
 From \pref{exp(y=0)} (with $\yy = 0$), we have
 $$ \Re u_{1,n+2}
 = - \left( \frac{1}{2}
|x_u|^2+ \Re(\phi_u \cjg{\eta_u}) \right)
  \asymp |x_u|^2+2\Re(\phi_u \cjg{\eta_u}) $$
 and
 $$ \Im u_{1,n+2}
 =  \xx_u .$$
 Then, from Step~\ref{SUegspf-Q=x2+phi2}, we see that $\Re u_{1,n+2} \asymp
|\phi_u|^2 + |x_u|^2$,
 so 
 $$ u_{1,n+2} \asymp |\Re u_{1,n+2}| + |\Im u_{1,n+2}|
 \asymp \bigl( |\phi_u|^2 + |x_u|^2 \bigr) + |\xx_u| ,$$
 as desired.

\begin{step}
 Completion of the proof.
 \end{step}
 From Step~\ref{SUegspf-u=x2+phi2+xx}, we have 
 $$ h_{1,n+2} = a_{1,1} u_{1,n+2}
  \asymp  a_{1,1}  \bigl( |x_u| + |\phi_u| \bigr)^2 + a_{1,1} |\xx_u|.$$
  Also, from \pref{exp(y=0)}, we have
  $$ h_{jk} =
  \begin{cases}
  O(1) & \mbox{if $j \neq 1$ and $k \neq n+2$} \\
  O \bigl( a_1 (|\phi_u| + |x_u|) \bigr) & \mbox{if $j = 1$ and $k \neq n+2$}
\\
  O \bigl( |\phi_u| + |x_u| \bigr) & \mbox{if $j \neq 1$ and $k = n+2$} .
  \end{cases}
  $$
  Thus, it is easy to see that
  $$ \rho(h) = O \bigl( a_{1,1} |\xx_u| + a_{1,1}  (|\phi_u| + |x_u|)^2 \bigr)
  = O(h_{1,n+2})
  = O(h) ,$$
 so the desired conclusion follows from Lemma~\fullref{Owalls}{linear}.

\pref{SUegs-tess} From \pref{SUegs-linear} and
\fullref{HsimSUorSp}{SU}, we see that $\Hc \sim
\SU(1,n)$. Then, because $d(\Hc) = 2n$ \see{d(Hc)},
Theorem~\ref{SUevenTessExists} implies that
$\SU(2,2m;\F)/\Hc$ has a tessellation.
 \end{proof}

\begin{rem}
 Proposition~\ref{HcUncountable} shows that if $\F = \complex$, then $\Hc$ is
not conjugate to $\Hcc{c'}$ unless $c = c'$ (for $c,c' \in (0,1]$). Thus,
Theorem~\fullref{SUegs}{tess} implies that, by varying $c$, one obtains
uncountably many nonconjugate subgroups $\Hc$, such that $\SU(2,2m)/\Hc$ has
a tessellation.
 \end{rem}

\subsection{The product of two rank-one groups}

\begin{prop} \label{G1xG2-tess}
 Let $G = G_1 \times G_2$ be the direct product of two connected, linear,
almost simple Lie groups $G_1$ and~$G_2$ of real rank one, with finite center,
and let $H$ be a nontrivial, closed, connected, proper subgroup of~$AN$.

The homogeneous space $G/H$ has a tessellation if and only if, perhaps after
interchanging $G_1$ and~$G_2$, there is a continuous homomorphism
$\sigma \colon
AN \cap G_1 \to AN \cap G_2$, such that
 $$ H = \{\, \bigl(h , \sigma(h) \bigr) \mid h \in AN \cap G_1 \,\} .$$
 \end{prop}

\begin{proof}
 ($\Rightarrow$) We may assume $d(G_1) \ge d(G_2)$ (by interchanging $G_1$
and $G_2$ if necessary).

\setcounter{case}{0}

\begin{case} \label{Hcapboth}
 Assume $H \cap G_1 \neq e$ and $H \cap G_2 \neq e$.
 \end{case}
 For $j = 1,2$, we know that $H \cap G_j$ is not
compact \fullsee{solvable}{nocpct}, so Corollary~\ref{Rrank1-CDS} implies that
there is a compact subset~$C_j$ of~$G_j$, such that  $C_j (H \cap G_j)C_j =
G_j$. Then, letting $C = C_1 C_2$, we have $CHC = G$, so
Proposition~\ref{CDS->notess} implies that $G/H$ does not have a
tessellation. This is a contradiction.

\begin{case}
 Assume $H \cap G_1 \neq e$ and $H \cap G_2 = e$.
 \end{case}
 From Corollaries~\ref{Rrank1-CDS} and~\ref{CDSvsmu}, we know that there is a
compact subset~$C$ of~$A \cap G_1$, such that $\mu(G_1) \subset \mu(H) C$.
Therefore, Corollary~\ref{noncpct-dim-notess} (with $G_1$ in the place
of~$H_1$) implies $d(H) \ge d(G_1) = \dim(G_1 \cap AN)$.
 Then, because $H \cap G_2 = e$ (and $H \subset AN$), we conclude that $H$ is
the graph of a homomorphism from $G_1 \cap AN$ to $G_2 \cap AN$, as desired.

\begin{case}
 Assume $H \cap G_1 = e$.
 \end{case}
 From Corollary~\ref{tess->dim>1,2}, we know that $\dim H \ge d(G_2)$. Then,
since $H \cap G_1 = e$, we conclude that $H$ is the graph of a homomorphism
from $G_2 \cap AN$ to $G_1 \cap AN$. Interchanging $G_1$ and~$G_2$ yields the
desired conclusion.

\medskip

($\Leftarrow$)
 We verify the hypotheses of Theorem~\ref{construct-tess}, with $G_2$ in the
role of~$L$.

 Let $\overline{H}$ be the image of~$H$ under the natural homomorphism $G \to
G/G_2$. Because $H \subset AN$, we know that $\overline{H}$ is closed
\fullsee{solvable}{H=Rn}. It is well known (and follows easily from
\pref{proper<>CHC}) that any closed subgroup acts properly on the ambient
group, so this implies that $\overline{H}$ acts properly on $G/G_2$. 
 From the definition of~$H$, we have $H \cap G_2 = e$, so we conclude that $H
\iso \overline{H}$ acts properly on $G/G_2$; equivalently, $G_2$ acts properly
on $G/H$ \cf{proper<>CHC}.

 Because $AN = (AN \cap G_1) \times (AN \cap G_2)$, we have
  $d(G) = d(G_1) + d(G_2)$.
 Also, we have $d(H) = \dim H$ \see{d(H)=dimH} and, from the definition
of~$H$, we have 
 $\dim H = \dim(AN \cap G_1) = d(G_1)$.
 Therefore
 $$ d(H) + d(G_2) = d(G_1) + d(G_2) = d(G) .$$

There is a cocompact lattice in~$G_2$ \fullcf{classical}{Borel}.

So Theorem~\ref{construct-tess} implies that $G/H$ has a tessellation.
 \end{proof}

\subsection{T.~Kobayashi's examples of higher real rank}

T.~Kobayashi observed that, besides the examples with $G = \SO(2,2n)$ or
$\SU(2,2n)$ \see{KulkarniEg}, Theorem~\ref{construct-tess} can also be used
to construct tessellations of some homogeneous spaces $G/H$ in which $G$
and~$H$ are simple Lie groups with $\Rrank G > 2$. He found one pair of
infinite families, and several isolated examples. 

\begin{thm}[{(Kobayashi \cite[Cor.~5.6]{Kobayashi-survey97})}]
\label{Kobayashibigeg}
 Each of the following homogeneous spaces has a tessellation:
 \begin{enumerate}
 \item \label{Kobayashibigeg-SO4/Sp1}
 $\SO(4,4n)/\Symp(1,n)$;
 \item $\SO(4,4n)/\SO(3,4n)$;
 \medskip
 \item $\SO(8,8)/\SO(8,7)$;
 \item $\SO(8,8)/\operatorname{Spin}(8,1)$;
 \medskip
 \item $\SO(4,4)/\SO(4,1)$;
 \item $\SO(4,4)/\operatorname{Spin}(4,3)$;
 \medskip
 \item $\SO(4,3)/\SO(4,1)$;
 \item $\SO(4,3)/G_{2(2)}$.
 \end{enumerate}
 \end{thm}

 It would be very interesting to find other examples of simple Lie groups~$G$
with reductive subgroups $H$ and~$L$ that satisfy the hypotheses of
Theorem~\ref{construct-tess}.

\begin{rem}
 Let $G = \SO(4,4n)$ and $H' = \Symp(1,n) \cap AN$.
From~\fullref{Kobayashibigeg}{SO4/Sp1}, we know that $G/H$ has a
tessellation. H.~Oh and D.~Witte \cite[Thm.~4.6(2)]{OhWitte-announce} pointed
out that the deformations $G/H_B$ (where $H_B$ is as in Theorem~\ref{HBthm},
with $\F = \complex$) also have tessellations, but it is not known whether
there are other deformations of $G/H'$ that also have tessellations. 

 It does not seem to be known whether the other examples in
Theorem~\ref{Kobayashibigeg} lead to nontrivial deformations, after
intersecting $H$ with $AN$.
 \end{rem}

\section{Large subgroups of $\SO(2,n)$ and $\SU(2,n)$} \label{SUFlargeSect}

This section presents a short proof of the results we need from
\cite{OhWitte-CDS} and \cite{IozziWitte-CDS}.  Those papers provide an
approximate calculation of $\mu(H)$, for every closed, connected subgroup~$H$
of $\SO(2,n)$ or $\SU(2,n)$, respectively, but here we consider only
subgroups of large dimension. Also, we do not need a complete description of
the entire set $\mu(H)$; we are only interested in whether or not there is a
curve~$h^t$, such that $\rho(h^t) \asymp \lVert h^t\rVert^k$, for some $k \in
\{1,2\}$. The main results of this section are Theorem~\ref{maxnolinear} (for
$k = 1$) and Theorem~\ref{bestnosquare} (for $k = 2$). They give a sharp
upper bound on $d(H)$, for subgroups~$H$ that fail to contain such a curve,
and, if $n$ is even, also provide a fairly explicit description of all the
subgroups of $AN$ that attain the bound.

 Because of the limited scope of this section, the proof
here is shorter than the previous work, and we are able to give a fairly
unified treatment of the two groups $\SO(2,n)$ and $\SU(2,n)$. The arguments
are elementary, but they involve case-by-case analysis and a lot of
details, so they are not pleasant to read. 

\begin{standing} \label{StandingSU2F}
 Throughout this section:
 \begin{enumerate}
 \item We use the notation of \S\ref{coordsSect}. (In particular, $\F = \real$
or~$\complex$, and $\df = \dimR \F$.)
 \item $G = \SU(2,n;\F)$.
 \item $n \ge 3$.
 \item $H$ is a closed, connected subgroup of $AN$ that is compatible
with~$A$ \see{compatibleDefn}, so $\dim H = d(H)$ \see{d(H)=dimH}.
 \item $U = H \cap N$. (Note that $U$ is connected \fullsee{solvable}{HcapL}.)
 \item $\Lie U_{\phi=0} = \{\, u \in \Lie U \mid \phi_u = 0 \,\}$.
 \item We use the notation of \S\ref{CalcSect}. (In particular,
$\lVert\rho(h)\rVert$ is defined in \pref{rhoDefn} and $\Delta(h)$ is defined
in~\pref{DeltaDefn}.)
 \item Except in Subsection~\ref{compatibleSubsect}, $H$~is
compatible with~$A$ \see{compatibleDefn}.
 \end{enumerate}
 \end{standing}

\subsection{Subgroups compatible with~$A$} \label{compatibleSubsect}

Recall that the Real Jordan Decomposition of an element of~$G$ is defined
in \pref{JordanDecomp}; any element~$g$ of $AN$ has a Real Jordan
Decomposition $g = au$ (with $c$ trivial). If $a$ is an element of~$A$,
rather than only conjugate to an element of~$A$, we could say that $g$ is
``compatible with~$A$." We now define a similar, useful notion for subgroups
of $AN$. Lemma~\ref{conjtocompatible} shows there is usually no loss of
generality in assuming that $H$ is compatible with~$A$, and
Lemma~\ref{not-semi} shows that the compatible subgroups can be described
fairly explicitly.

\begin{defn}[{\cite[Defn.~2.2]{OhWitte-CDS}}] \label{compatibleDefn}
 Let us say that $H$ is \emph{compatible} with~$A$ if $H
\subset T U C_N(T)$, where $T = A \cap (HN)$, $U = H \cap N$, and $C_N(T)$
denotes the centralizer of~$T$ in~$N$.
 \end{defn}

In preparation for the proofs of the main results, let us state a lemma that
records a few of the nice properties of Jordan components.

\begin{lem}[{(cf.\ \cite[Lem.~15.3, p.~99]{Humphreys-Algic})}]
\label{aucgood}
 Let $g = auc$ be the Real Jordan Decomposition of an element~$g$ of~$G$. Then
 \begin{enumerate}
 \item
 the Real Jordan Decomposition of $\Ad g$ is $\Ad g = (\Ad a) (\Ad u) (\Ad
c)$; and
 \item $a$, $u$, and~$c$ all belong to the Zariski closure of $\langle g
\rangle$.
 \end{enumerate}
 Therefore:
 \begin{enumerate} \renewcommand{\theenumi}{\alph{enumi}}
 \item \label{aucgood-norm}
 $a$, $u$, and~$c$ each normalize any connected
subgroup of~$G$ that is normalized by~$g$; and
 \item \label{aucgood-id}
 if $(\Ad g)v \in v + W$, for some $v \in \Lie G$ and
some $(\Ad g)$-invariant subspace~$W$ of~$\Lie G$, then 
 $$ \mbox{$(\Ad a)v$, $(\Ad u)v$, and $(\Ad c)v$ all belong to $v + W$} .$$
 \end{enumerate}
 \end{lem}

\begin{lem}[{\cite[Lem.~2.3]{OhWitte-CDS}}] \label{conjtocompatible}
 $H$ is conjugate, via an element of~$N$, to a subgroup that is compatible
with~$A$.
 \end{lem}

\begin{proof}[Idea of proof]

\setcounter{case}{0}

 \begin{case} \label{conjtocompatiblePf-algebraic}
 Assume, for the Real Jordan Decomposition $h = au$ of each element~$h$
of~$H$, that $a$ and~$u$ belong to~$H$.
 \end{case}
 Let $T$ be a maximal split torus of~$H$. (Recall that a split torus is a
subgroup consisting entirely of hyperbolic elements.) Then $T$ is
contained in some maximal split torus of~$G$, that is, in some subgroup
of~$G$ conjugate to~$A$; replacing $H$ by a conjugate, we may assume $T
\subset A$. In other words, we now know that $H \cap A$ is a maximal split
torus of~$H$.

Given $h \in H$, we have the Real Jordan Decomposition $h = au$. By
assumption, $a \in H$; thus, $a$ belongs to some maximal split torus~$T'$
of~$H$. A fundamental result of the theory of solvable algebraic groups
implies that all maximal split tori of~$H$ are conjugate via an
element of $H \cap N$ \cite[Thm.~4.21]{BorelTits-Reductive}, so there is some
$x \in H \cap N$, such that $x^{-1} a x \in A$. Then $\langle T, x^{-1} a x
\rangle$, being a subgroup of~$A$, is a split torus. Thus, the maximality
of~$T$ implies that $x^{-1} a x \in T$; let $t = x^{-1} a x$. Then 
 $$ h = au = x t x^{-1} u = t (t^{-1} x t) x^{-1} u \in T (H \cap N) .$$
 Since $h \in H$ is arbitrary, we conclude that 
 \begin{equation} \label{conjtocompatiblePf-H=TU}
 H = T (H \cap N) ,
 \end{equation}
 so $H$ is
compatible with~$A$.

 \begin{case}
 The general case.
 \end{case}
 Let 
 $$ \overline{H}
 =  \langle\, a, u \mid \mbox{$au = ua \in H$, $a$~hyperbolic, $u$~unipotent}
\,\rangle $$
 be the subgroup of~$AN$ generated by the Jordan Components of the elements
of~$H$. (Of course, since every element of~$H$ has a Jordan Decomposition, we
have $H \subset \overline{H}$.)
 Then Case~\ref{conjtocompatiblePf-algebraic} applies to~$\overline{H}$, so,
replacing $H$ by a conjugate, we may assume $\overline{H} = \overline{T} \,
\overline{U}$, where $\overline{T} = \overline{H} \cap A$ and $\overline{U} =
\overline{H} \cap N$ \see{conjtocompatiblePf-H=TU}.

Because $\overline{H}$ normalizes~$H$ \fullsee{aucgood}{norm}, we know that
$\Ad_G h$ acts as the identity on $\overline{\Lie H}/\Lie H$, for all $h \in
H$. Hence, Lemma~\fullref{aucgood}{id} implies that $\Ad_G \overline{h}$ acts
as the identity on $\overline{\Lie H}/\Lie H$, for all $\overline{h} \in
\overline{H}$; therefore $[\overline{H}, \overline{H}] \subset H$. 
 Also, we have $[\overline{H},\overline{H}] \subset [AN,AN] \subset N$. Thus,
letting $U = H \cap N$, we have 
 \begin{equation} \label{conjtocompatiblePf-comm}
 [\overline{H},\overline{H}] \subset H \cap N = U .
 \end{equation}
 Because $\overline{T} \subset A$ and $\overline{\Lie
U}$ is $(\Ad_G(T))$-invariant, the adjoint action
of~$\overline{T}$ on~$\overline{\Lie U}$ is completely reducible, so
\pref{conjtocompatiblePf-comm} implies that there is a subspace $\Lie C$
of~$\overline{\Lie U}$, such that $[\overline{T}, \Lie C] = 0$ and $\Lie U +
\Lie C = \overline{\Lie U}$. Therefore, $U \, C_{\overline{U}}(\overline{T}) =
\overline{U}$, so
 \begin{equation} \label{conjtocompatiblePf-Hbarcompat}
 \overline{H}
 = \overline{T} \, \overline{U}
 = \overline{T} U C_{\overline{U}}(\overline{T})
 \subset \overline{T} U C_N(\overline{T}) .
 \end{equation}

 Let $\pi \colon AN \to A$ be the projection with kernel~$N$, and let $T =
\pi(H)$. Then 
 $$T = \pi(H) \subset \pi(\overline{H}) = \overline{T} ,$$
 so $C_N(T) \supset C_N(\overline{T})$. For any $h \in H$, we know,
from~\pref{conjtocompatiblePf-Hbarcompat}, that there exist $t \in
\overline{T}$, $u \in U$ and $c \in C_N(\overline{T})$, such that $h = tuc$.
Because $uc \in N$, we must have $t = \pi(h) \in T$ and, because $C_N(T)
\supset C_N(\overline{T})$, we have $c \in C_N(T)$. Therefore, $h \in T U
C_N(T)$. We conclude that $H \subset T U C_N(T)$, so $H$ is compatible
with~$A$.
 \end{proof}

The preceding proposition shows that $H$ is conjugate to a subgroup~$H'$ that
is compatible with~$A$. The subgroup~$H'$ is usually not unique, however. The
following lemma provides one way to change~$H'$, often to an even better
subgroup.

\begin{lem} \label{conjUomega}
 Assume that $H$ is compatible with~$A$, and let $T = A \cap (HN)$. If $u \in
C_N(T)$, then $u^{-1} H u$ is compatible with~$A$.
 \end{lem}

\begin{proof}
 Let $H' = u^{-1} H u$, $T' = A \cap (H'N)$, and $U' = H' \cap N$. Because
$u$ centralizes~$T$, we have 
 $$u^{-1} T u = T .$$
 Also, because $u \in N$, and $N$ is normal, we have $u^{-1}HuN = HN$, so
 $$ u^{-1} T u = T = A \cap (HN) = A \cap (u^{-1}HuN) = A \cap (H'N) = T' .$$
 Since $u \in N$, we have $u^{-1} N u = N$, so
 $$u^{-1} U u = u^{-1} (H \cap N) u
 =  (u^{-1} H u) \cap (u^{-1} N u)
 = H' \cap N
 = U' ,$$
 and
 $$ u^{-1} C_N(T) u = C_{u^{-1}Nu}(u^{-1} T u) = C_N(T') .$$
 Thus,
 $$H' = u^{-1} H u
 \subset  u^{-1} T U C_N(T) u
 =  \bigl(u^{-1} T u \bigr) \bigl( u^{-1} U u \bigr) \bigl( u^{-1} C_N(T) u
\bigr)
 = T' U' C_N(T'),$$
 as desired.
 \end{proof}

\begin{lem}[{\cite[Lem.~2.4]{OhWitte-CDS}}] \label{not-semi}
 If $H$ is compatible with~$A$, then either
 \begin{enumerate}
 \item \label{not-semi-TU}
 $H = (H \cap A) \ltimes (H \cap N)$; or
 \item \label{not-semi-not}
 there is a positive root~$\omega$, a nontrivial group homomorphism
$\psi\colon \ker \omega \to \rsg_\omega \rsg_{2\omega}$, and a closed, connected
subgroup~$U$ of~$N$, such that
 \begin{enumerate}
 \item \label{not-semi-codim1}
 $H = \{\, a \, \psi(a) \mid a \in \ker \omega \,\} \, U$;
 \item \label{not-semi-normal}
 $U$ is normalized by both $\ker\omega$ and~$\psi( \ker\omega)$; and
 \item \label{not-semi-disjoint}
 $U \cap \psi( \ker\omega) = e$.
 \end{enumerate}
 \end{enumerate}
 \end{lem}

\begin{proof} 
 Because $H$ is compatible with~$A$, we have $H \subset T U C_N(T)$, where $T
= A \cap (HN)$ and $U = H \cap N$. We may assume that $H \neq TU$, for
otherwise \pref{not-semi-TU} holds. Therefore $C_N(T) \neq e$. Because $\Lie
n$ is a sum of root spaces, this implies that there is a positive
root~$\omega$, such that $T \subset \ker\omega$. Because $\Rrank G = 2$, we
have $\dim (\ker\omega) = 1$, so we must have $T =\ker\omega$ (otherwise
we would have $T = e$, so $H = U = TU$; hence \pref{not-semi-TU} holds).
Therefore, $C_N(T) = U_\omega U_{2\omega}$. 

Because $U \subset H \subset T U C_N(T)$, we have $H = U \bigl[H \cap \bigl(T
C_N(T) \bigr) \bigr]$, so there is a nontrivial one-parameter subgroup
$\{x^t\}$ in $H \cap \bigl(T C_N(T) \bigr)$ that is not contained in~$U$.
Because $T$ centralizes $C_N(T)$, we may write $x^t = a^t u^t$ where
$\{a^t\}$ is a one-parameter subgroup of~$T$ and $\{u^t\}$ is a one-parameter
subgroup of $C_N(T)$. Furthermore, this decomposition is unique, because $T
\cap C_N(T) = e$. (In fact, $x^t = a^t u^t$ is the Real Jordan Decomposition
of~$x^t$.) Define $\psi \colon \ker\omega \to U_\omega U_{2\omega}$ by
$\psi(a^t) = u^t$ for all $t \in \real$. 

\pref{not-semi-codim1} For all $t \in \real$, we have $a^t \psi(a^t) = a^t
u^t = x^t \in H$, which establishes one inclusion of~\pref{not-semi-codim1}.
The other will follow if we show that $\dim H - \dim U = 1$, so suppose $\dim
H - \dim U \ge 2$. Then Lemma~\fullref{dimT}{A} implies that $A \subset H$,
so it follows from Lemma~\ref{rootdecomp} (with $T = A$ and $\omega = 0$)
that $H = A \ltimes (H \cap N)$, contradicting our assumption that $H \neq
TU$.

 \pref{not-semi-normal} Because $x^t \in H$, we know that each of $a^t$
and~$u^t$ normalizes~$H$ \fullsee{aucgood}{norm}. Being in $AN$, they also
normalize~$N$. Therefore, they normalize $H \cap N = U$.

\pref{not-semi-disjoint} Suppose $U \cap \psi(\ker\omega) \neq e$. Because
the intersection $U \cap \psi(\ker\omega)$ is connected
(see~\ref{not-semi-normal} and~\fullref{solvable}{HcapL}), and
$\dim(\ker\omega) = 1$, we must have $\psi(\ker\omega) \subset U$. Therefore
$a^t = x^t u^{-t} \in HU = H$, so $T \subset H$. This contradicts our
assumption that $H \neq TU$.
 \end{proof}

\begin{cor} \label{TnormsU}
 If $H$ is compatible with~$A$, then $A \cap (HN)$ normalizes~$H \cap N$.
 \end{cor}

\begin{lem}[{\cite[Lem.~2.8]{OhWitte-CDS}}] \label{dimT}
 If $\dim \bigl( H/(H \cap N) \bigr) \ge \dim A$, then
 \begin{enumerate}
 \item \label{dimT-A}
 $H$ contains a conjugate of~$A$; and
 \item \label{dimT-CDS}
 $H$ is a Cartan-decomposition subgroup.
 \end{enumerate}
 \end{lem}

\begin{proof}
 \pref{dimT-A} Let $\pi \colon AN \to A$ be the projection with kernel~$N$,
and let $\overline{H}$ be the Zariski closure of~$H$. From the structure
theory of solvable algebraic groups \cite[Thm.~10.6(4),
pp.~137--138]{Borel-Algic}, we know that $\overline{H} = T \ltimes U$ is the
semidirect product of a torus~$T$ and and unipotent subgroup $U \subset N$.
Replacing $H$ by a conjugate under~$N$, we may assume that $T \subset A$.
Since
 $$ \dim A \le \dim \bigl( H/(H \cap N) \bigr) = \dim\bigl(\pi(H) \bigr) \le
\dim A ,$$
 we must have $\pi(H) = A$, so 
 $$A = \pi(H) \subset \pi(\overline{H}) = \pi(TU) = T \subset
\overline{H} $$
 normalizes~$H$ \see{Zar-norm}. Then, since $\pi(H) = A$, we conclude that $A
\subset H$ \see{rootdecomp}.

\pref{dimT-CDS} From~\pref{dimT-A}, we see that, by replacing $H$ with
a conjugate subgroup, we may assume $A \subset H$. Because $A$ is a
Cartan-decomposition subgroup \see{AisCDS}, this implies $H$ is a
Cartan-decomposition subgroup.
 \end{proof}

The following basic result was used twice in the
above arguments.

\begin{lem}[{(cf.~\cite[pf.~of Thm.~3.2.5, p.~42]{ZimmerBook})}]
\label{Zar-norm}
 If $H$ is a closed, connected subgroup of~$G$, then the Zariski closure
of~$H$ normalizes~$H$.
 \end{lem}

\subsection{Subgroups with no nearly linear curve}
 Our goal is to prove Theorem~\ref{maxnolinear}; we begin with some
preliminary results.

First, an observation that simplifies the calculations in some cases, by
allowing us to assume that $x_u = 0$.

\begin{lem} \label{x=0}
 Let $u \in \Lie U$. If $\dimF (\F x_u + \F y_u) \le 1$ and $y_u \neq 0$, then
there is some $g \in \rsg_\alpha$, such that 
 \begin{enumerate}
 \item $x_{g^{-1} u g} = 0$,
 \item $\phi_{g^{-1} u g} =  \phi_u$, and
 \item $y_{g^{-1} u g} = y_u$.
 \end{enumerate}
 \end{lem}

\begin{proof}
 Because $\dimF (\F x_u + \F y_u) \le 1$ and $y_u \neq 0$, there is some
$\lambda \in \F$, such that $x_u = \lambda y_u$. Let 
 \begin{itemize}
 \item $v$ be the element of~$\rsp_\alpha$ with $\phi_\alpha = -\lambda$, 
 \item $g = \exp(v) \in \rsg_\alpha$, and 
 \item $w = g^{-1} u g$.
 \end{itemize}
 From~\pref{conjugation}, we see that 
 \begin{itemize}
 \item $\phi_w = \phi_u$,
 \item $x_w = x_u + \phi_v y_u =
0$, and
 \item $y_w = y_u$,
 \end{itemize}
 as desired.
 \end{proof}

\begin{prop} \label{HinN-linear}
 If there does not exist a continuous curve $h^t \to \infty$ in~$U$, such that
$\rho(h^t) \asymp h^t$, then
 \begin{enumerate}
 \item \label{HinN-linear-phi=0&eta2=}
 for every nonzero element~$z$ of~$\Dh$, we have $|\eta_z|^2 + \xx_z \yy_z
\neq 0$; and
 \item \label{HinN-linear-phi=0&=0}
 for every element~$u$ of~$\Lie U_{\phi=0}$, such that $\dimF (\F x_u + \F
y_u) = 1$, we have
 $$ \xx_u |y_u|^2 + \yy_u |x_u|^2 + 2 \Im ( x_u y_u^{\dagger}
\eta_u^{\dagger}) \neq 0 .$$
 \end{enumerate}
 \end{prop}

\begin{proof}[Proof of the contrapositive]
 \pref{HinN-linear-phi=0&eta2=}
 Suppose there is a nonzero element~$z$ of~$\Dh$ with $\Delta(z) = 0$. Let
$h^t = \exp(tz) = \Id + tz$ \see{exp(phi=0)}. We have
 $$ h_{j,k} = 
 \begin{cases}
 O(t) & \mbox{for all $j,k$}, \\
 O(1) & \mbox{if $(j,k) \notin \{1,2\} \times \{n+1,n+2\}$}. 
 \end{cases}
 $$
 Then, because $\Delta(h^t) = 0$, it is easy to see that $\rho(h) \asymp t$.
Also, we have
  $h^t = \Id + tz \asymp t$, so $\rho(h^t) \asymp t \asymp h^t$, as desired.

 \pref{HinN-linear-phi=0&=0}
  Suppose there is an element~$u$ of~$\Lie U_{\phi=0}$, such that $\dimF (\F
x_u + \F y_u ) = 1$, and
 \begin{equation} \label{xx+yy=0}
 \xx_u |y_u|^2 + \yy_u |x_u|^2 + 2 \Im ( x_u y_u^{\dagger}  \eta_u^{\dagger})
= 0 .
 \end{equation}
Let $h = h^t = \exp(tu)$. 

\setcounter{case}{0}

\begin{case} \label{HinNlinearpf-x=0}
 Assume $x_u = 0$.
 \end{case}
 Because $\dimF (\F x_u + \F y_u ) = 1$, we must have $y_u \neq 0$. Then,
from~\eqref{xx+yy=0}, we know that $\xx_u = 0$. So, from~\eqref{exp(phi=0)},
we see that
 \begin{itemize}
 \item $h_{2,n+1} \asymp |y_u|^2 t^2 \asymp t^2$, 
 \item $h_{j,k} = O (t)$ whenever $(j,k) \neq (2,n+1)$, and
 \item $h_{j,k} = O(1)$ whenever $j \neq 2$ and $k \neq n+1$. 
 \end{itemize}
 This implies that $\rho(h) \asymp t^2 \asymp h$.

\begin{case} \label{HinNlinearpf-y=0}
 Assume $y_u = 0$.
 \end{case}
 This is similar to Case~\ref{HinNlinearpf-x=0}. (In fact, this can be
obtained as a corollary of Case~\ref{HinNlinearpf-x=0} by replacing~$H$ with
its conjugate under the Weyl reflection corresponding to the root~$\alpha$.)

\begin{case}
 Assume $y_u \neq 0$.
 \end{case}
 Because $\dimF (\F x_u + \F y_u ) = 1$, Lemma~\ref{x=0} implies there is
some $g \in \rsg_\alpha$, such that, letting $w = g^{-1} u g$, we have
$\phi_w = \phi_u = 0$, $x_w = 0$, and $y_w = y_u \neq 0$. We show below that
\eqref{xx+yy=0} is satisfied with~$w$ in the place of~$u$, so, from
Case~\ref{HinNlinearpf-x=0}, we conclude that $\rho \bigl( \exp(tw) \bigr)
\asymp \exp(tw)$. Thus, the desired conclusion follows from
Lemma~\ref{conjcurve} (with $k = 1$).

To complete the proof, we now show that \eqref{xx+yy=0} is satisfied with~$w$
in the place of~$u$. (This can be verified by direct calculation, but we give
a more conceptual proof.) Because $g^{-1} \in \rsg_\alpha$, multiplication by
$g^{-1}$ on the left performs a row operation on the first two rows of~$h$;
likewise, multiplication by $g$ on the right performs a column operation on
the last two columns of~$h$. These operations do not change the determinant
$\Delta(h)$: thus 
 $$\Delta \bigl( \exp(tw) \bigr) = \Delta \bigl( \exp(tu) \bigr) .$$

From \eqref{exp(phi=0)} and the definition of~$\Delta$, we see that
 $$ \Delta \bigl( \exp(tu) \bigr)
 = - \frac{1}{4} \bigl( |x_u|^2 |y_u|^2 - |x_u y_u^\dagger|^2 \bigr) t^4
 + \bigl( \xx_u |y_u|^2 + \yy_u |x_u|^2 + 2 \Im ( x_u y_u^{\dagger} 
\eta_u^{\dagger}) \bigr) t^3
 + O(t^2) .$$
 Because $\dimF (\F x_u + \F y_u) = 1$, we have $|x_u|^2 |y_u|^2 - |x_u
y_u^\dagger|^2 = 0$, so this simplifies to
 $$ \Delta \bigl( \exp(tu) \bigr)
 = \bigl( \xx_u |y_u|^2 + \yy_u |x_u|^2 + 2 \Im ( x_u y_u^{\dagger} 
\eta_u^{\dagger}) \bigr) t^3
 + O(t^2) .$$
 Thus, \eqref{xx+yy=0} is equivalent to the condition that $\Delta \bigl(
\exp(tu) \bigr) = O(t^2)$. Then, since 
 $$\Delta \bigl( \exp(tw) \bigr) = \Delta \bigl( \exp(tu) \bigr) = O(t^2) ,$$
 we conclude that \eqref{xx+yy=0} is also valid for~$w$.
 \end{proof}

\begin{lem} \label{dim(U+Z)<3}
  If there does not exist a continuous curve $h^t \to \infty$ in~$U$, such
that $\rho(h^t) \asymp h^t$, then 
 $$\dim \Dh + \dim \Lie U/ \Lie U_{\phi=0} \le 2\df-1 .$$

Furthermore, if equality holds, and $\F = \complex$, then $\Lie U = \Lie
U_{\phi=0}$ and $\dim \Dh = 3$.
 \end{lem}

\begin{proof}

\setcounter{case}{0}

\begin{case}
 Assume $\F = \complex$.
 \end{case}
 Because $|\eta_z|^2 + \xx_z \yy_z$ is a quadratic form of signature $(1,3)$
on~$\Dn$, we know, from \fullref{HinN-linear}{phi=0&eta2=}, that $\dim \Dh
\le 3 = 2\df-1$.

Thus, we may assume $\Lie U / \Lie U_{\phi=0} \neq 0$, so there is some $u \in
\Lie U$, such that $\phi_u \neq 0$. 

\begin{subcase} \label{dim(U+Z)<3pf-yy=0}
 Assume there exists $z \in \Dh$, such that $\yy_z = 0$ and $\eta_z \notin
\real \phi_u$.
 \end{subcase}
 From~\pref{[u,v]}, we see that $[u,z] \in \Dh$, with
 $\xx_{[u,z]} = -\Im (\phi_u \cjg{\eta_z}) \neq 0$
 and
 $\yy_{[u,z]} = \eta_{[u,z]} = 0$. This contradicts
\fullref{HinN-linear}{phi=0&eta2=}.

\begin{subcase}
 Assume there exists $z \in \Dh$, such that $\yy_z \neq 0$.
 \end{subcase}
 From~\pref{[u,v]}, we see that $[u,z]$ is an element of~$\Dh$, such that
 $\yy_{[u,z]} = 0$, and
 $\eta_{[u,z]} = \phi_u \yy_z$ is a purely imaginary multiple of~$\phi_u$. So
Subcase~\ref{dim(U+Z)<3pf-yy=0} applies (with $[u,z]$ in the place of~$z$).

\begin{subcase}
 Assume $\yy_z = 0$ and $\eta_z \in \real \phi_u$, for all $z \in \Dh$.
 \end{subcase}
 From~\fullref{HinN-linear}{phi=0&eta2=}, we see that $\Dh \cap
\rsp_{2\alpha+2\beta} = \{0\}$, so the assumption of this subcase implies
$\dim \Dh \le 1$. Thus,
 $$ \dim \Dh + \dim \Lie U/ \Lie U_{\phi=0}  \le 1 + 2 = 3 = 2\df-1 ,$$
 so the desired inequality holds. 

If equality holds, then $\dim \Dh = 1$ and $\dim \Lie U/ \Lie U_{\phi=0} =
2$. Thus, we may choose $z \in \Dh$, such that $z \neq 0$, and $u' \in
\Lie U$, such that $\real \phi_u + \real \phi_{u'} = \complex$. From the
assumption of this subcase, we know that $\eta_z \in \real \phi_u$; thus,
$\eta_z \notin \real \phi_{u'}$. Therefore, Subcase~\ref{dim(U+Z)<3pf-yy=0}
applies, with $u'$ in the place of~$u$.

\begin{case}
 Assume $\F = \real$.
 \end{case}
 Because $\dim \Dh \le \dim \rsp_{\alpha+2\beta} = 1$ and $\dim \Lie U /
\Lie U_{\phi=0} \le \dim \rsp_\alpha = 1$, the desired inequality holds
unless $\Dh \neq 0$ and $\Lie U / \Lie U_{\phi=0} \neq 0$. Thus, we may
assume there is some nonzero $z \in \Dh$ and some $u \in \Lie U$, such
that $\phi_u \neq 0$.

\begin{subcase}
 Assume $y_u = 0$.
 \end{subcase}
 We may assume  $|x_u|^2 + \phi_u \eta_u \neq 0$ (by replacing~$u$ with $u +
z$, if necessary). Let $h^t = \exp(tu)$. From~\eqref{exp(y=0)}, we see that 
 $h^t_{1,n+2} \asymp t^{2}$, but
  $$ h^t_{j,k} =
 \begin{cases}
 O(t) & \hbox{if $(j,k) \neq (1,n+2)$,} \\
 O(1) & \hbox{if $j \neq 1$ and $k \neq n+2$} .
 \end{cases}
 $$
 Therefore $\rho(h^t) = O(t^2) = O(h^t)$, so Lemma~\fullref{Owalls}{linear}
implies that $\rho(h^t) \asymp h^t$. This is a contradiction.

\begin{subcase}
 Assume $y_u \neq 0$.
 \end{subcase}
 Let $v$ be the element of $\rsp_\beta$ with $y_v = -(1/\phi_u) x_u$, and let
$w = \exp(-v) u \exp(v)$. Then $x_w = 0$ \seeand{conjugation}{[u,v]}. Thus,
by replacing $H$ with the conjugate $\exp(-v) H \exp(v)$ \see{conjcurve}, we
may assume $x_u = 0$. For any large real number~$t$, let $h = h^t$ be the
element of $\exp(t u + \rsp_{\alpha+2\beta})$ that satisfies $\eta_h =
-\phi_h \lVert y_h\rVert^2/{12}$. Then, from~\eqref{exp}, we see that
 $$ h =
 \begin{pmatrix}
 1 & \phi_h & \frac{1}{2} \phi_h y_h & - \frac{1}{4} \phi_h |y_h|^2
   & \frac{1}{8} \phi_h^2 |y_h|^2 \\
   &  1   & y_h                  & - \frac{1}{2} |y_h|^2
   &  \frac{1}{4} \phi |y_h|^2 \\
   &      &  1                 & -y_h^\dagger
   &  \frac{1}{2} \phi_h y_h^\dagger   \\
   &      &                    &  1                   &   -\phi_h \\
   &      &                    &                      &  1 \\
 \end{pmatrix}
 .$$
 Clearly, we have 
 $h \asymp \phi_h^2 |y_h|^2$.

A calculation shows that $\Delta(h) = 0$, and certain other $2 \times 2$
minors also have cancellation. With this in mind, it is not difficult to
verify that $\rho(h) \asymp \phi_h^2 |y_h|^2 \asymp h$ (see \cite[Case~3 of
pf.\ of 5.12($3 \Rightarrow 2$)]{OhWitte-CDS} for details).
This is a contradiction.
 \end{proof}

\begin{prop} \label{dim(U)<2n}
 If there does not exist a continuous curve $h^t \to \infty$ in~$U_{\phi=0}$,
such that $\rho(h^t) \asymp h^t$, then
 $$ \dim \Lie U_{\phi=0}/ \Dh \le
 \begin{cases}
 \df(n-2) & \mbox{if $n$ is even} \\
 \df(n-3) & \mbox{if $n$ is odd and $n \ne 3$} \\
  \df-1   & \mbox{if $n = 3$} .
 \end{cases}
 $$
 Furthermore, 
 \begin{enumerate} \renewcommand{\theenumi}{\alph{enumi}}
 \item \label{dim(U)<2n-eqeven}
 if equality holds, and $n$ is even, then $\dimF (\F x_u + \F y_u) = 2$, for
every $u \in \Lie U_{\phi=0} \smallsetminus \Dh$; 
 \item \label{dim(U)<2n-eq=3}
 if equality holds, and $n = 3$, then  $\dim \Dh \le \df$.
 \end{enumerate}
 \end{prop}

\begin{proof}
 By passing to a subgroup, we may assume $\Lie U = \Lie U_{\phi=0}$.
  Let $V$ be the projection of $\Lie U$ to $\rsp_\beta + \rsp_{\alpha+\beta}$;
then $\dim V = \dim \Lie U/ \dim \Dh$.

\setcounter{case}{0}

\begin{case} \label{dim(U)<2n-<xy>pf}
 Assume $\dimF (\F x_u + \F y_u ) = 2$ for every $u \in \Lie U \smallsetminus \Dh$.
 \end{case}

\begin{subcase}
 Assume $n$ is even.
 \end{subcase}
 From Theorem~\fullref{HinN-linear}{phi=0&eta2=}, we know that $V$ does not
intersect $\rsp_\beta$ (or $\rsp_{\alpha+\beta}$, either, for that matter),
 so 
 $$\dim V + \dim \rsp_\beta
 \le \dim (\rsp_\beta + \rsp_{\alpha+\beta})
 = \dim \rsp_\beta + \dim \rsp_{\alpha+\beta} .$$
 Therefore
 $$ \dim \Lie U/ \dim \Dh = \dim V \le \dim \rsp_{\alpha+\beta} = \df(n-2)
,$$
 as desired.
 (If equality holds, then we have Conclusion~\pref{dim(U)<2n-eqeven}.)

\begin{subcase} \label{dim(U)<2npf-<xy>-n=2}
 Assume $n$ is odd.
 \end{subcase}

\begin{stepinsubcase} \label{dim(U)<2npf-U/Z}
 We have $\dim \Lie U / \Dh \le \df(n-3)$.
 \end{stepinsubcase}
 Suppose not: then
 $$\dim V \ge \df(n-3) + 1 .$$
 (This will lead to a contradiction.) Let $X = \{x_v \mid v \in V\}$, so $X$
is a $\real$-subspace of~$\F^{n-2}$. For each $x \in X$, there is some $v =
v(x) \in V$, such that $x_v = x$; define $f(x) = y_{v(x)}$. By the assumption
of this case, we know 
 $$V \cap \rsp_\beta = \{0\} ,$$
 so $v(x)$ is uniquely determined by~$x$; thus, $f \colon X \to \F^{n-2}$ is
a well-defined $\real$-linear map. Also, again from the assumption of this
case, we know that 
 \begin{equation} \label{TxnotinFx}
 \mbox{$f(x) \notin \F x$ for every nonzero $x \in X$.} 
 \end{equation}
 Because $V \cap \rsp_\beta = 0$, we have 
 $$\dim X = \dim V \ge \df(n-3) + 1= \dim \F^{n-2} - (\df-1) .$$

 If $\F = \real$ (that is, if $\df = 1$), this implies $X = \real^{n-2}$, so
$f$ is defined on all of $\real^{n-2}$. Because $n$~is odd, this implies that
$f$ has a real eigenvalue, which contradicts \eqref{TxnotinFx}.

We may now assume $\F = \complex$. 
 Let 
 \begin{itemize}
 \item $E = (X \times \complex^{n-2})/ {\equiv}$, where $(x,v) \equiv
(-x,-v)$,
 \item $\projX$ be the projective space of the real vector space~$X$, and
 \item  $\zeta(x,v) = [x] \in \projX$, for $(x,v) \in E$,
 \end{itemize}
 so $(E, \zeta)$ is a vector bundle over $\projX$. 

Define
 $g \colon X \to \complex^{n-2}$ by $g(x) = i x$.
 Any $\real$-linear transformation $Q \colon X \to \complex^{n-2}$ is a
continuous function, such that $Q(-x) = - Q(x)$ for all $x \in X$; that is, a
section of $(E, \zeta)$. Thus, $\Id$, $f$, and~$g$ each define a section of
$(E, \zeta)$. Furthermore, these three sections are pointwise linearly
independent over~$\real$, because \eqref{TxnotinFx} implies that $x$, $f(x)$,
and $ix$ are linearly independent over~$\real$, for every nonzero $x \in X$.
On the other hand, the theory of characteristic classes \cite[Prop.~4,
p.~39]{MilnorStasheff} implies that $(E, \zeta)$ does not have three
pointwise $\real$-linearly independent sections (see
\cite[Lem.~8.2]{IozziWitte-CDS} for details). This is a contradiction.

\begin{stepinsubcase}
 Completion of the proof of Subcase~\ref{dim(U)<2npf-<xy>-n=2}.
 \end{stepinsubcase}
 From Step~\ref{dim(U)<2npf-U/Z}, we see that the desired inequality holds.

We may now assume $n = 3$ and $\dim \Lie U/\Dh =
\df - 1$. Since $\dim \Lie U/\Dh \le \df(n-3) = 0$, we must have
$\df = 1$, so $\F = \real$. Therefore $\rsp_{2\alpha} = \rsp_{2\beta} = 0$,
so 
 $$ \dim \Dh \le \dim \rsp_{\alpha+2\beta} = \df ,$$
 as desired.

\begin{case}
 Assume there is some $v \in \Lie U \smallsetminus \Dh$, such that $\dimF (\F
x_v + \F y_v) = 1$.
 \end{case}

\begin{subcase} \label{dimUpf-x=0}
 Assume $x_v = 0$. 
 \end{subcase}
 Since $v \notin \Dh$, we must have $y_v \neq 0$. Then 
 $\xx_{v+z} \neq 0$ for every $z \in \Dh$ (otherwise
\fullref{HinN-linear}{phi=0&=0} yields a contradiction); this implies 
 $$\xx_v \neq 0 ,$$
 and 
 $$ \mbox{$\xx_z = 0$ for every $z \in \Dh$.} $$
 Because $\xx_v \neq 0$, we know that $\F \neq \real$; so
 $$\F = \complex .$$
 Since $\xx_z = 0$ for every $z \in \Dh$, but $\Dh \cap \rsp_{2\beta} =
0$ \fullsee{HinN-linear}{phi=0&eta2=}, we must have $\eta_z \neq 0$
for every $z \in \Dh$. Therefore
 $$ \dim \Dh \le \dim \rsp_{\alpha+2\beta} = \df = 2 .$$

Let $p \colon V \to \rsp_{\alpha+\beta}$ be the natural projection. Note that 
 $$\dim \ker p = 1 .$$
 (If $v' \in \Lie U$, with $x_{v'} = 0$, then there is some $t \in \real$,
such that $\xx_{v' + tv} = \xx_{v'} + t \xx_v = 0$. We also have $x_{v' + tv}
= 0$, so, from~\fullref{HinN-linear}{phi=0&=0}, we see that $v' + tv \in \Dh$. Thus $v' \in \real v + \Dh$. So $\ker p = (\real v + \Dh)/\Dh$
is 1-dimensional.)

Because $\xx_z = 0$ for every $z \in \Dh$, and $\Lie U$ is a Lie algebra,
we see, from~\pref{[u,v]}, that $p(V)$ must be a totally isotropic subspace
for the symplectic form $i \Im(x \tilde x^{\dagger})$, so 
 $$\dim p(V) \le \frac{1}{2} \dim \rsp_{\alpha+\beta} = n-2 .$$
Therefore 
 $$ \dim \Lie U/ \Dh
 = \dim V
 = \dim p(V) + \dim \ker p
 \le (n-2) + 1
 = n-1 .$$
 This completes
the proof if $n \neq 3$:
 \begin{itemize}
 \item If $n$ is even, then, because $n \ge 4$, we have $n - 1 <
2(n-2) = \df(n-2)$.
 \item If $n > 3$ is odd, then $n \ge 5$, so $n-1 \le 2(n-3) = \df(n-3)$.
 \end{itemize}

Now let $n = 3$, and suppose $\dim V = 2$. (This will lead to a
contradiction.) Because equality is attained in the proof above, we must have
$\dim p(V) = n-2 = 1$, so there exists
$w \in \Lie U$ with $x_w \neq 0$. For $t \in \real$, let $w_t = w + tv$. Then
 \begin{align*}
 \xx_{w_t} |y_{w_t}|^2 + \yy_{w_t}|x_{w_t}|^2 + 2 \Im(
x_{w_t}y_{w_t}^{\dagger} \eta_{w_t})
 &= t^3 \xx_v |y_v|^2 + O(t^2) \\
 &\to
 \begin{cases}
 +\xx_v \infty & \text{as $t \to \infty$} \\
 -\xx_v \infty & \text{as $t \to -\infty$}
 .
 \end{cases}
 \end{align*}
 Thus, this expression changes sign, so it must vanish for some~$t$. On the
other hand, since $n = 3$, we have $\dimC (\complex x + \complex y) \le 1$
for every $x,y \in \complex^{n-2} = \complex$. Thus
\fullref{HinN-linear}{phi=0&=0} yields a contradiction.

\begin{subcase} \label{dimUpf-y=0}
 Assume $y_v = 0$.
 \end{subcase}
 This is similar to Subcase~\ref{dimUpf-x=0}. (In fact, this can be
obtained as a corollary of Subcase~\ref{dimUpf-x=0} by replacing~$H$ with
its conjugate under the Weyl reflection corresponding to the root~$\alpha$.)

\begin{subcase}
 Assume $y_v \neq 0$.
 \end{subcase}
 Because $\dimF (\F x_v + \F y_v ) = 1$, Lemma~\ref{x=0} implies there is
some $g \in \rsg_\alpha$, such that, letting $w = g^{-1} v g$, we have
$\phi_w = \phi_u = 0$, $x_w = 0$, and $y_w = y_v \neq 0$. There is no harm in
replacing $H$ with $g^{-1} H g$ \see{conjcurve}. Then
Subcase~\ref{dimUpf-x=0} applies (with $w$ in the place of~$v$).
 \end{proof}

\begin{thm} \label{maxnolinear}
 Recall that Assumptions~\ref{StandingSU2F} are in effect.

 If there does not exist a continuous curve $h^t \to \infty$ in~$H$, such that
$\rho(h^t) \asymp h^t$, then
 \begin{equation} \label{maxnolinear-<} 
 \dim H \le
 \begin{cases}
 \df n & \hbox{if $n$ is even} \\
 \df (n-1) & \hbox{if $n$ is odd.}
 \end{cases}
 \end{equation}
 Furthermore, if equality holds, and $n$~is even, then
 \begin{enumerate}
 \item \label{maxnolinear-TU}
 $\Lie H = (\ker \alpha) \ltimes \Lie U$;
 \item \label{maxnolinear-phi0}
 $\phi_u = 0$ for every $u \in \Lie U$;
 \item \label{maxnolinear-<xy>}
 $\dimF ( \F x_u + \F y_u ) = 2$, for every $u \in \Lie U \smallsetminus \Dh$;
 \item \label{maxnolinear-Z}
 $|\eta_z|^2 + \xx_z \yy_z \neq 0$ for every
nonzero $z \in \Dh$;
 \item \label{maxnolinear-dimU/Z}
 $\dim \Lie U / \Dh = \df(n-2)$; and
 \item \label{maxnolinear-dimZ}
 $\dim \Dh = 2\df - 1$.
 \end{enumerate}
 \end{thm}

\begin{proof}
  Let
 $$m =
 \begin{cases}
 \df(n-2) & \hbox{if $n$ is even} \\
 \df(n-3) & \hbox{if $n \ge 5$ is odd} \\
 \df-1 & \hbox{if $n = 3$}
 .
 \end{cases}
 $$ 
 From Lemmas~\fullref{dimT}{A} and~\ref{dim(U+Z)<3}, and
Proposition~\ref{dim(U)<2n}, we have
 \begin{equation} \label{<2d+m}
 \begin{split}
 \dim H
 &\le \dim \Lie H/\Lie U + (\dim \Lie U/\Lie U_{\phi=0} + \dim \Dh)
 + \dim \Lie U_{\phi=0}/ \Dh \\
 &\le 1 + (2\df-1) + m \\
 &= 2\df+m .
 \end{split}
 \end{equation}
 This implies the desired inequality, unless $n = 3$, $\F = \complex$, and we
have equality in both Lemma~\ref{dim(U+Z)<3} and Proposition~\ref{dim(U)<2n}.
This is impossible, because equality in Lemma~\ref{dim(U+Z)<3} requires $\dim
\Dh = 3$, but Proposition~\fullref{dim(U)<2n}{eq=3} implies $\dim \Dh
\le 2$.

\begin{pfassump}
 In the remainder of the proof, we assume that equality holds in
\eqref{maxnolinear-<}, and that $n$ is even.
 \end{pfassump}
 Proposition~\fullref{HinN-linear}{phi=0&eta2=} implies~\pref{maxnolinear-Z}.

\setcounter{case}{0}

\begin{case} \label{maxnolinearpf-F=C}
 Assume $\F = \complex$.
 \end{case}
 Because equality holds, Lemma~\ref{dim(U+Z)<3}
implies~\pref{maxnolinear-phi0} and~\pref{maxnolinear-dimZ}. Then
Proposition~\fullref{dim(U)<2n}{eqeven} implies~\pref{maxnolinear-<xy>}
(because $\Lie U = \Lie U_{\phi=0}$). Since $\Lie U = \Lie U_{\phi=0}$
\see{maxnolinear-phi0} and equality holds in \eqref{maxnolinear-<}, we have
 \begin{equation} \label{maxnolinearpf-U/Z}
 \dim \Lie U/ \Dh
 = \dim \Lie U_{\phi=0}/ \Dh
 = m
 = \df(n-2)
 \end{equation}
 and
  \begin{equation} \label{maxnolinearpf-dimZ}
 \dim \Dh
 = \dim \Lie U/ \Lie U_{\phi=0} + \dim \Dh
 = 2\df - 1 ,
 \end{equation}
 so \pref{maxnolinear-dimU/Z} and~\pref{maxnolinear-dimZ} hold.

Let $T = A \cap (HN)$. Corollary~\ref{TnormsU} implies that $T$
normalizes~$\Lie U$, so, from~\pref{maxnolinear-<xy>} and
Lemma~\ref{rootdecomp}, we see that $T \subset \ker \alpha$. On the other
hand, $\dim T = \dim \Lie H / \Lie U$, so, from equality in~\eqref{<2d+m}, we
conclude that $\dim T = 1$. Therefore $T = \ker \alpha$.

Suppose $\psi \colon \ker \alpha \to \rsg_\alpha$ is any continuous group
homomorphism, such that $\psi( \ker \alpha)$ normalizes~$U$.
From~\pref{maxnolinear-<xy>} and~\pref{[u,v]}, we see that
$N_{\rsg_\alpha}(U) = e$, so $\psi$ must be trivial. This implies that
\fullref{not-semi}{not} cannot apply here, so \fullref{not-semi}{TU}
yields~\pref{maxnolinear-TU}.

\begin{case}
 Assume $\F = \real$.
 \end{case}
 Proposition~\fullref{dim(U)<2n}{eqeven} implies that $\dimF (\F x_u + \F y_u)
= 2$ for every $u \in \Lie U_{\phi=0} \smallsetminus \Dh$.

Suppose \pref{maxnolinear-phi0} is false. Then there is some $u \in \Lie U$,
such that $\phi_u \neq 0$. Also, because $\dim \Lie U_{\phi=0}/\Dh = m >
0$, we may fix some $v \in \Lie U_{\phi=0} \smallsetminus \Dh$. Then, letting
$w = [u,v]$, we see, from~\pref{[u,v]}, that $y_w = 0$ and $x_w \neq 0$, so
$\dimF (\F x_w + \F y_w) = 1$. This contradicts the conclusion of the
preceding paragraph.

Conclusion \pref{maxnolinear-<xy>} follows from~\pref{maxnolinear-phi0}
and~\fullref{dim(U)<2n}{eqeven}.

Conclusion~\pref{maxnolinear-TU} can be established by arguing as in the last
two paragraphs of Case~\ref{maxnolinearpf-F=C}.

Equations~\pref{maxnolinearpf-U/Z} and~\pref{maxnolinearpf-dimZ} establish
\pref{maxnolinear-dimU/Z} and~\pref{maxnolinear-dimZ}.
 \end{proof}

\subsection{Subgroups with no nearly quadratic curve}
 Our goal is to prove Theorem~\ref{bestnosquare}; we start with two
preliminary results.

\begin{lem} \label{HinN-square}
 If there does not exist a continuous curve $h^t \to \infty$ in~$U$, such that
$\rho(h^t) \asymp \lVert h^t\rVert^2$, then
 \begin{enumerate}
 \item \label{HinN-square-indep}
 for every element~$u$ of~$\Lie U_{\phi=0}$, we have $\dimF( \F
x_u + \F y_u) \le 1$;
 \item \label{HinN-square-eta2neq}
 for every element~$z$ of~$\Dh$, we have $|\eta_z|^2 + \xx_z \yy_z = 0$;
and
 \item \label{HinN-square-y=0+yy=0}
 for every  element~$u$ of~$\Lie U$, such that $\phi_u \neq 0$, $y_u = 0$,
and $\yy_u = 0$, we have $|x_u|^2 + 2 \Re(\phi_u \cjg{\eta_u}) \neq 0$.
 \end{enumerate}
 \end{lem}

\begin{proof}[Proof of the contrapositive.]
 \pref{HinN-square-indep}
 Suppose there is an element~$u$ of~$\Lie U_{\phi=0}$, such
that $\dimF( \F x_u + \F y_u) = 2$.
 Let $h^t = \exp(t u)$. Then,
from \eqref{exp(phi=0)}, we see that 
 $h^t = O(t^2)$. Furthermore,
 $$ \Delta(h^t)
 = \det
 \begin{pmatrix}
 \eta_u t - \frac{1}{2} x_u y_u^\dagger t^2
 & \xx_u t - \frac{1}{2} |x_u|^2
t^2 \\
 \yy_u t - \frac{1}{2} |y_u|^2 t^2 &
 -\cjg{\eta_u} t - \frac{1}{2} y_u x_u^\dagger t^2
 \end{pmatrix}
 = \frac{1}{4} \bigl| |x_u| |y_u|^2 - |x_u y_u^\dagger|^2 \bigr| t^4
 + O(t^3)
 . $$
 Because $\dimF( \F x_u + \F y_u) = 2$,
 we have $|x_u| |y_u| > |x_u y_u^\dagger|$, 
 so $|x_u|^2 |y_u|^2 - |x_u y_u^\dagger|^2 \neq 0$;
 therefore $\Delta(h^t) \asymp t^4$, so
 $$ \lVert h^t\rVert^2 = O(t^4) = O \bigl( \Delta(h^t) \bigr) = O\bigl(
\rho(h^t) \bigr), $$
 so Lemma~\fullref{Owalls}{square} implies that
 $\rho(h^t) \asymp \lVert h^t\rVert^2$,
 as desired.

\pref{HinN-square-eta2neq}
 Suppose there is an element~$z$ of~$\Dh$, such that $|\eta_z|^2 + \xx_z
\yy_z \neq 0$; in other words, we have $\Delta(z) \neq 0$.
 Let $h^t = \exp(t z) = \Id + tz$ \see{exp(phi=0)}. Then $h^t = O(t)$ and
 $$t^2 \asymp \Delta(z) t^2 = \Delta(h^t) = O \bigl( \rho(h^t) \bigr) ,$$
 so 
 $$ \lVert h^t\rVert^2 = O(t^2) = O \bigl( \rho(h^t) \bigr) ,$$
 so Lemma~\fullref{Owalls}{square} implies that
 $\rho(h^t) \asymp \lVert h^t\rVert^2$,
 as desired.

\pref{HinN-square-y=0+yy=0}
 Suppose there is an element~$u$ of~$\Lie U$, such that
 $\phi_u \neq  0$,
 $y_u = 0$,
 $\yy_u = 0$, and
 $|x_u|^2 + 2 \Re(\phi_u \cjg{\eta_u}) = 0$.
 Let $h^t = \exp(tu)$. From~\pref{exp(y=0)}, we see that $h^t = \Id + tu$
(note that, because $|x_u|^2 + 2 \Re(\phi_u \cjg{\eta_u}) = 0$, we have $\Re
h^2_{1,n+2} = 0$). Then
 $ h^t = O(t) $
 and
 $$ \lVert \rho(h^t) \rVert
 \ge \left| \det
 \begin{pmatrix}
 h^t_{1,2} & h^t_{1,n+2} \\
 h^t_{n+1,2} & h^t_{n+1,n+2}
 \end{pmatrix}
 \right|
 = \left| \det
 \begin{pmatrix}
 t\phi_u & * \\
 0 & - t\phi_u^\dagger
 \end{pmatrix}
 \right|
 \asymp t^2 .$$
 So $\lVert h^t\rVert^2 = O(t^2) = O \bigl( \rho(h^t) \bigr)$. Thus,
Lemma~\fullref{Owalls}{square} implies that
 $\rho(h^t) \asymp \lVert h^t\rVert^2$,
 as desired.
 \end{proof}

The following lemma obtains a dimension bound from
Condition~\fullref{HinN-square}{indep}.

\begin{lem} \label{dimV<d(n-2)}
 If $V$ is a $\real$-subspace of $\F^{n-2} \oplus \F^{n-2}$, such that
 $\dimF( \F x + \F y ) \le 1$ for every $(x,y) \in V$,
 then either
 \begin{enumerate}
 \item \label{dimV<d(n-2)-n>3}
 $\dim V \le \df(n-2)$; or
 \item \label{dimV<d(n-2)-n=3}
 $n = 3$ and $\dim V \le 2\df$.
 \end{enumerate}
 \end{lem}

\begin{proof}
 Because $\dimR \F^{n-2} = \df(n-2)$, we may assume that there exist nonzero
$x_0,y_0 \in \F^{n-2}$, such that $(x_0,0) \in V$ and $(0,y_0) \in V$
(otherwise, the projection to one of the factors of $\F^{n-2} \oplus \F^{n-2}$
is injective when restricted to~$V$, so \pref{dimV<d(n-2)-n>3} holds). Then
$(x_0,y_0) \in V$, so, by assumption, we have $\dimF(\F x_0 + \F y_0) \le
1$. Because $x_0$ and~$y_0$ are nonzero, this implies $\F x_0 = \F y_0$.

\setcounter{step}{0}

\begin{step} \label{dimV<d(n-2)Pf-yinFx}
 For all $(x,y) \in V$, we have $y \in \F x_0$.
 \end{step}
 We may assume $y \neq 0$ (otherwise the desired conclusion is obvious).
Then, since $\dimF(\F x + \F y ) \le 1$, we conclude that $x \in \F y$.
Similarly, because 
 $$(x + x_0, y) = (x,y) + (x_0,0) \in V + V = V ,$$
 we must have $x+x_0 \in \F y$. Therefore
 $$ x_0 = (x+x_0) - x \in \F y - \F y = \F y .$$
 Since $x_0 \neq 0$, this implies $\F x_0 = \F y$, so $y \in \F x_0$, as
desired.

\begin{step} \label{dimV<d(n-2)Pf-V=FxF}
 We have $V \subset \F y_0 \oplus \F x_0$.
 \end{step}
 Given $(x,y) \in V$, Step~\ref{dimV<d(n-2)Pf-yinFx} asserts that $y \in \F
x_0$. By symmetry (interchanging the two factors of $\F^{n-2} \oplus
\F^{n-2}$), we must also have $x \in \F y_0$. So $(x,y) \in \F y_0 \oplus \F
x_0$, as desired.

\begin{step}
 Completion of the proof.
 \end{step}
 From Step~\ref{dimV<d(n-2)Pf-V=FxF}, we have
 $$ \dim V \le \dimR(\F y_0 \oplus \F x_0) = 2 \df .$$
 If $n \ge 4$, then \pref{dimV<d(n-2)-n>3} holds; otherwise,
\pref{dimV<d(n-2)-n=3} holds.
 \end{proof}

\begin{thm} \label{bestnosquare}
 Recall that Assumptions~\ref{StandingSU2F} are in effect.

 If there does not exist a continuous curve $h^t \to \infty$ in~$H$, such that
$\rho(h^t) \asymp \lVert h^t\rVert^2$, then $\dim H \le \df n$.

Furthermore, if equality holds, then $H$ is of the form
 $H = T \ltimes U$, where
 \begin{enumerate}
 \item \label{bestnosquare-T}
 $T = \ker \beta$,
  \item \label{bestnosquare-U}
 $\Lie U = \bigl( (\rsp_\alpha + \rsp_{\alpha+\beta} + \rsp_{\alpha+2\beta} )
\cap \Lie U \bigr)
 + \rsp_{2\alpha+2\beta}$,
 and
 \item \label{bestnosquare-x}
 $|x_u|^2 + 2\Re(\phi_u \cjg{\eta_u}) \neq 0$ for every $u \in \Lie U
\smallsetminus \rsp_{2\alpha+2\beta}$.
 \end{enumerate}
 \end{thm}

\begin{proof}
 Note that
 $$\dim \Lie H / \Lie U \le 1$$
 \fullsee{dimT}{A} and
 $$\dim \Lie U/\Lie U_{\phi=0} \le \dim \rsp_\alpha = \df .$$

\setcounter{step}{0}

\begin{step}
 We have $\dim \Lie U_{\phi=0} / \Dh \le \df(n-2)$. 
 \end{step}
 Suppose not. Let $V$ be the projection of~$\Lie U_{\phi=0}$ to $\rsp_{\beta}
+ \rsp_{\alpha+\beta}$. We have 
 $$\dim V = \dim \Lie U_{\phi=0} / \Dh > \df(n-2) ,$$
 and, for every $u \in \Lie U_{\phi=0}$ with $x_u \neq 0$, we have
 $\dimF ( \F x_u + \F y_u ) \le 1$
 (see~\fullref{HinN-square}{indep}), so Lemma~\ref{dimV<d(n-2)} implies that
$n = 3$. Therefore $\dim \rsp_\beta = \dim \rsp_{\alpha+\beta} = \df$. Then,
because $\dim V > \df(n-2) = \df$, we know that $V \cap \rsp_{\beta} \neq 0$
and $V \cap \rsp_{\alpha+\beta} \neq 0$; thus, there exist $u,v \in \Lie
U_{\phi=0}$, such that
 \begin{itemize} 
 \item $x_u = 0$, $y_u \neq 0$; and
 \item $x_v \neq 0$, $y_v = 0$.
 \end{itemize}
 Therefore $[u,v]$ is a nonzero element of~$\rsp_{\alpha+2\beta}$
\see{[u,v]}, so $\Delta \bigl( [u,v] \bigr) \neq 0$. This contradicts
Lemma~\fullref{HinN-square}{eta2neq}.

\begin{step}
 We have $\dim \Dh \le \df-1$. 
 \end{step}
 Suppose not: then, because $\dim \rsp_{2\alpha+2\beta} = \df-1$, there is
some $u \in \Dh \smallsetminus \rsp_{2\alpha+2\beta}$, and, because $\dim
\rsp_{2\beta} = \df-1$, there is some nonzero $v \in \Dh$, such that
$\yy_v = 0$. We must have $\eta_v = 0$ (otherwise
\fullref{HinN-square}{eta2neq} yields a contradiction); thus $v \in
\rsp_{2\alpha+2\beta}$. We must have $\yy_u \neq 0$ (otherwise
\fullref{HinN-square}{eta2neq} yields a contradiction). Thus, we see that
 $$| \eta |^2 + \xx_{u+tv} \yy_{u+tv}
 = |\eta_u|^2 + (\xx_u + t \xx_v) (t \yy_u) $$
 is nonconstant as a function of $t \in \real$, so
\fullref{HinN-square}{eta2neq} yields a contradiction.

\begin{step}
 The desired inequality. 
 \end{step}
We have
 \begin{align*}
 \dim \Lie H
 &\le \dim \Lie H/\Lie U + \dim \Lie U/\Lie U_{\phi=0} + \dim \Lie
U_{\phi=0}/ \Dh + \dim \Dh \\
 &\le 1 + \df + \df(n-2) + (\df-1) \\
 &= \df n ,
 \end{align*}
 as desired.

\begin{pfassump}
 In the remainder of the proof, we assume that $\dim H = \df n$. 
 \end{pfassump}
 We must have equality throughout the preceding paragraphs.

\begin{step} \label{Vinalpha+beta}
 We have $V \subset \rsp_{\alpha+\beta}$.
 \end{step}
 Suppose not: then there is some $v \in \Lie U_{\phi=0}$, such that $y_v \neq
0$. Let $u \in \Lie U \smallsetminus \Lie U_{\phi=0}$ and $w = [u,v]$. Then,
from~\eqref{[u,v]}, we see that $y_w = 0$ and $x_w \neq 0$, and that
 $[v,w] \in \rsp_{\alpha+2\beta} + \rsp_{2\alpha+2\beta}$.
 From~\fullref{HinN-square}{indep}, we have 
 $x_v \in \F y_v$ and $x_{v+2} \in \F y_{v+w} = \F y_v$,
 so 
 $$ x_w = x_{v+w} - x_v \in \F y_v - \F y_v = \F y_v .$$
 Therefore $x_w y_w ^\dagger \neq 0$, so
 $\eta_{[v,w]} \neq 0$ \see{[u,v]}, so \fullref{HinN-square}{eta2neq} yields
a contradiction.

\begin{step}
 We have $\Dh = \rsp_{2\alpha+2\beta}$. 
 \end{step}
  From Step~\ref{Vinalpha+beta}, together with the fact that
 $$\dim V = \dim \Lie U_{\phi=0} / \Dh = \df(n-2) = \dim
\rsp_{\alpha+\beta} ,$$
 we conclude that $V = \rsp_{\alpha+\beta}$. Therefore,
 $$ \Lie U_{\phi=0} + \Dn 
 = V + \Dn
 = \rsp_{\alpha+\beta} + \Dn ,$$
 so
 \begin{align*}
 \Lie U
 &\supset [\Lie U_{\phi=0}, \Lie U_{\phi=0}] \\
 &= [\Lie U_{\phi=0} + \Dn, \Lie U_{\phi=0} + \Dn] \\
 &= [\rsp_{\alpha+\beta} + \Dn, \rsp_{\alpha+\beta} + \Dn] \\
 &= [\rsp_{\alpha+\beta},\rsp_{\alpha+\beta}] \\
 &= \rsp_{2\alpha+2\beta} .
 \end{align*}
 Because $\dim \Dh = \df-1 = \dim \rsp_{2\alpha+2\beta}$, we must have
$\Dh = \rsp_{2\alpha+2\beta}$.

\begin{step}
 We have $\Dh = \rsp_{2\alpha+2\beta}$. 
 Let $T = (HN) \cap A$ be the projection of~$H$ to~$A$. 
 Then there exists $\sigma \in \{\beta, \alpha+\beta, \alpha+2\beta\}$, such
that $T = \ker(\alpha- \sigma)$, and, in the notation of
Lemma~\ref{rootdecomp}, we have
 $$\Lie U = (\Lie U \cap \Lie N^{=\alpha}) + (\Lie U \cap \Lie N^{\neq
\alpha}) .$$
 \end{step}%
 Because $T$ normalizes~$\Lie U$ \see{TnormsU}, we know, from
Lemma~\ref{rootdecomp}, that
 $\Lie U = (\Lie U \cap \Lie N^{=\alpha}) + (\Lie U \cap \Lie N^{\neq
\alpha})$. 
 Since $\dim \Lie U / \Lie U_{\phi=0} = \df$, we know that $\Lie U \cap \Lie
N^{=\alpha}$ projects nontrivially (in fact, surjectively) to~$\rsp_\alpha$.
On the other hand, we know that $\Lie U \cap \rsp_{\alpha} = 0$ (otherwise
\fullref{HinN-square}{y=0+yy=0} yields a contradiction). Therefore $\Lie
N^{=\alpha} \neq \rsp_\alpha$, so there must be a positive root~$\sigma \neq
\alpha$, such that $\sigma|_T = \alpha|_T$. Then $T \subset
\ker(\alpha-\sigma)$; since $\dim T = \dim H/U = 1$, we must have $T =
\ker(\alpha-\sigma)$.

 Because 
 $\Lie U \cap \Lie N^{\neq\alpha} \subset \Lie U_{\phi=0}$,
 we have
 $$\dim (\Lie U \cap \Lie N^{=\alpha})
 \ge \dim \frac{\Lie U \cap \Lie N^{=\alpha}}{\Lie U_{\phi=0} \cap \Lie
N^{=\alpha}}
 = \dim \frac
 {\Lie U / (\Lie U \cap \Lie N^{\neq\alpha})}
 {\Lie U_{\phi=0} / (\Lie U \cap \Lie N^{\neq\alpha})}
 = \dim \frac{\Lie U}{\Lie U_{\phi=0}}
 = \df .$$
 Then, since $\Lie U \cap \rsp_\alpha = 0$, we must have
 $$\dim \Lie N^{=\alpha}
 \ge \dim (\Lie U \cap \Lie N^{=\alpha}) + \dim \rsp_\alpha
 \ge \df + \df
 = 2\df .$$
 By inspection, we see that this implies $\sigma \notin \{2\alpha, 2\beta\}$,
so we conclude that $\sigma \in \{\beta, \alpha+\beta, \alpha+2\beta\}$, as
desired.

\begin{step} \label{bestnosquarepf-sigma}
 We have $\sigma \in \{ \alpha+\beta, \alpha+2\beta\}$, $T = \ker \beta$, and
 $\Lie N^{=\alpha} = \rsp_\alpha + \rsp_{\alpha+\beta} +
\rsp_{\alpha+2\beta}$.
 \end{step}
 Since $\ker \beta = \ker 2\beta$, it suffices to show $\sigma \neq \beta$.
Thus, let us suppose $\sigma = \beta$. (This will lead to a contradiction.)
 We have $\Lie N^{=\alpha} = \rsp_\alpha + \rsp_\beta$ (and recall that
$\Lie U \cap \rsp_\alpha = \{0\}$),
 so there is some $u \in \Lie U$, such that $\phi_u \neq 0$ and $y_u \neq 0$.
Because $V = \rsp_{\alpha+\beta}$, we have
 $$ \dim \{\, v \in V \mid x_v \in \F y_u \,\} = \df > \dim \Fim ,$$
 so there is some $v \in \Lie U_{\phi=0}$, such that $0 \neq x_v \in \F y_u$
and $\yy_v = 0$. Then $[u,v] \in \rsp_{\alpha+2\beta} + \rsp_{2\alpha +
2\beta}$, with $\eta_{[u,v]} \neq 0$ \see{[u,v]}, so
\fullref{HinN-square}{eta2neq} yields a contradiction.

\begin{step} \label{bestnosquarepf-TU}
 We have $H = (H \cap A) \ltimes ( H \cap N)$.
 \end{step}
 Suppose not: because $T = \ker \beta$, we conclude that there is some nonzero
$w \in \rsp_\beta + \rsp_{2\beta}$, such that $w$ normalizes~$\Lie U$
\see{not-semi}.

 If $y_w \neq 0$, then, because $V = \rsp_{\alpha+\beta}$, there is some $v
\in \Lie U_{\phi=0}$, such that $y_w \in \F x_v$ and $\yy_v = 0$. Then $[w,v]
\in \rsp_{\alpha+2\beta} + \rsp_{2\alpha + 2\beta}$, with $\eta_{[w,v]} \neq
0$ \see{[u,v]}, so \fullref{HinN-square}{eta2neq} yields a contradiction.

 If $y_w = 0$, then, since $w \neq 0$, we must have $\yy_w \neq 0$. There is
some $v \in \Lie U$ with $\phi_v \neq 0$. Then $[w,v] \in
\rsp_{\alpha+2\beta} + \rsp_{2\alpha + 2\beta}$, with $\eta_{[w,v]} \neq 0$
\see{[u,v]}, so \fullref{HinN-square}{eta2neq} yields a contradiction.

\begin{step}
 Completion of the proof.
 \end{step}
 \pref{bestnosquare-T} From Step~\ref{bestnosquarepf-TU}, We know that $H = T
\ltimes U$, and, from Step~\ref{bestnosquarepf-sigma}, that $T = \ker \beta$.

 \pref{bestnosquare-U} Since $\Lie N^{=\alpha} = \rsp_\alpha +
\rsp_{\alpha+\beta} + \rsp_{\alpha+2\beta}$, it suffices to show 
 $\Lie U \cap \Lie N^{\neq \alpha} = \rsp_{2\alpha+2\beta}$: given $v \in
\Lie U \cap \Lie N^{\neq \alpha}$, we wish to show $v \in
\rsp_{2\alpha+2\beta}$. Because $V = \rsp_{\alpha+\beta}$, we know that $y_v
= 0$. Thus, all that remains is to show that $\yy_v = 0$. If not, then
choosing $u \in \Lie U$ with $\phi_u \neq 0$, we see that $\eta_{[u,v]} \neq
0$ \see{[u,v]}. So \fullref{HinN-square}{eta2neq} yields a contradiction.

\pref{bestnosquare-x} From Lemma~\fullref{HinN-square}{y=0+yy=0}, we know that
 $|x_u|^2 + 2 \Re(\phi_u \cjg{\eta_u}) \neq 0$ for every $u \in \Lie U
\smallsetminus \rsp_{2\alpha+2\beta}$.
 \end{proof}

\section{Homogeneous spaces of $\SO(2,n)$ and $\SU(2,n)$}
\label{ProofSect}

This section proves two main results. Both assume that $G$ is
either $\SO(2,n)$ or $\SU(2,n)$.
 \begin{enumerate}
 \item Theorem~\ref{SUF->complete} shows that if $n$~is odd, and one or two
specific homogeneous spaces of~$G$ do not have tessellations, then no
interesting homogeneous space of~$G$ has a tessellation.
 \item Theorem~\ref{SUFevenTess} shows that if $n$~is even, then certain
deformations of the examples found by R.~Kulkarni and T.~Kobayashi
\see{KulkarniEg} are essentially the only interesting homogeneous spaces
of~$G$ that have tessellations.
 \end{enumerate}
 The classification results of~\S\ref{SUFlargeSect} (specifically,
Theorems~\ref{bestnosquare} and~\ref{maxnolinear}) play a crucial
role in the proofs.

 We use the notation $\SU(2,n;\F)$ of Section~\ref{coordsSect}, to provide a
fairly unified treatment of $\SO(2,n)$ and $\SU(2,n)$.

 \begin{thm} \label{SUF->complete}
 Assume $G = \SU(2,2m+1; \F)$ with $m \ge 1$,
 and let $H$ be any closed, connected subgroup of $G$, such that neither $H$
nor~$G/H$ is compact.

 If Conjecture~\ref{notessSU/Sp} is true, then $G/H$ does not have a
tessellation.
 \end{thm}

\begin{proof}
 Assume Conjecture~\ref{notessSU/Sp} is true, and suppose $\Gamma$ is a
crystallographic group for $G/H$. (This will lead to a contradiction.)
 Let 
 $$ \mbox{$H_1 = \SU(1,2m+1;\F)$ and $H_2 = \Symp(1,m;\F)$} $$
 \see{SU1nDefn}. From~\pref{d(Sp)}, we have $d(H_1) = \df(2m+1)$ and $d(H_2)
= \df(2m)$, where $\df = \dimR \F$.
 We may assume that $H \subset AN$ \see{HcanbeAN}, and that $H$ is compatible
with~$A$ \see{conjtocompatible}.

Because $H$ is not a Cartan-decomposition subgroup \see{CDS->notess}, the
contrapositive of Proposition~\ref{CDS<>h_m} implies, for some $k \in
\{1,2\}$, that there does \textbf{not} exist a continuous curve $h^t \to
\infty$ in~$H$, such that $\rho(h^t) \asymp \lVert h^t\rVert^k$. Therefore,
either Theorem~\ref{maxnolinear} (if $k = 1$) or Theorem~\ref{bestnosquare}
(if $k = 2$) implies that $d(H) \le \df(2m+1) = d(H_1)$.

We consider two cases.

\setcounter{case}{0}

\begin{case}
 Assume that $\Gamma$ acts properly discontinuously on $G/H_1$.
 \end{case}
  Theorem~\fullref{noncpctdim}{tess} (combined with the fact that $d(H) \le
d(H_1)$) implies that $G/H_1$ has a tessellation. This contradicts either
Theorem~\ref{SO2n/SO1odd-notess} (if $\F = \real$) or
Conjecture~\ref{notessSU/Sp}\ref{notessSU/Sp-SU2/SU1} (if $\F = \complex$).

\begin{case}
 Assume that $\Gamma$ does not act properly discontinuously on $G/H_1$.
 \end{case}
  From Lemma~\ref{mu(SUorSp)}, we know that $\mu(H_1)$ and $\mu(H_2)$ are
the two walls of~$A^+$, so Corollary~\ref{tess->missHk} (combined with
the assumption of this case) implies that $\Gamma$ acts properly
discontinuously on $G/H_2$. Therefore, since
Conjecture~\ref{notessSU/Sp}\ref{notessSU/Sp-SO2/SU1}\ref{notessSU/Sp-SU2/Sp1}
asserts that $G/H_2$ does not have a tessellation, the contrapositive of
Theorem~\fullref{noncpctdim}{tess} (with $H_2$ in the role of~$H_1$) implies
that $d(H) > d(H_2) = \df(2m)$. Hence, the contrapositive of
Theorem~\ref{maxnolinear} implies there is a continuous curve $h^t \to
\infty$ in~$H$, such that $\rho(h^t) \asymp h^t$. Thus, there is a compact
subset~$C$ of~$G$, such that $H_1 \subset CHC$ \see{SU1inH}. Since $\Gamma$
acts properly discontinuously on $G/H$, this implies that $\Gamma$ acts
properly discontinuously on $G/H_1$ \see{CHCproper}. This contradicts the
assumption of this case.
 \end{proof}

\begin{thm} \label{SUFevenTess}
 Assume $G = \SU(2,2m;\F)$ with $m \ge 2$, and let $H$ be a closed, connected,
nontrivial, proper subgroup of $AN$.

 The homogeneous space $G/H$ has a tessellation if and only if either
 \begin{enumerate}
 \item \label{SUFevenTess-Sp}
 there is an $\real$-linear map $B \colon \F^{n-2} \to \F^{n-2}$, such that
 \begin{enumerate}
 \item \label{SUFevenTess-Sp-xB}
 $\Im \bigl( (vB)(wB)^\dagger \bigr) = -\Im (v w^\dagger)$ for every $v,w \in
\F^{n-2}$ \see{Bsymplectic}, and
 \item \label{SUFevenTess-Sp-<x,y>}
 $xB \notin \F x$, for every nonzero $x \in \F^{n-2}$ \see{xBnotinFx}, and
 \item \label{SUFevenTess-Sp-HB}
 $H$ is conjugate to~$H_B$ {\upshape(}see~\ref{HB-defn}
and~\ref{HB=Sp1m}{\upshape)}; or
 \end{enumerate}
 \item \label{SUFevenTess-SUR}
 $\F = \real$ and $H$ is conjugate to $\SU(1,2m;\real) \cap AN$
\see{SU1nDefn}; or
 \item \label{SUFevenTess-SUC}
 $\F = \complex$ and there exists $c \in (0,1]$, such that $H$ is conjugate
to $\Hc$ \see{SUegsDefn}.
 \end{enumerate}
 \end{thm}

\begin{proof}
 ($\Leftarrow$) See
\pref{SUFevenTess-Sp}~Theorem~\fullref{HBthm}{tess},
\pref{SUFevenTess-SUR}~Theorem~\ref{SUevenTessExists} (and~\ref{d(Sp)}), or
\pref{SUFevenTess-SUC}~Theorem~\fullref{SUegs}{tess}.

 ($\Rightarrow$)
 Let $n = 2m$, so $G = \SU(2,n;\F)$.
 By combining Remark~\ref{d(H)=dimH}, Corollary~\ref{tess->dim>1,2},
Lemma~\ref{mu(SUorSp)}, and Remark~\ref{d(Sp)}, we see that 
 $$\dim H = d(H)
 \ge \min \bigl\{\, d \bigl( \SU(1,n; \F) \bigr), d \bigl( \Symp(1,m; \F) \bigr)
\, \bigr\}
 =\df n .$$

Also, we may assume $H$ is compatible with~$A$ \see{conjtocompatible}.
 Because $H$ is not a Cartan-decomposition subgroup \see{CDS->notess},
Proposition~\ref{CDS<>h_m} implies that one of the following two cases
applies.

\setcounter{case}{0}

\begin{case} \label{SUFevenTessPf-nolinear}
 Assume there does not exist a continuous curve $h^t \to \infty$ in~$H$, such
that $\rho(h^t) \asymp h^t$.
 \end{case}
 Since $\dim H \ge \df n$, Theorem~\ref{maxnolinear} implies that $\dim H
= \df n$, and that $H$ is of the form $H=T\ltimes U$ (with $U \subset N$),
where
 \begin{enumerate}
 \renewcommand{\theenumi}{\roman{enumi}}
 \item \label{SUFevenTessPf-TU}
 $T = \ker \alpha$;
 \item \label{SUFevenTessPf-phi0}
 $\phi_u = 0$ for every $u \in \Lie U$;
 \item \label{SUFevenTessPf-<xy>}
 $\dimF ( \F x_u + \F y_u ) = 2$, for every $u \in \Lie U \smallsetminus \Dh$;
 \item \label{SUFevenTessPf-Z}
 $|\eta_z|^2 + \xx_z \yy_z \neq 0$ for every
nonzero $z \in \Dh$;
 \item \label{SUFevenTessPf-dimU/Z}
 $\dim \Lie U / \Dh = \df(n-2)$; and
 \item \label{SUFevenTessPf-dimZ}
 $\dim \Dh = 2\df - 1$.
 \end{enumerate}

\setcounter{step}{0}
 \renewcommand{\thestep}{\thecase.\arabic{step}}

\begin{step} \label{SUFevenTessSpPf-xx=yy}
 We may assume that
 $\Dh = \{\, z \in \Dn \mid \xx_z = - \yy_z \,\}$.
 \end{step}
 Because $\dim \Dh = 2\df - 1$ \see{SUFevenTessPf-dimZ}, it
suffices to show that $\xx_z = - \yy_z$ for all $z \in \Dh$. This is
trivially true if $\F=\real$, as $\xx_z,\yy_z \in \Fim =\{0\}$ in this case.
Thus, we assume $\F = \complex$.

For any $z \in \Dh$ with $\eta_z = \yy_z = 0$, we know,
from~\pref{SUFevenTessPf-Z}, that $z = 0$; therefore,
Lemma~\fullref{O(linear)}{O} implies there exist $\real$-linear maps $R
\colon \complex \to i\real$ and $S \colon i \real \to i \real$, such that
$\xx_z = R(\eta_z) + S(\yy_z)$ for all $z \in \Dh$. More concretely, we
may say that there exist $\lambda \in \complex$ and $c \in \real$, such that
$\xx_z = \Im(\lambda \eta_z) + c \yy_z$ for all $z \in \Dh$.

Let $v$ be the element of~$\rsp_\alpha$ with $\phi_v = \cjg{\lambda}/2$, and
let $H^* = \exp(-v) H \exp(v)$ be the conjugate of~$H$ by $\exp(v)$. Then
$H^*$ satisfies the conditions imposed on~$H$ (note that $H^*$, like~$H$, is
compatible with~$A$ \see{conjUomega}), so there exist $\lambda^* \in
\complex$ and $c^* \in \real$, such that $\xx_{z^*} = \Im(\lambda^*
\eta_{z^*}) + c^* \yy_{z^*}$ for all $z^* \in \Dh^*$. Given $z^* \in \Dh^*$ with $\yy_{z^*} = 0$, let $z = \exp(v) z \exp(-v)$. Because $\yy_{z^*} =
0$, we have $\bigl[ [z^*,-v], -v \bigr] = 0$, so, from
Remark~\ref{conjugation} and \eqref{[u,v]}, we see that
 $$\yy_z = \yy_{z^*} = 0 ,$$
 $$\eta_z = \eta_{z^*} - \phi_{-v} \yy_{z^*} = \eta_{z^*} ,$$
 and
 $$ \xx_z
 = \xx_{z^*} + 2 \Im(\phi_{-v} \cjg{\eta_{z^*}})
 = \xx_{z^*} + 2\Im \bigl( (-\cjg{\lambda}/2) \cjg{\eta_z} \bigr)
 = \xx_{z^*} + \Im( \lambda \eta_z )
 = \xx_{z^*} + \xx_z .$$
 Therefore 
 $$ 0 = \xx_{z^*} = \Im(\lambda^* \eta_{z^*}) + c^* \yy_{z^*} = \Im(\lambda^*
\eta_{z^*}) .$$
 Since $\eta_{z^*}$ is arbitrary, this implies $\lambda^* = 0$. Thus, by
replacing $H$ with~$H^*$, we may assume that $\lambda = 0$. This means that
$\yy_z = c \xx_z$ for all $z \in \Dh$.

From~\pref{SUFevenTessPf-dimZ} (and because $\F = \complex$, so $\df = 2$), we
know that $\dim \Dh = 3 > 1$, so there is some nonzero $w \in \Dh$, such
that $\yy_w = 0$. (So $\xx_w = c \yy_w = 0$.) Then 
 $|\eta_w|^2 + \xx_w \yy_w = |\eta_w|^2 > 0$, so we see,
from~\pref{SUFevenTessPf-Z}, that $|\eta_z|^2 + \xx_z \yy_z  > 0$ for every
nonzero $z \in \Dh$. Now, since 
 $$\dim \Dh = 3 > 2 = \dim \rsp_{\alpha+2\beta} ,$$
 there is some nonzero $z \in \Dh$, such that $\eta_z = 0$. We have
 $$ 0 < |\eta_z|^2 + \xx_z \yy_z = 0 + c \yy_z^2 .$$
 Because $\yy_z$ is pure imaginary, we know that $\yy_z^2 < 0$, so this
implies that $c < 0$. Thus, replacing $H$ by a conjugate under a diagonal
matrix, we may assume $c = -1$, as desired.

\begin{step} \label{SUFevenTessPf-U'+Z}
 Setting $\Lie U' = (\rsp_\alpha + \rsp_{\alpha+\beta}) \cap \Lie U$, we have
$\Lie U = \Lie U' + \Dh$.
 \end{step}
 Since $T = \ker \alpha$ \see{SUFevenTessPf-TU}, we have
 \begin{itemize}
 \item $\beta|_T = (\alpha+\beta)|_T$,
 \item $2\beta|_T = (\alpha+2\beta)|_T = (2\alpha+2\beta)|_T$, and 
 \item $\beta|_T \neq 2\beta|_T$.
 \end{itemize}
 Thus, in the notation of Lemma~\ref{rootdecomp}, we have $\Lie N^{=\beta}
\cap \Lie U = \Lie U'$ and $\Lie N^{\neq \beta} \cap \Lie U = \Dh$, so
$\Lie U = \Lie U' \oplus \Dh$, as desired. (Note that this is a direct sum
of vector spaces, not of Lie algebras: we have $[\Lie U',\Lie U'] \subset
\Dh$.)

\begin{step}
 Completion of the proof of Case~\ref{SUFevenTessPf-nolinear}.
 \end{step}
 For any $u \in \Lie U'$ with $x_u = 0$, we have
 $$\dimF (\F x_u + \F y_u) = \dimF \F y_u \le 1  < 2 ,$$
 so $u \in \Lie U' \cap \Dh = \{0\}$ \see{SUFevenTessPf-<xy>}; therefore,
Lemma~\fullref{O(linear)}{O} implies there is a $\real$-linear map $B \colon
\F^{n-2} \to \F^{n-2}$, such that $y_u = x_u B$ for all $u \in \Lie U'$.
Then, because
 $$\dim \Lie u_{\alpha+\beta} = \dimR \F^{n-2} = \df(n-2) = \dim \Lie U'$$
 \see{SUFevenTessPf-dimU/Z}, we must have
 $$ \Lie U' = \{\, u \in \rsp_\beta + \rsp_{\alpha+\beta} \mid y_u = x_u B
\,\} .$$
 Combining this with \pref{SUFevenTessPf-TU} and the conclusions of
Steps~\ref{SUFevenTessSpPf-xx=yy} and~\ref{SUFevenTessPf-U'+Z}, we see that
$\Lie H = \Lie H_B$. Therefore $H = H_B$, so
Conclusion~\pref{SUFevenTess-Sp-HB} holds.

From~\pref{SUFevenTessPf-<xy>}, we see that
Conclusion~\pref{SUFevenTess-Sp-<x,y>} holds.

Letting $z = [u,v]$, for any $u,v \in \Lie U'$, we see, from \eqref{[u,v]},
that
 $$ \xx_z = -2 \Im (x_u x_v^\dagger)$$
 and
 $$ \yy_z = -2 \Im (y_u y_v^\dagger)
 = -2 \Im \bigl( (x_u B)(x_v B)^\dagger \bigr) .$$
 From Step~\ref{SUFevenTessSpPf-xx=yy}, we know that $\yy_z = - \xx_z$, so
this implies that Conclusion~\pref{SUFevenTess-Sp-xB} holds.

\begin{case} \label{SUFevenTessPf-nosquare}
 Assume there does not exist a continuous curve $h^t \to \infty$ in~$H$, such
that $\rho(h^t) \asymp \lVert h^t\rVert^2$.
 \end{case}
 Since $\dim H \ge \df n$, Theorem~\ref{bestnosquare} implies that $\dim H
= \df n$, and that $H$ is of the form $H=T\ltimes U$, where
 \begin{enumerate}
 \renewcommand{\theenumi}{\roman{enumi}}
 \item $T=\ker \beta$,
 \item $\Lie U=((\rsp_\alpha+\rsp_{\alpha+\beta}+
  \rsp_{\alpha+2\beta})\cap\Lie U)+\rsp_{2\alpha+2\beta}$,
  and
 \item\label{pdqf}
  $|x_u|^2+2\Re(\phi_u \cjg{\eta_u})\neq0$
  for every $u\in\Lie U\smallsetminus\rsp_{2\alpha+2\beta}$.
 \end{enumerate}

Let 
 $$\Lie U' = (\rsp_\alpha+\rsp_{\alpha+\beta}+ \rsp_{\alpha+2\beta})\cap\Lie
U$$
 (so $\Lie U = \Lie U' \oplus \rsp_{2\alpha+2\beta}$).  
 Let $Q$ be the sesquilinear form (or bilinear form, if $\F = \real$) on
 $\F \oplus \F^{n-2}\oplus\F$ defined by
 $$ Q \bigl( (\phi_1, x_1, \eta_1) , (\phi_2, x_2, \eta_2) \bigr)
 = \phi_1 \cjg{\eta_2} + x_1 x_2^\dagger + \eta_1 \cjg{\phi_2} .$$
 Let
 $$ V_{\Lie H} = \{\, (\phi_u, x_u, \eta_u) \in \F \oplus \F^{n-2} \oplus \F
\mid
 u \in \Lie U' \,\}.$$
 From~\pref{pdqf}, we see that the restriction of $\Re Q$ to~$V_{\Lie H}$ is a
(positive-definite) inner product.

Let $V_{\Lie H}^\perp$ be the $(\Re Q)$-orthogonal complement to~$V_{\Lie H}$. As a form
over~$\F$, $Q$ has signature $(1,n-1)$. Thus, as a form over~$\real$, $\Re Q$
has signature $\bigl( \df, \df(n-1) \bigr)$. Since 
 $$\dim V_{\Lie H}
 = \dim \Lie H - \dim \Lie T - \dim \rsp_{2\alpha+2\beta}
 = \df n - 1 - (\df-1)
 = \df(n-1) ,$$
 we conclude that $V_{\Lie H}^\perp$ is a $\df$-dimensional $\real$-subspace
on which $\Re Q$ is negative-definite.

 Choose some nonzero $u \in V_{\Lie H}^\perp$. Multiplying by a real scalar
to normalize, we may assume $Q(u,u) = -2$. Because $\SU(1,n-1)$ is
transitive on the vectors of norm~$-1$, there is some $g \in \SU(\Re Q)$,
such that $g(u) = (1,0,-1)$. Thus, letting
 \begin{equation} \label{ghat}
 \hat g = 
  \begin{pmatrix}
  1&0&0\\
  0&g&0\\
  0&0&1
  \end{pmatrix}
  \in\SU(2,n;\F) ,
 \end{equation}
 and $\Lie H^\sharp = {\hat g}^{-1} \Lie H g$, we have $(1,0,-1) \in V_{\Lie
H^\sharp}$, so, by replacing $\Lie H$ with the conjugate $\Lie H^\sharp$, we
may assume $u = (1,0,-1)$.

Then
 \begin{equation} \label{uperp}
 \begin{split}
 V_{\Lie H} 
 &\subset u^\perp \\
 &= (1,0,-1)^\perp \\
 &= \{\, (\phi, x, \eta) \mid
 \Re Q \bigl( (\phi, x, \eta), (1,0,-1) \bigr) = 0 \,\} \\
 &= \{\, (\phi, x, \eta) \mid
 \Re \bigl( \phi (-1) + x (0^\dagger) + \eta (1) \bigr) = 0 \,\} \\
 &= \{\, (\phi, x, \eta) \mid \Re \eta = \Re \phi \,\}
 .
 \end{split}
 \end{equation}

\begin{subcase}
 Assume $\F = \real$. 
 \end{subcase}
 By comparing \pref{uperp} and \pref{SU1nAN} (with $\F = \real$), we conclude
that
 $$\Lie H \subset \SU(1,n;\real) \cap (\Lie A + \Lie N) .$$
 By comparing dimensions, we see that equality must hold; this establishes
Conclusion~\pref{SUFevenTess-SUR}.

\begin{subcase}
 Assume $\F = \complex$.
 \end{subcase}
 Choose some nonzero~$v \in V_{\Lie H}^\perp$, such that $v$ is $(\Re
Q)$-orthogonal to~$u$.  Multiplying by a real scalar
to normalize, we may assume $Q(v,v) = -2$. By replacing $v$ with~$-v$ if
necessary, we may assume $\bigl( \Im Q(u,v) \bigr)/i \ge 0$.

Because $v$ is $(\Re Q)$-orthogonal to $u = (1,0,-1)$, we have $\Re \eta_v =
\Re \phi_v$ \see{uperp}. Let $s = (\Im \phi_v)/i$ and $t = (\Im \eta_v)/i$.
Then
 \begin{align*}
  0
 &\le \bigl( \Im Q(u,v) \bigr)/i \\
 &= \Bigl( \Im \bigl( (1) \cjg{\eta_v}  + 0 (x_v^\dagger) + (-1) \cjg{\phi_v}
\bigr) \Bigr)/i \\
 &= \bigl( - \Im \eta_v +  \Im \phi_v \bigr)/i \\
 &= -t + s,
 \end{align*}
 so 
 $$ \bigl( \Im Q(u,v) \bigr)/i = |s-t| .$$
 Also,
 \begin{align*}
 -2
 &= Q(v,v) \\
 & = |x_v|^2 + 2 \Re(\phi_v \cjg{\eta_v}) \\
 &= |x_v|^2 + 2 (\Re \phi_v)^2 - 2 (\Im \phi_v) (\Im \eta_v) \\
 &\ge - 2 (\Im \phi_v) (\Im \eta_v) \\
 &= 2st
 ,
 \end{align*}
 so $st \le -1$. Thus, $s$ and~$t$ are of opposite signs so, because $|s| \,
|t| \ge 1$, we have 
 $$ \bigl( \Im Q(u,v) \bigr) /i = |s-t| = |s| + |t| \ge 2 .$$
 Therefore, we may choose $c \in
(0,1]$, such that 
 $$ \Im Q(u,v) = i \left( c + \frac{1}{c} \right) .$$
 Let 
 $$ w = \left( \frac{i}{c}, 0, - ic \right) .$$
 Then 
 $$ Q(w,w) 
 = |x_w|^2 + 2 \Re( \phi_w \cjg{\eta_w})
 = 0^2 + 2 \bigl( i / c \bigr) \bigl( i c \bigr)
 = -2
 = Q(v,v) ,$$
 and
 $$ Q(u,w)
 = \phi_u \cjg{\eta_w} + x_u x_w^\dagger + \eta_u \cjg{\phi_w}
 = (1)(i c) + 0 + (-1)(-i / c )
 = i \left( c + \frac{1}{c} \right)
 = \Im Q(u,v) . $$
 Hence, there is some $h \in \SU(Q)$, such that $h(u) = u$ and $h(v) = w$.
Thus, replacing $\Lie H$ with the conjugate ${\hat h}^{-1} \Lie H \hat h$
\cf{ghat}, we may assume $v = w$.

Therefore 
 \begin{align*}
 V_{\Lie H} 
 &\subset v^\perp \\
 &= (cv)^\perp \\
 &= (i,0,-ic^2)^\perp \\
 &= \{\, (\phi, x, \eta) \mid
 \Re Q \bigl( (\phi, x, \eta), (i,0,-ic^2) \bigr) = 0 \,\} \\
 &= \{\, (\phi, x, \eta) \mid
 \Re \bigl( \phi (ic^2) + x (0^\dagger) + \eta (-i) \bigr) = 0 \,\} \\
 &= \{\, (\phi, x, \eta) \mid \Im \eta = c^2 \Im \phi \,\}
 .
 \end{align*}
 By combining this with \eqref{uperp} and comparing with \pref{SUegsDefn}
(with $\F = \complex$), we conclude that
 $\Lie H \subset \LieHcc{c^2} $.
 By comparing dimensions, we see that equality must hold; this establishes
Conclusion~\pref{SUFevenTess-SUC} (because $0 < c^2 \le 1$).
 \end{proof}

Theorem~\ref{SUevenTess} can be restated in the following more elementary (but
less precise) form.

\begin{cor} \label{SUF-known}
 Let $H$ be a closed, connected subgroup of $G = \SU(2,2m; \F)$ with $m \ge
2$, such that neither $H$ nor $G/H$ is compact, and let $\df = \dimR \F$.

 The homogeneous space $G/H$ has a tessellation if and only if
 \begin{enumerate}
 \item \label{SUF-known-d(H)}
 $d(H) = 2\df m$; and
 \item \label{SUF-known-sim}
 either $H \sim \SU(1,2m;\F)$ or $H \sim \Symp(1,m;\F)$.
 \end{enumerate}
 \end{cor}

\begin{proof}
 ($\Leftarrow$) This is Theorem~\ref{SUevenTessExists}.

 ($\Rightarrow$) Theorem~\ref{SUFevenTess}($\Rightarrow$) provides us with
three possibilities. 

\pref{SUF-known-d(H)} In each case, we have $d(H) = 2 \df m$ (see
\ref{d(HB)}, \ref{d(Sp)}, and~\ref{d(Hc)}). 

\pref{SUF-known-sim} In each case, there is some $k \in \{1,2\}$, such that
$\rho(h) \asymp \lVert h\rVert^k$ for $h \in H$ (see \fullref{HBthm}{mu},
\fullref{mu(SUorSp)}{SU}, and~\fullref{SUegs}{linear}). Then
Corollary~\ref{HsimSUorSp} implies either that $H \sim \SU(1,2m;\F)$ (if $k =
1$) or that $H \sim \Symp(1,m;\F)$ (if $k = 2$).
 \end{proof}

The following proposition shows that no further restriction can be placed
on~$c$ in the statement of Theorem~\fullref{SUFevenTess}{SUC}.

\begin{prop} \label{HcUncountable}
 If $\F = \complex$, then $\Hc$ is not conjugate to $\Hcc{c'}$, unless $c =
c'$ {\upshape(}for $c,c' \in (0,1]${\upshape)}. 
 \end{prop}

\begin{proof}
 Suppose $g^{-1} \Hc g = \Hcc{c'}$, for some $g\in G = \SU(2,2m)$. Because all
maximal split tori in $\Hcc{c'}$ are conjugate, we may assume that $g$
normalizes $\ker\beta$. Since all roots of $\ker\beta$ on both $\LieHc$
and~$\LieHcc{c'}$ are positive, $g$~cannot invert $\ker \beta$, so we conclude
that $g$ centralizes $\ker \beta$; that is, $g \in C_G(\ker \beta)$.

In the notation of Case~\ref{SUFevenTessPf-nosquare} of the proof of
Theorem~\ref{SUFevenTess}, define
 $$ S = \{\, \hat h \mid h \in \SU(Q) \,\}$$
 \cf{ghat}. Then
 $C_G( \ker \beta ) = (\ker \beta) S$, so we may assume $g \in S$ (because
$\ker \beta$, being a subgroup of~$\Hc$, obviously normalizes~$\Hc$). Write
$g = \hat h$. Then, because $g^{-1} \Hc g = \Hcc{c'}$, we must have
$h(V_{\LieHc}) = V_{\LieHcc{c'}}$; hence $h(V_{\LieHc}^\perp) =
V_{\LieHcc{c'}}^\perp$. 

For any basis $\{u,v\}$ of $V_{\LieHc}^\perp$ with $Q(u,u) = Q(v,v) = -2$ and
$\Re Q(u,v) = 0$, we have 
 $$\Im Q(u,v) = \pm i \bigl( c + (1/c) \bigr) .$$
 Similarly for any $(\Re Q)$-orthonormal basis $\{u',v'\}$ of
$V_{\LieHcc{c'}}^\perp$. Because $h \in \SU(Q)$, this implies
 $c + (1/c) = c' + (1/c')$. Because $c,c' \in (0,1]$, we conclude that $c =
c'$, as desired.
 \end{proof}


\begin{thebibliography}{GWn}

\bibitem[Abe]{Abels}
 H.~Abels,
 Properly discontinuous groups of affine transformations,
 (preprint).
 \bibURL
 {http://www.mathematik.uni-bielefeld.de/sfb343/preprints/pr00025.pdf.gz}

\bibitem[AMS]{AbelsMargulisSoifer}
 H.~Abels, G.~A.~Margulis and G.~A.~Soifer,
 Properly discontinuous groups of affine transformations with orthogonal
linear part,
 \emph{Comptes Rendus Acad. Sci. Paris} 324~I (1997) 253--258.

\bibitem[Ben]{Benoist}
 Y.~Benoist, Actions propres sur les espaces homog\`enes r\'eductifs,
 \emph{Ann. Math.} 144 (1996) 315--347.

\bibitem[BL]{BenoistLabourie}
 Y.~Benoist and F.~Labourie,
 Sur les espaces homog\`enes mod\`eles de vari\'et\'es compactes,
 \emph{Publ. Math. I.H.E.S} 76 (1992) 99--109.

\bibitem[Br1]{Borel-CK}
 A.~Borel,
 Compact Clifford-Klein forms of symmetric spaces,
 \emph{Topology} 2 (1963), 111--122.

\bibitem[Br2]{Borel-Algic}
 A.~Borel,
 \emph{Linear Algebraic Groups, 2nd~ed,}
 Springer, New York, 1991.

\bibitem[BT]{BorelTits-Reductive}
 A.~Borel and J.~Tits,
 Groupes r\'eductifs,
 \emph{Publ. Math. IHES} 27 (1965) 55--150.

\bibitem[Cha]{CharlapBook}
 L.~S.~Charlap,
 \emph{Bieberbach Groups and Flat Manifolds,}
 Springer, New York, 1986.

\bibitem[Coh]{Cohen-coho1}
 D.~E.~Cohen,
 \emph{Groups of Cohomological Dimension One,}
 Lecture Notes in Math. \#245,
 Springer, New York, 1972.

\bibitem[Cow]{Cowling}
 M.~Cowling, Sur les coefficients des repr\'esentations unitaires des
groupes de Lie simple, in:
 P.~Eymard, J.~Faraut, G.~Schiffmann, and R.~Takahashi, eds.: 
 \emph{Analyse Harmonique sur les Groupes de Lie~II}
 (S\'eminaire Nancy-Strasbourg 1976--78), Lecture Notes in Math. \#739,
Springer, New York, 1979, pp.~132--178.

\bibitem[Dix]{Dixmier-exp}
 J.~Dixmier,
 L'application exponentielle dans les groupes de Lie r\'esolubles,
 \emph{Bull. Soc. Math. France} 85 (1957) 113--121.

\bibitem[Dol]{Dold}
 A.~Dold,
 \emph{Lectures on Algebraic Topology, 2nd ed.,}
 Springer, New York, 1980.

 \bibitem[EN]{EllisNerurkar}
 R.~Ellis and M.~Nerurkar,
 Enveloping semigroup in ergodic theory and a proof of
Moore's ergodicity theorem,
 in: J.~C.~Alexander., ed.,
 \emph{Dynamical systems (College
Park, MD, 1986--87)},  Lecture Notes in Math.
\#1342,  Springer, New York, 1988, pp.~172--179.

\bibitem[FG]{FriedGoldman}
 D.~Fried and W.~M.~Goldman,
 Three-dimensional affine crystallographic groups,
 \emph{Adv. Math.} 47 (1983) 1--49.

\bibitem[Gol]{Goldman-nonstandard}
 W.~Goldman,
 Nonstandard Lorentz space forms,
 \emph{J. Diff. Geom.} 21 (1985) 301--308.

\bibitem[GdW]{GoodmanWallach}
 R.~Goodman and N.~Wallach,
 \emph{Representations and Invariants of the Classical Groups,}
 Cambridge U. Press, Cambridge, 1998.

\bibitem[GtW]{GotoWang}
 M.~Goto and H.--C.~Wang,
 Non-discrete uniform subgroups of semisimple Lie groups,
 \emph{Math. Ann.} 198 (1972) 259--286.

\bibitem[Gro]{Gromov-asymptotic}
 M.~Gromov, 
 Asymptotic invariants of infinite groups,
 in: G.~A.~Niblo and M.~A.~Roller, eds.,
 \emph{Geometric group theory, Vol.~2} (Sussex, 1991), 
 London Math. Soc. Lecture Notes \#182, 
 Cambridge Univ. Press, Cambridge, 1993, pp.~1--295.

\bibitem[Hel]{HelgasonBook}
 S.~Helgason,
 \emph{Differential Geometry, Lie Groups, and Symmetric Spaces,}
 Academic Press, New York, 1978.

\bibitem[Hc1]{Hochschild-Algic}
 G.~P.~Hochschild: 
 \emph{Basic Theory of Algebraic Groups and Lie Algebras,}
 Springer, New York, 1981.

\bibitem[Hc2]{Hochschild-Lie}
 G.~P.~Hochschild,
 \emph{The Structure of Lie Groups},
 Holden-Day, San Francisco, 1965.

\bibitem[How]{Howe}
 R.~Howe,
 A notion of rank for unitary representations of the classical groups, in:
 A.~Fig\`a Talamanca, ed.:
 \emph{Harmonic Analysis and Group Representations},
 (CIME 1980),
 Liguori, Naples, 1982, pp.~223--331.

\bibitem[HM]{HoweMoore}
 R.~Howe and C.~C.~Moore,
 Asymptotic properties of unitary representations,
 \emph{J. Func. Anal.} 32 (1979) 72--96.

\bibitem[Hm1]{Humphreys-Algic}
 J.~E.~Humphreys, \emph{Linear Algebraic Groups,}
 Springer, New York, 1975.

\bibitem[Hm2]{Humphreys-Lie}
 J.~E.~Humphreys,
 \emph{Introduction to Lie Algebras and Representation Theory,}
 Springer, New York, 1980.

\bibitem[Hus]{Husemoller}
 D.~Husemoller,
 \emph{Fibre Bundles, 2nd ed.,}
 Springer, New York, 1966.

\bibitem[IW]{IozziWitte-CDS}
 A.~Iozzi and D.~Witte,
 Cartan-decomposition subgroups of $\SU(2,n)$,
 \emph{J. Lie Theory} (to appear).
 \goodbreak
 {\tt http://front.math.ucdavis.edu/math.RT/0007039}

\bibitem[Iwa]{Iwasawa}
 K.~Iwasawa, On some types of topological groups,
 \emph{Ann. Math.} 50 (1949) 507--558.

\bibitem[Jac]{JacobsonBook}
 N.~Jacobson,
 \emph{Lie Algebras,}
 Dover, New York, 1979.

\bibitem[KS]{KatokSpatzier}
 A.~Katok and R.~Spatzier,
 First cohomology of Anosov actions of higher rank abelian groups and
applications to rigidity,
 \emph{IHES Publ. Math.} 79 (1994) 131--156.

\bibitem[KN]{KobayashiNomizu1}
 S.~Kobayashi and K.~Nomizu,
 \emph{Foundations of Differential Geometry, vol.~1},
 Interscience, New York, 1963.

\bibitem[Kb1]{Kobayashi-properaction}
 T.~Kobayashi,
 Proper action on a homogeneous space of reductive type,
 \emph{Math. Ann.} 285 (1989), 249--263.

\bibitem[Kb2]{Kobayashi-necessary}
 T.~Kobayashi, A necessary condition for the existence of compact
Clifford-Klein forms of homogeneous spaces of reductive type,
 \emph{Duke Math. J.} 67 (1992) 653--664.

\bibitem[Kb3]{Kobayashi-isotropy}
 T.~Kobayashi, On discontinuous groups acting on homogeneous spaces with
non-compact isotropy groups,
 \emph{J. Geom. Physics} 12 (1993) 133--144.

\bibitem[Kb4]{Kobayashi-criterion}
 T.~Kobayashi, Criterion of proper actions on homogeneous spaces of
reductive groups,
 \emph{J. Lie Th.} 6 (1996) 147--163.

\bibitem[Kb5]{Kobayashi-survey97}
 T.~Kobayashi,
 Discontinuous groups and Clifford-Klein forms of
pseudo-Riemannian homogeneous manifolds,
 in: B.~\O rsted and H.~Schlichtkrull, eds.,
 \emph{Algebraic and Analytic Methods in Representation Theory,}
 Academic Press, New York, 1997, pp.~99--165.

\bibitem[Kb6]{Kobayashi-deformation}
 T.~Kobayashi, Deformation of compact Clifford-Klein forms of
indefinite Riemannian homogeneous manifolds,
 \emph{Math. Ann.} 310 (1998) 395--409.

\bibitem[Kb7]{Kobayashi-survey00}
 T.~Kobayashi,
 Discontinuous groups for non-Riemannian homogeneous spaces,
 in: B.~Engquist and W.~Schmidt, eds.,
 \emph{Mathematics Unlimited---2001 and Beyond,}
 Springer, New York, 2001, pp.~723--747.

\bibitem[KO]{KobayashiOno}
 T.~Kobayashi and K.~Ono,
 Note on Hirzebruch's proportionality principle,
 \emph{J. Fac. Sci. Univ. Tokyo, Math.} 37 (1990) 71--87.

\bibitem[Kos]{Kostant}
 B.~Kostant, On convexity, the Weyl group, and the Iwasawa decomposition,
 \emph{Ann. Sc. ENS.} 6 (1973) 413--455.

\bibitem[Kul]{Kulkarni}
 R.~Kulkarni,
 Proper actions and pseudo-Riemannian space forms,
 \emph{Adv. Math.} 40 (1981) 10--51.

\bibitem[Lab]{Labourie-survey}
 F.~Labourie,
 Quelques r\'esultats r\'ecents sur les espaces localement
homog\`enes compacts,
 in: P.~de Bartolomeis, F.~Tricerri and E.~Vesentini, eds.,
 \emph{Manifolds and Geometry},
 Symposia Mathematica, v.~XXXVI,
 Cambridge U.~Press, 1996.

\bibitem[Mr1]{Margulis-Auslander}
 G.~A.~Margulis,
 Free properly discontinuous groups of affine transformations,
 \emph{Dokl. Akad. Nauk SSSR} 273 (1983) 937--940.

\bibitem[Mr2]{Margulis-CK}
 G.~A.~Margulis,
 Existence of compact quotients of homogeneous
spaces, measurably proper actions, and decay of matrix coefficients,
 \emph{Bull. Soc. Math. France} 125 (1997) 447--456.

\bibitem[MS]{MilnorStasheff}
 J.~W.~Milnor and J.~D.~Stasheff,
 \emph{Characteristic Classes,}
 Princeton U.\ Press, Princeton, 1974.

\bibitem[Mos]{Mostow-FSS}
 G.~D.~Mostow, Factor spaces of solvable groups,
 \emph{Ann. Math.} 60 (1954) 1--27.

\bibitem[Oh]{Oh-tempered}
 H.~Oh,
 Tempered subgroups and representations with minimal decay of matrix
coefficients,
 \emph{Bull. Soc. Math. France} 126 (1998) 355--380.

\bibitem[OW1]{OhWitte-CDS}
 H.~Oh and D.~Witte,
 Cartan-decomposition subgroups of $\SO(2,n)$,
 \emph{Trans. Amer. Math. Soc.}
 (to appear).
 \bibURL{http://front.math.ucdavis.edu/math.RT/9902049}

\bibitem[OW2]{OhWitte-announce}
 H.~Oh and D.~Witte,
 New examples of compact Clifford-Klein forms of homogeneous spaces of
$\SO(2,n)$,
 \emph{Internat. Math. Res. Not.}
 2000 (8 March 2000), no. 5, 235--251.

\bibitem[OW3]{OhWitte-CK}
 H.~Oh and D.~Witte,
 Compact Clifford-Klein forms of homogeneous spaces of $\SO(2,n)$,
 \emph{Geometriae Dedicata}
 (to appear).
 \bibURL{http://front.math.ucdavis.edu/math.RT/9902050}

\bibitem[Pal]{Palais-Slice}
 R.~S.~Palais,
 On the existence of slices for actions of non-compact Lie groups,
 \emph{Ann. Math.} 73, no.~2 (1961) 295--323.

\bibitem[PR]{PlatonovRapinchuk}
 V.~Platonov and A.~Rapinchuk,
 \emph{Algebraic Groups and Number Theory,}
 Academic Press, Boston, 1994.

\bibitem[Rag]{Raghunathan}
 M.~S.~Raghunathan,
 \emph{Discrete Subgroups of Lie Groups,}
 Springer, New York, 1972.

\bibitem[Sai]{Saito2}
 M.~Saito,
 Sur certains groupes de Lie r\'esolubles~II,
 \emph{Sci. Papers Coll. Gen. Ed. Univ. Tokyo} 7 (1957) 157--168.

\bibitem[Sal]{Salein}
 F.~Salein,
 Vari\'et\'es anti-deSitter de dimension 3 poss\'edant un champ de Killing
non trivial,
 \emph{Comptes Rendus Acad. Sci. Paris} 324~I (1997) 525--530.

\bibitem[Tom]{Tomanov}
 G.~Tomanov,
 The virtual solvability of the fundamental group of a generalized Lorentz
space form,
 \emph{J. Diff. Geom.} 32 (1990) 539--547.

\bibitem[Var]{Varadarajan}
 V.~S.~Varadarajan, 
 \emph{Lie Groups, Lie Algebras, and their Representations},
 Springer, New York, 1984. 

\bibitem[Whi]{WhiteheadBook}
 G.~W.~Whitehead,
 \emph{Elements of Homotopy Theory,}
 Springer, New York, 1978. 

\bibitem[Wit]{Witte-Solvtess}
 D.~Witte,
 Tessellations of solvmanifolds,
 \emph{Trans. Amer. Math. Soc.} 350 (1998), no. 9, 3767--3796.

\bibitem[ZS]{ZariskiSamuel1}
 O.~Zariski and P.~Samuel,
 \emph{Commutative Algebra, vol.~1,}
 Springer, New York, 1958.

\bibitem[Zeg]{Zeghib-deSitter}
 A.~Zeghib,
 On closed anti-deSitter spacetimes,
 \emph{Math. Ann.} 310 (1998), no.~4, 695--716.

\bibitem[Zm1]{Zimmer-orbitspace}
 R.~J.~Zimmer,
 Orbit spaces of unitary representations, ergodic theory, and simple Lie
groups,
 \emph{Ann. Math.} 106 (1977) 573--588.

\bibitem[Zm2]{ZimmerBook}
 R.~J.~Zimmer,
 \emph{Ergodic Theory and Semisimple Groups},
 Birkh\"auser, Boston, 1984.

\end{thebibliography}
\end{document}